\documentclass[a4paper,11pt]{article}
\usepackage{amsmath}
\usepackage{amssymb}
\usepackage{amsthm}
\usepackage{enumerate}
\usepackage{hyperref,url}
\usepackage{bm} 
\usepackage{color}
\usepackage[margin=30mm]{geometry}

\newcommand{\N}{\mathbb{N}}

\renewcommand{\P}{\mathbb{P}}
\newcommand{\R}{\mathbb{R}}
\newcommand{\dHaus}{\, d\mathcal{H}^{d-1}}
\newcommand{\dx}{\, dx}
\newcommand{\ds}{\, ds}
\newcommand{\dt}{\, dt}
\newcommand{\dd}{\, d}
\newcommand{\der}{\mathrm{D}}
\newcommand{\pd}{\partial}
\newcommand{\abs}[1]{\left| #1 \right|}
\newcommand{\norm}[1]{\| #1 \|}
\newcommand{\bignorm}[1]{\left\| #1 \right\|}
\newcommand{\inner}[2]{\langle #1 , #2 \rangle}

\newcommand{\eps}{\varepsilon}
\newcommand{\Laplace}{\Delta}

\newcommand{\md}{\pd^{\bullet}_{t}}

\renewcommand{\div}{\mathrm{div} \, }

\newcommand{\mean}[1]{\overline{ #1 }}
\newcommand{\Id}{\mathrm{I}}

\newcommand{\jump}[1]{\left [#1 \right ]_{1}^{2}}
\newcommand{\velo}{\mathcal{V}}
\newcommand{\pdnu}{\partial_{\bm{\nu}}}
\newcommand{\AInv}{\mathbb{A}}
\newcommand{\Stokes}{\mathbb{S}}
\theoremstyle{plain}
\newtheorem{thm}{Theorem}[]
\newtheorem{lemma}{Lemma}[section]

\newtheorem{remark}{Remark}[section]

\numberwithin{equation}{section}

\allowdisplaybreaks

\title{Thermodynamically Consistent Navier--Stokes--Cahn--Hilliard Models with Mass Transfer and Chemotaxis}

\author{Kei Fong Lam \footnotemark[1] \and Hao Wu \footnotemark[2]}

\date{ }

\begin{document}

\maketitle

\footnotetext[1]{Fakult\"at f\"ur Mathematik, Universit\"at Regensburg, 93040 Regensburg, Germany
({\tt Kei-Fong.Lam@mathematik.uni-regensburg.de}).}
\footnotetext[2]{School of Mathematical Sciences and Shanghai Key Laboratory for Contemporary Applied Mathematics, Fudan University, 220 Han Dan Road, Shanghai 20043, China ({\tt haowufd@yahoo.com}).}

\begin{abstract}
We derive a class of Navier--Stokes--Cahn--Hilliard systems that models two-phase flows with mass transfer coupled to the process of chemotaxis.
These thermodynamically consistent models can be seen as the natural Navier--Stokes analogues of earlier Cahn--Hilliard--Darcy models proposed for modeling tumor growth, and are derived based on a volume-averaged velocity. Then we perform mathematical analysis on the simplified model variant with zero excess of total mass and equal densities.
We establish the existence of global weak solutions in two and three dimensions for prescribed mass transfer terms.
Under additional assumptions, we prove the global strong well-posedness in two dimensions with variable fluid viscosity and mobilities, which also includes a continuous dependence on initial data and mass transfer terms for the chemical potential and the order parameter in strong norms.

\end{abstract}



\section{Introduction}
Biological phenomena such as tumor growth often involve complex interactions between various actors that take place over multiple spatial and temporal scales \cite{CLLW,CL,Foty,Frieboes,GLSS,Hawkins,OHP,Sitka}.  While discrete models that track the evolution of individual cells are able to capture the biophysical rules proposed from biological intuition, the rapid increase in computational costs with the number of cells and the difficulties encountered with model calibration place limitations on discrete models.  On the other hand, continuum models, which group multiple characteristic properties into one or more phenomenological parameters, tend to be too coarse at resolving the microstructures occurring at the level of individual cells, but in exchange they offer a tractable phenomenological description of the key dynamics for which mathematical analysis and numerical simulations can be carried out.

Within the class of continuum models, we focus on the category of diffuse interface models, also commonly known as phase field models \cite{Anderson}.  Although it is natural to view interfaces separating different components of matter as idealized hypersurfaces with zero thickness, these so-called sharp interface models break down when the interface experiences topological changes.  The diffuse interface models replace this hypersurface description of the interface with a thin layer where microscopic mixing of the macroscopically distinct components of matter are allowed.  This not only yields systems of equations that are better amenable to further analysis, but topological changes of the interface can also be handled naturally (see e.g., \cite{Anderson,LT}).

In continuum models, the ensemble of cells are often assumed to be tightly packed and move together with some averaged mixture velocity.  Grounded in the observation that the morphology of certain cells and tissues behave akin to viscous fluids \cite{Alber,Beysens,Foty,Ranft}, many models treat the cells as viscous incompressible fluids, leading to a mathematical consideration with the Navier--Stokes equations or Darcy's law for the macroscopic cellular velocity.  The latter has been a popular choice for tumor growth (see for example \cite{CL,Frieboes,GLSS,Wise} and the references therein), which treats the cells as viscous inertia-less fluids, similar in spirit to the approach taken in \cite{Ambrosi}.

In this work, we include inertia effects into the modeling, and obtain a Navier--Stokes system for the mixture velocity.  More precisely, we derive a thermodynamically consistent model for a two-component fluid mixture, which allows for mass transfer between the two components, in the presence of a chemical species that is subjected to diffusion and other transport mechanisms such convection and chemotaxis.  The general model is a Navier--Stokes--Cahn--Hilliard system of the following form
\begin{subequations}\label{Intro:NSCH}
\begin{alignat}{2}
\pd_{t} \rho + \div \left ( \rho \bm{v} - \frac{\overline{\rho}_{2} - \overline{\rho}_{1}}{2} m(\varphi) \nabla \mu \right )  & = \Gamma_{1} + \Gamma_{2}, \label{Intro:mass} \\
\div \bm{v}  & =  \frac{\Gamma_{1}}{\overline{\rho}_{1}} + \frac{\Gamma_{2}}{\overline{\rho}_{2}}, \label{Intro:div} \\
\pd_{t} (\rho \bm{v}) + \div \left (\rho \bm{v} \otimes \bm{v} - \frac{\overline{\rho}_{2} - \overline{\rho}_{1}}{2} m(\varphi) \bm{v} \otimes \nabla \mu \right )  & =  - \nabla p + \div \left ( 2 \eta(\varphi) \der \bm{v} \right ) \label{Intro:mom} \\
\notag &\quad \ - \div \left ( B \nabla \varphi \otimes \nabla \varphi \right ), \\
\pd_{t} \varphi + \div (\varphi \bm{v} - m(\varphi) \nabla \mu) & =  \frac{\Gamma_{2}}{\overline{\rho}_{2}} -  \frac{\Gamma_{1}}{\overline{\rho}_{1}}, \label{Intro:varphi} \\
\mu & = A \Psi'(\varphi) - B \Laplace \varphi + N_{\varphi}, \label{Intro:mu} \\
\pd_{t} \sigma + \div (\sigma \bm{v} - n(\varphi) \nabla N_{\sigma}) & = S. \label{Intro:sigma}
\end{alignat}
\end{subequations}
Here, $\overline{\rho}_{i}$, $i = 1,2$, is the constant mass density of a pure component of fluid $i$, $\bm{v}$ is the volume-averaged velocity of the mixture, $\der \bm{v} = \frac{1}{2}(\nabla \bm{v} + (\nabla \bm{v})^{\top})$ is the rate of deformation tensor, $p$ denotes the pressure, and $\rho, \eta$ denote the density and the viscosity of the mixture, respectively.  The order parameter $\varphi$ is the difference in volume fractions, with $\{\varphi = -1\}$ representing fluid 1 and $\{\varphi = 1\}$ representing fluid 2, and associated to $\varphi$ is the chemical potential $\mu$.  The function $\Psi'$ is the derivative of a potential $\Psi$ that has two equal minima at $\pm 1$.  The constants $A, B > 0$ are related to the surface tension and the thickness of the interfacial layers, while $m(\varphi)$ denotes a non-negative mobility.  The functions $\Gamma_{i}$, $i = 1,2$, models the mass transfer of fluid $i$.  The density of a chemical species is given by $\sigma$, and $S$ models a reaction term for $\sigma$.  The mobility of $\sigma$ is given by $n(\varphi) \geq 0$, and the functions $N_{\varphi}$ and $N_{\sigma}$ are the partial derivatives of the chemical free energy density $N(\varphi, \sigma)$ with respect to $\varphi$ and $\sigma$, respectively.  In the analysis below we will take $N(\varphi, \sigma)$ as
\begin{align}\label{Intro:Nutrient}
N(\varphi, \sigma) = \frac{1}{2} \abs{\sigma}^{2} + \chi \sigma(1-\varphi) \quad \Longrightarrow \quad N_{\varphi} = -\chi \sigma, \quad N_{\sigma} = \sigma + \chi(1-\varphi),
\end{align}
for a non-negative constant $\chi$ (see, e.g., \cite{GLSS}).

Equations \eqref{Intro:mass} and \eqref{Intro:mom} are the mass and momentum balance for the fluid mixture, while \eqref{Intro:div} relates the divergence of the volume-averaged velocity to the mass transfer terms.  Equations \eqref{Intro:varphi} and \eqref{Intro:mu} constitute a convective Cahn--Hilliard system for the order parameter $\varphi$, and equation \eqref{Intro:sigma} is a convection-reaction-diffusion equation for the chemical density $\sigma$.  Under suitable boundary conditions, the above model satisfies an energy identity of the form
\begin{align}\label{Intro:Energy:Id}
\notag & \frac{\dd}{\dt} \int_{\Omega} e \dx + \int_{\Omega} \left ( m(\varphi) \abs{\nabla \mu}^{2} + n(\varphi) \abs{\nabla N_{\sigma}}^{2} + 2 \eta(\varphi) \abs{\der \bm{v}}^{2} \right ) \dx \\
\notag & \quad = - \int_{\Omega} U_{v} \left ( \varphi \mu + \sigma N_{\sigma} + \frac{\overline{\rho}_{2} + \overline{\rho}_{1}}{2} \frac{\abs{\bm{v}}^{2}}{2} - p - A \Psi(\varphi) - \frac{B}{2} \abs{\nabla \varphi}^{2} - N(\varphi, \sigma) \right ) \dx \\
& \quad \quad + \int_{\Omega} \left[\left ( \mu - \frac{\overline{\rho}_{2} - \overline{\rho}_{1}}{2} \frac{\abs{\bm{v}}^{2}}{2} \right ) U_{\varphi} + N_{\sigma} S \right]\dx,
\end{align}
where the total energy density
$$e = \frac{\rho}{2} \abs{\bm{v}}^{2} + A \Psi(\varphi) + \frac{B}{2} \abs{\nabla \varphi}^{2} + N(\varphi, \sigma)$$ is the sum of the kinetic energy density $\frac{\rho}{2} \abs{\bm{v}}^{2}$, the chemical free energy density $N(\varphi, \sigma)$ and the Ginzburg--Landau energy density $A \Psi(\varphi) + \frac{B}{2} \abs{\nabla \varphi}^{2}$, and
\begin{align*}
U_{v} := \frac{\Gamma_{1}}{\overline{\rho}_{1}} + \frac{\Gamma_{2}}{\overline{\rho}_{2}}, \quad U_{\varphi} := \frac{\Gamma_{2}}{\overline{\rho}_{2}} - \frac{\Gamma_{1}}{\overline{\rho}_{1}}.
\end{align*}

The system \eqref{Intro:NSCH} is a generalization of the model derived by Abels, Garcke and Gr\"{u}n \cite{AGG} for two-phase flow with unmatched densities, where we account for mass transfer between the fluids and coupling to the dynamics of a chemical species.  Furthermore \eqref{Intro:NSCH} can also be seen as the Navier--Stokes analogue of the general Cahn--Hilliard--Darcy model derived in Garcke et al \cite{GLSS}.  In Section \ref{sec:Darcy} we briefly discuss how to obtain Darcy's law from the momentum equation by an averaging procedure (see \cite{DGL} for a simpler situation), which allows us to formally recover Cahn--Hilliard--Darcy models from \eqref{Intro:NSCH}.  To the authors' best knowledge, a similar Navier--Stokes--Cahn--Hilliard system was first derived by Sitka \cite{Sitka} for modeling tumor growth with chemotaxis, which can be obtained from the more general model \eqref{Intro:NSCH}--\eqref{Intro:Nutrient} with the following specific choice of mass transfer terms:
\begin{align}\label{zero:excess}
\Gamma_{2} = -\Gamma_{1},
\end{align}
so that
\begin{align}\label{zero:excess:source}
U_{v} = \alpha \Gamma_{2}, \ U_{\varphi} = \overline{\rho}_{S} \Gamma_{2}\ \text{ with } \alpha := \frac{1}{\overline{\rho}_{2}} - \frac{1}{\overline{\rho}_{1}}, \ \overline{\rho}_{S} := \frac{1}{\overline{\rho}_{2}} + \frac{1}{\overline{\rho}_{1}}.
\end{align}
The choice $\Gamma_{2} = - \Gamma_{1}$ stands for the case of zero excess of total mass, as any mass loss by fluid 1 is equal to the mass gained by fluid 2, and vice versa.  Furthermore, by integrating \eqref{Intro:mass} over the physical domain $\Omega$ and under suitable boundary conditions, the total mass $\int_{\Omega} \rho \dx$ is conserved.  For the modeling of tumor growth, Sitka \cite{Sitka}, and also \cite{CLLW,CL,Frieboes,GLSS,OHP,Wise} considered the following biologically relevant choices:
\begin{align*}
\Gamma_{2} = \frac{1}{2} \left ( \mathcal{G} \sigma - \mathcal{A} \right )(\varphi + 1), \quad S = -\frac{1}{2} \mathcal{C} \sigma (\varphi + 1),
\end{align*}
where the non-negative constants $\mathcal{G}$, $\mathcal{A}$, $\mathcal{C}$ correspond to the constant proliferation, apoptosis and nutrient consumption rates.  These mass transfer terms model the evolution of a tumor at its early stage, whose growth is proportional to the local chemical density, and experiences cell death at a constant rate.  The prefactor $(\varphi + 1)$ in $\Gamma_{2}$ and $S$ ensure that these source term only affects the tumor, which is given as the region $\{\varphi = 1\}$, and not the healthy tissues $\{\varphi = -1 \}$.  Next, let us mention that for the choice \eqref{Intro:Nutrient} of the chemical free energy density $N$, the fluxes for $\varphi$ and $\sigma$ are given by
\begin{align*}
\bm{q}_{\varphi} & := - m(\varphi) \nabla \mu = m(\varphi)\chi \nabla \sigma - m(\varphi) \nabla \left (A \Psi'(\varphi) - B \Laplace \varphi \right ),\\
\bm{q}_{\sigma} & := - n(\varphi) \nabla N_{\sigma} = n(\varphi) \chi \nabla \varphi - n(\varphi) \nabla \sigma.
\end{align*}
As pointed out in \cite{GLSS}, the parameter $\chi$ is related to transport mechanisms such as chemotaxis (movement of fluid 2 towards high regions of $\sigma$) and active transport (establishing a persistent concentration difference between the different fluid components even against the chemical concentration gradient).  The former is due to the fact that the term $m(\varphi) \chi \nabla \sigma$ in the flux $\bm{q}_{\varphi}$ points in the direction of increasing $\sigma$, and the latter is due to the fact that the term $n(\varphi) \chi \nabla \varphi$ in the flux $\bm{q}_{\sigma}$ points in the direction of increasing $\varphi$.  Note that $\nabla \varphi$ is only non-zero in the vicinity of the interface between the two fluids, as the order parameter takes distinct constant values in the fluid regions, and hence, we expect that the term $\div (n(\varphi) \chi \nabla \varphi)$ in \eqref{Intro:sigma}, with the choice \eqref{Intro:Nutrient} for $N$, attempts to drive the chemical species towards fluid 2, i.e., the region $\{\varphi = 1\}$, leading to a higher chemical density in fluid 2 than in fluid 1 near the interface.  We refer the reader to \cite[\S 5, Fig. 9 and 10]{GLSS} for numerical simulations investigating this mechanism, and to \cite{GLNeumann,GLSS} for a discussion on how to decouple the chemotaxis and active transport mechanisms as at present they are both related to the parameter $\chi$.

As a first step toward the mathematical analysis of system \eqref{Intro:NSCH}, we will study the special case of matched densities $\overline{\rho}_{1} = \overline{\rho}_{2} = 1$ and zero excess of total mass $\Gamma_{1} = -\Gamma_{2}$, for which the system \eqref{Intro:NSCH} reduces to
\begin{subequations}\label{Intro:ZEED}
\begin{align}
\pd_{t} \bm{v} + (\bm{v} \cdot \nabla ) \bm{v}  & = - \nabla p + \div (2 \eta(\varphi) \der \bm{v} - B \nabla \varphi \otimes \nabla \varphi), \\
\div \bm{v} & = 0, \\
\pd_{t} \varphi + \bm{v} \cdot \nabla \varphi & = \div (m(\varphi) \nabla \mu) + 2 \Gamma_{2}, \\
\mu & = A \Psi'(\varphi) - B \Laplace \varphi - \chi \sigma, \\
\pd_{t} \sigma + \bm{v} \cdot \nabla \sigma & = \div (n(\varphi) \nabla (\sigma - \chi \varphi)) + S.
\end{align}
\end{subequations}
 We remark that obtaining useful a priori estimates from the energy identity \eqref{Intro:Energy:Id} for the general model \eqref{Intro:NSCH} turns out to be rather complicated due to the source terms $U_{v} \abs{\bm{v}}^{2}$, $U_{\varphi} \abs{\bm{v}}^{2}$ involving the velocity and $U_{v} p$ involving the pressure.  Unlike in the analytical treatment of the Cahn--Hilliard--Darcy variant with source terms performed in \cite{GLDarcy,GLRome,JWZ}, useful estimates for the pressure $p$ involving the left-hand side of \eqref{Intro:Energy:Id} seem not available a priori in the case with the Navier--Stokes equations.  Hence, as a first contribution to the mathematical treatment of Navier--Stokes--Cahn--Hilliard systems with mass transfer and chemotaxis, we focus on the case of zero excess of total mass \eqref{zero:excess} and matched densities, so that $\alpha = 0$ and $U_{v} = 0$ in \eqref{zero:excess:source}.  We leave the analysis of the general model for future research.

The above approach to model chemotaxis bears both similarities and differences to the classical approach of Keller and Segel \cite{Keller}.  In the context of coupling Navier--Stokes flow with chemical chemotaxis, the mathematical model of by Tuval et al. \cite{Tuval} for oxygen-driven bacteria swimming in viscous incompressible fluids has been intensively studied by many authors, see for example \cite{Chae,Chae2,Francesco,Duan,JWZChemo,LiuLorz,Lorz,Winkler1,Winkler2,Winkler3,ZhangLi,ZhangZheng}.  For the oxygen density $c$, bacteria density $n$, fluid velocity $\bm{u}$ and pressure $p$, the model in \cite{Tuval} reads as
\begin{equation}\label{Intro:Tuval}
\begin{aligned}
\pd_{t} \bm{u} + (\bm{u} \cdot \nabla ) \bm{u} & = - \nabla p + \eta \Laplace \bm{u} + \bm{g},\\
\div \bm{u} & = 0, \\
\pd_{t} c + \bm{u} \cdot \nabla c & = D_{c} \Laplace c - n f(c), \\
\pd_{t} n + \bm{u} \cdot \nabla n & = D_{n} \Laplace n - \div (r(c) n \nabla c),
\end{aligned}
\end{equation}
where $D_{c}$, $D_{n}$ are the diffusivities of $c$ and $n$, respectively, $r(c)$ denotes the chemotactic sensitivity, $f(c)$ is a consumption rate of the chemical by the cells, and $\bm{g}$ accounts for external forces such as gravity or buoyancy forces.

Writing \eqref{Intro:varphi}, \eqref{Intro:mu} and \eqref{Intro:sigma} (with the specific choice \eqref{Intro:Nutrient}) as
\begin{equation}\label{Intro:sub:ZEED}
\begin{aligned}
\pd_{t}\varphi + \div (\varphi \bm{v}) & = \div (m(\varphi) \nabla (A \Psi'(\varphi) - B \Laplace \varphi)) - \div (m(\varphi) \chi \nabla \sigma) + U_{\varphi}, \\
\pd_{t}\sigma + \div (\sigma \bm{v}) & = \div (n(\varphi) \nabla \sigma) - \div (n(\varphi) \chi \nabla \varphi) + S,
\end{aligned}
\end{equation}
then, we have a correspondence between $\sigma$ and $c$ as the chemoattractant, between $\varphi$ and $n$ as the variable exhibiting the chemotactic movement, between $\bm{v}$ and $\bm{u}$ as the fluid velocity, and between the terms $\div (m(\varphi) \chi \nabla \sigma)$ and $\div (r(c) n \nabla c)$ as the mechanism for chemotaxis.  The main differences between \eqref{Intro:ZEED} (with the choice \eqref{Intro:Nutrient}) and \eqref{Intro:Tuval} are the appearance of a cross-diffusion type term $-\div (n(\varphi) \chi \nabla \varphi)$ in the equation of $\sigma$ accounting for the active transport mechanism, and the fact that $\varphi$ is subject to a fourth-order diffusion operator. Besides, \eqref{Intro:ZEED} is essentially a two-phase model.

Lastly, let us mention that the models derived here, and also in \cite{AGG,GLSS}, utilizes a volume-averaged velocity, which is in contrast to the models of \cite{CLLW,CL,Frieboes,OHP,Wise} that employ a mass-averaged velocity.  In the two-phase model of Lowengrub and Truskinovsky \cite{LT}, the mass-averaged velocity is not solenodial, whereas in the two-phase model of Abels, Garcke and Gr\"{u}n \cite{AGG}, the volume-averaged velocity is divergence-free, which seems to be an advantageous property in the analytical treatment as the pressure variable can be eliminated.  In the present setting with mass transfer between the two components in the Navier--Stokes flow, the velocity is not divergence-free in general.  Thus, it would be interesting to compare the analogous Navier--Stokes--Cahn--Hilliard models derived using a mass-averaged velocity with the models presented below.  We leave this comparison for future research.

The outline and contributions of the present work are as follows:
\begin{itemize}
\item In Section \ref{sec:volave}, we derive the general model \eqref{Intro:NSCH} from classical conservation laws.  Under specific choices of the source terms and free energy we obtain the Navier--Stokes analogue of the models in Garcke et al. \cite{GLSS}.  The sharp interface limit $\eps \to 0$ of \eqref{Intro:NSCH} with the choice $A = \frac{\beta}{\eps}$, $B = \beta \eps$ is also stated, which can be deduced following the formally matched asymptotic analysis presented in \cite{AGG,DGL,GLSS}.
\item In Section \ref{sec:Analysis}, we study the simplified model variant \eqref{Intro:ZEED} and introduce preliminary results that will be used in our analysis, as well as a summary of our main mathematical results.
\item In Section \ref{sec:weak} we establish the global existence of weak solutions in two and three dimensions under general assumptions on the source terms $\Gamma$ and $S$, the potential $\Psi$, the viscosity $\eta$ and the mobilities $m$ and $n$.  The assertion is given in Theorem \ref{thm:ZEED:weaksoln}.
\item In Section \ref{sec:ZEED:2D:strong}, under additional assumptions on the initial conditions, the source term $\Gamma$, the mobilities $m$ and  $n$ and the viscosity $\eta$, we prove the global existence of strong solutions to \eqref{Intro:ZEED} in two dimensions.  The assertion is given in Theorem \ref{thm:2D:strong}.
\item In Section \ref{sec:ctsdep2D}, we show the continuous dependence of strong solutions to \eqref{Intro:ZEED} in two dimensions on the initial data and source terms.  This is stated in Theorem \ref{thm:2D:ctsdep}, and the proof is based on a similar procedure used in \cite{GLDirichlet,GLNeumann}, which allow us to deduce continuous dependence for the chemical potential $\mu$ in $L^{2}(0,T;L^{2}(\Omega))$ and for the order parameter $\varphi$ in $L^{2}(0,T;H^{2}(\Omega))$ as well.  In particular, thanks to the regularities of the strong solutions, the continuous dependence result is valid for variable fluid viscosity and mobilities.
\end{itemize}

\section{Model derivation}\label{sec:volave}

\subsection{Balance laws}
We begin with the following modeling assumption:  in a bounded domain $\Omega \subset \R^{d}$, $d = 2,3$, there is a two-component mixture of immiscible fluids and a chemical species.  Let $\overline{\rho}_{i}$, $i = 1,2$, denote the constant mass density a pure component $i$, while $\rho_{i}$ denote the actual mass of the component $i$ per volume in the mixture.  We define the volume fraction $u_{i} := \rho_{i} / \overline{\rho}_{i}$ for $i = 1,2$, which lies in between the physical interval $[0,1]$.  Assuming the excess volume of mixing is zero, i.e., the fluids behave like a simple mixture, it holds that
\begin{align}\label{sum:ui}
u_{1} + u_{2} = 1.
\end{align}
We define the mass density $\rho$ of the mixture as $\rho = \rho_{1} + \rho_{2}$, and introduce the volume-averaged velocity
\begin{align}\label{defn:volavervelo}
\bm{v} := u_{1} \bm{v}_{1} + u_{2} \bm{v}_{2},
\end{align}
where $\bm{v}_{i}$, $i = 1,2$, is the individual velocity for component $i$.  The balance laws for mass and linear momentum are given as
\begin{subequations}
\begin{alignat}{3}
\pd_{t} \rho_{i} + \div (\rho_{i} \bm{v}_{i}) - \Gamma_{i} & = 0, \label{proto:density} \\
\pd_{t}(\rho \bm{v}) + \div (\rho \bm{v} \otimes \bm{v}) - \div (\bm{S} - p \Id) & = 0, \label{proto:mom}
\end{alignat}
\end{subequations}
where $\Gamma_{i}$ denotes a source term for component $i$, $\bm{S}$ is a symmetric tensor (due to conservation of angular momentum) yet to be identified, $p$ denotes the pressure, $\Id$ is the identity tensor, and for two vectors $\bm{a}, \bm{b} \in \R^{d}$, the tensor product $\bm{a} \otimes \bm{b}$ is defined as $\bm{a} \otimes \bm{b} = (a_{i} b_{j})_{1 \leq i,j \leq d}$.  Furthermore the divergence of a tensor $\bm{A}$ is taken row-wise.  Note that upon dividing \eqref{proto:density} by $\overline{\rho}_{i}$ and using \eqref{sum:ui} and \eqref{defn:volavervelo}, we obtain
\begin{align}\label{proto:div}
\div \bm{v} = \frac{\Gamma_{1}}{\overline{\rho}_{1}} + \frac{\Gamma_{2}}{\overline{\rho}_{2}} =: U_{v}.
\end{align}
Introducing the fluxes
\begin{align}\label{defn:flux}
\bm{J}_{i} := \rho_{i} (\bm{v}_{i} - \bm{v}), \quad \bm{J} := \bm{J}_{1} + \bm{J}_{2}, \quad \bm{J}_{\varphi} := \frac{\bm{J}_{2}}{\overline{\rho}_{2}} - \frac{\bm{J}_{1}}{\overline{\rho}_{1}}, \quad \bm{K} := \frac{\bm{J}_{1}} {\overline{\rho}_{1}} + \frac{\bm{J}_{2}}{\overline{\rho}_{2}}.
\end{align}
Then, by the definition of $\bm{K}$, $\bm{J}_{i}$ and the volume-averaged velocity $\bm{v}$, it holds that
\begin{align}\label{K=0}
\bm{K} = u_{1}(\bm{v}_{1} - \bm{v}) + u_{2}(\bm{v}_{2} - \bm{v}) = \bm{v} - (u_{1} + u_{2}) \bm{v} = \bm{0}.
\end{align}
Adding \eqref{proto:density} and using \eqref{defn:flux} yields the balance law for the mass density $\rho$:
\begin{align}\label{proto:massdensity}
\pd_{t}\rho + \div (\rho \bm{v}) + \div \bm{J} = \Gamma_{1} + \Gamma_{2} =: \Theta.
\end{align}
From \eqref{proto:density} we also deduce that
\begin{align}\label{equ:ui}
\pd_{t}u_{i} + \div (u_{i} \bm{v}) + \frac{1}{\overline{\rho}_{i}} \div \bm{J}_{i} = \frac{\Gamma_{i}}{\overline{\rho}_{i}}.
\end{align}
Defining the order parameter $\varphi$ as the difference in volume fractions
\begin{align}\label{defn:varphi}
\varphi := u_{2} - u_{1},
\end{align}
so that the fluid $2$ is represented by the region $\{ \varphi = 1\}$ and fluid $1$ is represented by the region $\{\varphi = -1 \}$, we can deduce from subtracting \eqref{equ:ui} for $i = 1$ from \eqref{equ:ui} for $i = 2$ that
\begin{align}\label{proto:varphi}
\pd_{t}\varphi + \div (\varphi \bm{v}) + \div \bm{J}_{\varphi} = \frac{\Gamma_{2}}{\overline{\rho}_{2}} - \frac{\Gamma_{1}}{\overline{\rho}_{1}} =: U_{\varphi}.
\end{align}
Furthermore, the mass density $\rho$ can be expressed as an affine function of $\varphi$:
\begin{align}\label{massdensity:varphi}
\rho(\varphi) = \overline{\rho}_{1} u_{1} + \overline{\rho}_{2} u_{2} = \frac{\overline{\rho}_{2} - \overline{\rho}_{1}}{2} \varphi + \frac{\overline{\rho}_{2} + \overline{\rho}_{1}}{2}.
\end{align}
The density of the chemical species present in $\Omega$, denoted by $\sigma$, satisfies the following balance law
\begin{align}\label{proto:sigma}
\pd_{t} \sigma + \div (\sigma \bm{v}) + \div \bm{J}_{\sigma} = S,
\end{align}
for some source term $S$ and flux $\bm{J}_{\sigma}$ which we will determine later.  The prototype model consists of equations \eqref{proto:mom}, \eqref{proto:massdensity}, \eqref{proto:varphi} and \eqref{proto:sigma}, and we will employ the Lagrange multiplier method of M\"{u}ller and Liu (see \cite[Chap. 7]{Liu} and \cite[\S 2.2]{AGG}) to determine the constitutive assumptions on $\bm{S}$, $\bm{J}_{\varphi}$ and $\bm{J}_{\sigma}$ so that the resulting model is thermodynamically consistent.  Note that by \eqref{defn:flux} and \eqref{K=0}
\begin{align*}
\bm{J} = \bm{J}_{1} + \bm{J}_{2} = \frac{\overline{\rho}_{2} - \overline{\rho}_{1}}{2} \bm{J}_{\varphi} + \frac{\overline{\rho}_{2} + \overline{\rho}_{1}}{2} \bm{K} = \frac{\overline{\rho}_{2} - \overline{\rho}_{1}}{2} \bm{J}_{\varphi},
\end{align*}
so that we can rewrite \eqref{proto:massdensity} as
\begin{align}\label{alt:proto:massdensity}
\pd_{t}\rho + \div (\rho \bm{v}) + \frac{\overline{\rho}_{2} - \overline{\rho}_{1}}{2} \div \bm{J}_{\varphi} = \Theta.
\end{align}
Thus, it is sufficient to deduce $\bm{J}_{\varphi}$ in order to obtain the equation for the mass density $\rho$.

\subsection{Energy inequality}\label{sec:EnergyIneq}
We postulate a general energy density of the form
\begin{align}\label{energy}
e = \frac{\rho}{2} \abs{\bm{v}}^{2} + A \Psi(\varphi) + \frac{B}{2} \abs{\nabla \varphi}^{2} + N(\varphi, \sigma),
\end{align}
where $A,B$ are positive constants and $\Psi$ is a potential with equal minima at $\pm 1$.  The first term of $e$ is the kinetic energy, the second and third term of $e$ form the so-called Ginzburg--Landau energy which accounts for the interfacial energy and unmixing tendencies.  The last term of $e$ is the free energy of the chemical species, which we incorporate possible interactions with the fluids by introducing a dependence on $\varphi$.  In contrast to \cite{GLSS}, here we include the effects of inertia that leads to the appearance of the kinetic energy density in \eqref{energy}.  We mention that this effect has also been considered in the thesis \cite{Sitka}, but below we will derive a more general model.

We now derive the model based on a dissipation inequality for balance laws with source terms that has been used in \cite{GLSS,Gurtin89,Gurtin96,Podio}, see also \cite[Chap. 62]{Gurtinbook}.  The second law of thermodynamics in the isothermal situation can be formulated as a free energy inequality, which states that for all volumes $V(t) \subset \Omega$ transported by the fluid velocity, the following inequality has to hold
\begin{align*}
\frac{\dd}{\dt} \int_{V(t)} e \dx + \int_{\pd V(t)} \bm{J}_{e} \cdot \bm{\nu} \dHaus - \int_{V(t)} c_{\varphi} U_{\varphi} + c_{v} U_{v} + c_{S} S \dx \leq 0,
\end{align*}
where $\bm{\nu}$ is the outer unit normal on $\pd V(t)$ and $\bm{J}_{e}$ is an energy flux yet to be determined.  In the above, $\dx$ and $\dHaus$ denote integration with respect to the $d$ dimensional Lebesgue measure and $(d-1)$ dimensional Hausdorff measure, respectively.  The source terms $U_{\varphi}$, $U_{v}$ and $S$ carry with them an energy contribution of the form
\begin{align*}
\int_{V(t)} c_{\varphi} U_{\varphi} + c_{v} U_{v} + c_{S} S \dx,
\end{align*}
where the prefactors $c_{\varphi}, c_{v}, c_{S}$ will be determined below.

By the Reynold transport theorem and divergence theorem \cite{Gurtinbook}, this leads to the following local inequality
\begin{align}\label{energy:diss:local}
\pd_{t} e + \div (e \bm{v}) + \div \bm{J}_{e} - c_{\varphi} U_{\varphi} - c_{v} U_{v} - c_{S} S \leq 0.
\end{align}
Instead of asking \eqref{energy:diss:local} to hold only for variables $(\varphi, \sigma, \bm{v}, U_{\varphi}, U_{v}, S)$ that satisfy the prototype model equations \eqref{proto:mom}, \eqref{proto:varphi}, and \eqref{proto:sigma}, the Lagrange multiplier method of M\"{u}ller and Liu  relaxes the constraint on the model equations by introducing Lagrange multipliers (see, e.g., \cite[Chap. 7]{Liu}), in our case $\lambda_{\varphi}$ and $\lambda_{\sigma}$ for \eqref{proto:varphi} and \eqref{proto:sigma}, and instead asks that
\begin{equation}\label{ML:energy}
\begin{aligned}
-\mathcal{D} & :=  \pd_{t} e + \div (e \bm{v}) + \div \bm{J}_{e} - c_{\varphi} U_{\varphi} - c_{v} U_{v} - c_{S} S \\
&\quad \  - \lambda_{\varphi} \left ( \pd_{t} \varphi + \div (\varphi \bm{v}) + \div \bm{J}_{\varphi} - U_{\varphi} \right ) \\
&\quad \  - \lambda_{\sigma} \left ( \pd_{t} \sigma + \div (\sigma \bm{v}) + \div \bm{J}_{\sigma} - S \right ) \leq 0
\end{aligned}
\end{equation}
holds for arbitrary $(\varphi, \sigma, \bm{v}, U_{\varphi}, U_{v}, S, \md \varphi, \md \sigma)$, where the material derivative $\md f$ of a function $f$ is defined as
\begin{align*}
\md f := \pd_{t} f + \bm{v} \cdot  \nabla f.
\end{align*}
We introduce the notations
\begin{align*}
N_{\varphi} := \frac{\pd N}{\pd \varphi}, \quad N_{\sigma} := \frac{\pd N}{\pd \sigma},\quad \mu := A \Psi'(\varphi) - B \Laplace \varphi + N_{\varphi},
\end{align*}
and recall the identities
\begin{align}
\md \varphi \nabla \varphi & = \pd_{t} \varphi \nabla \varphi +  (\nabla \varphi \otimes \nabla \varphi)\bm{v},\\
\nabla \varphi \cdot \md (\nabla \varphi) & = \div (\md \varphi \nabla \varphi) - \md \varphi \Laplace \varphi - \nabla \bm{v} : (\nabla \varphi \otimes \nabla \varphi),\label{md:grad:exchange}
\end{align}
where for vectors $\bm{a}, \bm{b} \in \R^{d}$, $\nabla \bm{a} : \nabla \bm{b} := \sum_{i,j=1}^{d} \pd_{i} a_{j} \pd_{i} b_{j}$.
Then using \eqref{proto:mom}, \eqref{proto:div} and \eqref{proto:massdensity}, we compute that
\begin{equation}\label{kinetic:part}
\begin{aligned}
& \md \left ( \frac{\rho}{2} \abs{\bm{v}}^{2} \right ) 
= (\pd_{t}(\rho \bm{v}) + \div (\rho \bm{v} \otimes \bm{v}) - \rho (\div \bm{v}) \bm{v} ) \cdot \bm{v} - \frac{\abs{\bm{v}}^{2}}{2} \md \rho \\
& \quad = \div (\bm{S} - p \Id) \cdot \bm{v} + \frac{\abs{\bm{v}}^{2}}{2} \left ( \div \bm{J} - \Theta - \rho U_{v} \right) \\
& \quad = \div \left ((\bm{S}^{\top} - p \Id) \bm{v} + \frac{\abs{\bm{v}}^{2}}{2} \bm{J} \right )  - \nabla \bm{v} : \left ( \bm{S} + \bm{v} \otimes \bm{J} \right )
- \frac{\abs{\bm{v}}^{2}}{2} \left ( \Theta + \rho U_{v} \right ) + p U_{v} .
\end{aligned}
\end{equation}
Furthermore, using the relations
\begin{align*}
\Gamma_{1} = \frac{\overline{\rho}_{1}}{2} (U_{v} - U_{\varphi}), \quad \Gamma_{2} = \frac{\overline{\rho}_{2}}{2} (U_{v} + U_{\varphi}) \quad \text{such that} \quad \Theta = \frac{\overline{\rho}_{2} + \overline{\rho}_{1}}{2} U_{v} + \frac{\overline{\rho}_{2} - \overline{\rho}_{1}}{2} U_{\varphi},
\end{align*}
we can express \eqref{kinetic:part} as
\begin{align*}
\md \left ( \frac{\rho}{2} \abs{\bm{v}}^{2} \right ) & =  \div \left ((\bm{S} - p \Id)^{\top} \bm{v} + \frac{\abs{\bm{v}}^{2}}{2} \bm{J} \right ) - \nabla \bm{v} : \left ( \bm{S}  + \bm{v} \otimes \bm{J} \right ) \\
& \quad - \left ( \left ( \rho + \frac{\overline{\rho}_{2} + \overline{\rho}_{1}}{2} \right ) \frac{\abs{\bm{v}}^{2}}{2} - p \right )U_{v} - \frac{\abs{\bm{v}}^{2}}{2} \frac{\overline{\rho}_{2} - \overline{\rho}_{1}}{2} U_{\varphi}.
\end{align*}
Then after a long calculation, we obtain from \eqref{ML:energy} that
\begin{equation}\label{Dissipation}
\begin{aligned}
-\mathcal{D} & =  \div \left ( \bm{J}_{e} - \lambda_{\varphi} \bm{J}_{\varphi} - \lambda_{\sigma} \bm{J}_{\sigma} + (\bm{S}^{\top} - p \Id) \bm{v} + \frac{\abs{\bm{v}}^{2}}{2} \bm{J} + B \md \varphi \nabla \varphi \right ) + \bm{J}_{\varphi} \cdot \nabla \lambda_{\varphi} \\
& \quad + \bm{J}_{\sigma} \cdot \nabla \lambda_{\sigma} + \left ( \mu - \lambda_{\varphi} \right ) \md \varphi + \left ( N_{\sigma} - \lambda_{\sigma} \right ) \md \sigma + (\lambda_{\sigma} - c_{S}) S \\
& \quad + \left (\lambda_{\varphi} - c_{\varphi} - \frac{\overline{\rho}_{2} - \overline{\rho}_{1}}{2} \frac{\abs{\bm{v}}^{2}}{2} \right ) U_{\varphi} - \nabla \bm{v} : \left ( \bm{S} + \bm{v} \otimes \bm{J} + B \nabla \varphi \otimes \nabla \varphi \right ) \\
& \quad + U_{v} \left ( f(\varphi, \nabla \varphi, \sigma) - c_{v} - \lambda_{\varphi} \varphi - \lambda_{\sigma} \sigma - \frac{\overline{\rho}_{2} + \overline{\rho}_{1}}{2} \frac{\abs{\bm{v}}^{2}}{2} + p \right ),
\end{aligned}
\end{equation}
where we recall that $f(\varphi, \nabla \varphi, \sigma) =  A \Psi(\varphi) + \frac{B}{2} \abs{\nabla \varphi}^{2} + N(\varphi, \sigma)$ denotes the free energy.

\subsection{Constitutive assumptions and the general model}
In order for $-\mathcal{D} \leq 0$ to hold for arbitrary $(\varphi, \sigma, \bm{v}, U_{\varphi}, U_{v}, S, \md \varphi, \md \sigma)$, where $-\mathcal{D}$ is given in \eqref{Dissipation}, we make the following constitutive assumptions:
\begin{subequations}\label{constitutive}
\begin{align}
\lambda_{\varphi} & = \mu, \quad \lambda_{\sigma} = N_{\sigma}, \label{con:LM} \\
\bm{J}_{e} & = \mu \bm{J}_{\varphi} + N_{\sigma} \bm{J}_{\sigma} - \bm{S}^{\top} \bm{v} + p \bm{v} - \frac{\abs{\bm{v}}^{2}}{2} \frac{\overline{\rho}_{2} - \overline{\rho}_{1}}{2} \bm{J}_{\varphi} - B \md \varphi \nabla \varphi, \label{con:Je} \\
\bm{S} & = 2 \eta(\varphi) \der \bm{v} - \bm{v} \otimes \bm{J} - B \nabla \varphi \otimes \nabla \varphi, \label{con:S}\\
\bm{J}_{\varphi} & = -m(\varphi) \nabla \mu, \quad \bm{J}_{\sigma} = -n(\varphi) \nabla N_{\sigma}, \label{con:flux} \\
c_{S} & = N_{\sigma}, \quad c_{\varphi} = \mu - \frac{\overline{\rho}_{2} - \overline{\rho}_{1}}{2} \frac{\abs{\bm{v}}^{2}}{2}, \label{con:cS} \\
c_{v} & = A \Psi(\varphi) + \frac{B}{2} \abs{\nabla \varphi}^{2} + N(\varphi, \sigma) - \mu \varphi - N_{\sigma} \sigma - \frac{\overline{\rho}_{2} + \overline{\rho}_{1}}{2} \frac{\abs{\bm{v}}^{2}}{2} + p, \label{con:cv}
\end{align}
\end{subequations}
where $\der \bm{v} := \frac{1}{2} (\nabla \bm{v} + (\nabla \bm{v})^{\top})$ is the rate of deformation tensor, $\eta(\varphi)$ is the mixture viscosity, and $m(\varphi), n(\varphi)$ are positive mobilities.  Equations \eqref{con:LM}--\eqref{con:cv} yield that the right-hand side of \eqref{Dissipation} is non-positive, namely,
\begin{align*}
-\mathcal{D} = -m(\varphi) \abs{\nabla \mu}^{2} - n(\varphi) \abs{\nabla N_{\sigma}}^{2} - 2 \eta(\varphi) \abs{\der \bm{v}}^{2} \leq 0.
\end{align*}
We remark that the constitutive assumption \eqref{con:Je} for the energy flux $\bm{J}_{e}$ is chosen so that the divergence term in \eqref{Dissipation} vanishes.  The term $(p \Id - \bm{S}^{\top}) \bm{v}$ accounts for energy change due to work done by macroscopic stresses (see \cite{AGG}), while energy flux due to mass diffusion are described by the terms $\mu \bm{J}_{\varphi}$ and $N_{\sigma} \bm{J}_{\sigma}$.  Changes of kinetic energy due to mass diffusion is given by $\frac{1}{2} \frac{\overline{\rho}_{2} - \overline{\rho}_{1}}{2} \abs{\bm{v}}^{2} \bm{J}_{\varphi} = \frac{1}{2} \abs{\bm{v}}^{2} \bm{J}$, and the term $\md \varphi \nabla \varphi$ arises from the moving phase boundaries.

The above constitutive assumptions lead to the following Navier--Stokes--Cahn--Hilliard model with mass transfer and chemical coupling:
\begin{subequations}\label{NSCH}
\begin{alignat}{2}
\pd_{t} \rho + \div \left ( \rho \bm{v} - \frac{\overline{\rho}_{2} - \overline{\rho}_{1}}{2} m(\varphi) \nabla \mu \right )  & = \Theta = \Gamma_{1} + \Gamma_{2}, \label{mass} \\
\div \bm{v}  & = U_{v} = \frac{\Gamma_{1}}{\overline{\rho}_{1}} + \frac{\Gamma_{2}}{\overline{\rho}_{2}}, \label{div} \\
\pd_{t} (\rho \bm{v}) + \div \left (\rho \bm{v} \otimes \bm{v} - \frac{\overline{\rho}_{2} - \overline{\rho}_{1}}{2} m(\varphi) \bm{v} \otimes \nabla \mu \right )  & =  - \nabla p + \div \left ( 2 \eta(\varphi) \der \bm{v} \right ) \label{mom} \\
\notag & \quad \ - \div \left ( B \nabla \varphi \otimes \nabla \varphi \right ), \\
\pd_{t} \varphi + \div (\varphi \bm{v} - m(\varphi) \nabla \mu) & =  U_{\varphi} = \frac{\Gamma_{2}}{\overline{\rho}_{2}} -  \frac{\Gamma_{1}}{\overline{\rho}_{1}}, \label{varphi} \\
\mu & = A \Psi'(\varphi) - B \Laplace \varphi + N_{\varphi}, \label{mu} \\
\pd_{t} \sigma + \div (\sigma \bm{v} - n(\varphi) \nabla N_{\sigma}) & = S. \label{sigma}
\end{alignat}
\end{subequations}
Formally, we can obtain an energy identity by integrating \eqref{Dissipation} over $\Omega$, where we recall \eqref{ML:energy}, the constitutive assumptions \eqref{con:LM}--\eqref{con:cv}, and now $\varphi$ and $\sigma$ satisfy \eqref{proto:varphi} and \eqref{proto:sigma}, which is given by
\begin{equation}\label{Energy:Id}
\begin{aligned}
& \frac{\dd}{\dt} \int_{\Omega} \frac{\rho}{2} \abs{\bm{v}}^{2} + A \Psi(\varphi) + \frac{B}{2} \abs{\nabla \varphi}^{2} + N(\varphi, \sigma) \dx \\
& \ + \int_{\Omega} m(\varphi) \abs{\nabla \mu}^{2} + n(\varphi) \abs{\nabla N_{\sigma}}^{2} + 2 \eta(\varphi) \abs{\der \bm{v}}^{2} \dx \\
& \ + \int_{\Omega} U_{v} \left ( \varphi \mu + \sigma N_{\sigma} + \frac{\overline{\rho}_{2} + \overline{\rho}_{1}}{2} \frac{\abs{\bm{v}}^{2}}{2} - p - A \Psi(\varphi) - \frac{B}{2} \abs{\nabla \varphi}^{2} - N(\varphi, \sigma) \right ) \dx \\
& \ - \int_{\Omega} \left ( \mu - \frac{\overline{\rho}_{2} - \overline{\rho}_{1}}{2} \frac{\abs{\bm{v}}^{2}}{2} \right ) U_{\varphi} + N_{\sigma} S \dx \\
& \ + \int_{\pd \Omega} \left (\frac{\rho}{2} \abs{\bm{v}}^{2}  + p + A \Psi(\varphi) + \frac{B}{2} \abs{\nabla \varphi}^{2} + N(\varphi, \sigma) \right ) \bm{v} \cdot \bm{\nu} - n(\varphi) N_{\sigma} \pdnu N_{\sigma}  \dHaus \\
& \ - \int_{\pd \Omega} m(\varphi) \left ( \mu + \frac{\abs{\bm{v}}^{2}}{2}  \frac{\overline{\rho}_{2} - \overline{\rho}_{1}}{2} \right ) \pdnu \mu + 2 \eta(\varphi) (\der \bm{v}) \bm{v} \cdot \bm{\nu} + B \pd_t \varphi \pdnu \varphi \dHaus = 0,
\end{aligned}
\end{equation}
where $\pdnu f := \nabla f \cdot \bm{\nu}$ is the normal derivative of $f$ on $\pd \Omega$.
\subsection{Modified pressure and reformulation of the momentum equation}\label{sec:Reformulation}
We now present three reformulations of the pressure and the corresponding momentum equation \eqref{mom}.
\begin{itemize}
\item Define $q := p + A \Psi(\varphi) + \frac{B}{2} \abs{\nabla \varphi}^{2}$, so that
\begin{align*}
- \nabla p - B \div (\nabla \varphi \otimes \nabla \varphi) = - \nabla q + (A\Psi'(\varphi) - B \Laplace \varphi) \nabla \varphi = - \nabla q + (\mu - N_{\varphi}) \nabla \varphi,
\end{align*}
and \eqref{mom} becomes
\begin{align}\label{alt:mom:q}
\pd_{t} (\rho \bm{v}) + \div \left (\rho \bm{v} \otimes \bm{v} - \frac{\overline{\rho}_{2} - \overline{\rho}_{1}}{2} m(\varphi) \bm{v} \otimes \nabla \mu \right )  = & - \nabla q + \div \left ( 2 \eta(\varphi) \der \bm{v} \right ) \\
\notag & + (\mu - N_{\varphi}) \nabla \varphi,
\end{align}
while the prefactor $\lambda_{v} := - c_{v}$ that is multiplied with $U_{v}$ in the energy identity \eqref{Energy:Id} now reads as
\begin{align}\label{alt:cv:q}
\lambda_{v} = \varphi \mu + \sigma N_{\sigma} +  \frac{\overline{\rho}_{2} + \overline{\rho}_{1}}{2}  \frac{\abs{\bm{v}}^{2}}{2} - q - N(\varphi, \sigma).
\end{align}
\item Define $r := p + A \Psi(\varphi) + \frac{B}{2} \abs{\nabla \varphi}^{2} + N(\varphi, \sigma)$, so that
\begin{align*}
-\nabla p - B \div (\nabla \varphi \otimes \nabla \varphi) = - \nabla r + \mu \nabla \varphi + N_{\sigma} \nabla \sigma,
\end{align*}
and \eqref{mom} becomes
\begin{align}\label{alt:mom:r}
\pd_{t} (\rho \bm{v}) + \div \left (\rho \bm{v} \otimes \bm{v} - \frac{\overline{\rho}_{2} - \overline{\rho}_{1}}{2} m(\varphi) \bm{v} \otimes \nabla \mu \right )  = & - \nabla r + \div \left ( 2 \eta(\varphi) \der \bm{v} \right ) \\
\notag & + \mu  \nabla \varphi + N_{\sigma} \nabla \sigma,
\end{align}
while the prefactor $\lambda_{v}$ becomes
\begin{align}\label{alt:cv:r}
\lambda_{v} = \varphi \mu + \sigma N_{\sigma} + \frac{\overline{\rho}_{2} + \overline{\rho}_{1}}{2}   \frac{\abs{\bm{v}}^{2}}{2} - r.
\end{align}
\item Define $s := p + A \Psi(\varphi) + \frac{B}{2} \abs{\nabla \varphi}^{2} + N(\varphi, \sigma) - \varphi \mu - N_{\sigma} \sigma$, so that
\begin{align*}
-\nabla p - B \div (\nabla \varphi \otimes \nabla \varphi) = - \nabla s - \varphi \nabla \mu - \sigma \nabla N_{\sigma},
\end{align*}
and \eqref{mom} becomes
\begin{align}\label{alt:mom:s}
\pd_{t} (\rho \bm{v}) + \div \left (\rho \bm{v} \otimes \bm{v} - \frac{\overline{\rho}_{2} - \overline{\rho}_{1}}{2} m(\varphi) \bm{v} \otimes \nabla \mu \right )  = & - \nabla s + \div \left ( 2 \eta(\varphi) \der \bm{v} \right ) \\
\notag & - \varphi \nabla \mu - \sigma \nabla N_{\sigma},
\end{align}
while the prefactor $\lambda_{v}$ becomes
\begin{align}\label{alt:cv:s}
\lambda_{v} = \frac{\overline{\rho}_{2} + \overline{\rho}_{1}}{2}   \frac{\abs{\bm{v}}^{2}}{2} - s.
\end{align}
\end{itemize}
\begin{remark}
In our theoretical analysis below, we will employ the reformulation \eqref{alt:mom:q} for the momentum equation.
\end{remark}

\subsection{Reduction to special models}
\subsubsection{Absence of the chemical substance}
In the absence of chemical substances, by setting $\sigma = 0$, $N(\varphi,\sigma) = 0$, and $S = 0$, we obtain a Navier--Stokes--Cahn--Hilliard system with source terms that is composed of \eqref{mass}--\eqref{varphi} and
\begin{align*}
\mu = A \Psi'(\varphi) - B \Laplace \varphi,
\end{align*}
which can be seen as the Navier--Stokes analogue of \cite[(2.35)]{GLSS}.  Furthermore, if $\Gamma_{1} = \Gamma_{2} = 0$, then we obtain the Navier--Stokes--Cahn--Hilliard model of Abels, Garcke and Gr\"{u}n \cite{AGG}. For the existence of global weak solutions to the Navier--Stokes--Cahn--Hilliard models, we refer to, for instances \cite{Abels,ADG1,ADG2,AGW,Boyer,GG10, GG10a,Weber,ZWH}, where the reformulation \eqref{alt:mom:r} of the momentum equation has been used.

\subsubsection{Zero excess of total mass}
We consider the case $$\Gamma_{2} = - \Gamma_{1} =: \Gamma.$$
Define
\begin{align*}
\alpha = \frac{1}{\overline{\rho}_{2}} - \frac{1}{\overline{\rho}_{1}}, \quad \rho_{S} := \frac{1}{\overline{\rho}_{2}} + \frac{1}{\overline{\rho}_{1}} \, \Longrightarrow \,
\Theta = 0, \quad U_{v} = \alpha \Gamma, \quad U_{\varphi} = \rho_{S} \Gamma.
\end{align*}
This leads to the Navier--Stokes analogue of \cite[(2.33)]{GLSS} which consists of equations \eqref{mom}, \eqref{mu}, \eqref{sigma} and
\begin{subequations}\label{Zero:Excess}
\begin{align}
\pd_{t} \rho + \div \left ( \rho \bm{v} - \frac{\overline{\rho}_{2} - \overline{\rho}_{1}}{2} m(\varphi) \nabla \mu \right ) & = 0, \\
\div \bm{v} & = \alpha \Gamma, \\
\pd_{t} \varphi + \div (\varphi \bm{v} - m(\varphi) \nabla \mu) & = \rho_{S} \Gamma.
\end{align}
\end{subequations}
Furthermore, in the case of equal densities $\overline{\rho}_{1} = \overline{\rho}_{2} = \rho_{*}$ so that $\rho = \rho_{*}$ is constant, then
\begin{align*}
\alpha = 0,\quad \rho_{S} = \frac{2}{\rho_{*}},
\end{align*}
and the resulting system becomes
\begin{subequations}\label{ZeroExcess:equal:densities}
\begin{align}
\rho_{*} \left ( \pd_{t} \bm{v} + (\bm{v} \cdot \nabla ) \bm{v} \right ) & = - \nabla p + \div (2 \eta(\varphi) \der \bm{v} - B \nabla \varphi \otimes \nabla \varphi), \\
\div \bm{v} & = 0, \\
\pd_{t} \varphi + \bm{v} \cdot \nabla \varphi & = \div (m(\varphi) \nabla \mu) + \frac{2}{\rho_{*}} \Gamma, \\
\mu & = A \Psi'(\varphi) - B \Laplace \varphi + N_{\varphi}, \\
\pd_{t} \sigma + \bm{v} \cdot \nabla \sigma & = \div (n(\varphi) \nabla N_{\sigma}) + S,
\end{align}
\end{subequations}
which is the Navier--Stokes analogue of \cite[(2.34)]{GLSS}.  If in addition, we set $\bm{v} = \bm{0}$, then neglecting the first two equations of \eqref{ZeroExcess:equal:densities} leads to \cite[(2.36)]{GLSS}.

\subsubsection{Scaled zero excess of total mass}
Alternatively, we consider the case
\begin{align*}
\frac{\Gamma_{2}}{\overline{\rho}_{2}} = -\frac{\Gamma_{1}}{\overline{\rho}_{1}},
\end{align*}
then it holds that
\begin{align*}
U_{v} = 0, \quad U_{\varphi} = \frac{2}{\overline{\rho}_{2}} \Gamma_{2}, \quad \Theta = \left ( 1 - \frac{\overline{\rho}_{1}}{\overline{\rho}_{2}} \right ) \Gamma_{2},
\end{align*}
which leads to the model consisting of \eqref{mom}, \eqref{mu}, \eqref{sigma} and
\begin{subequations}\label{Scaled:ZeroExcess}
\begin{align}
\pd_{t} \rho + \div \left (\rho \bm{v} - \frac{\overline{\rho}_{2} - \overline{\rho}_{1}}{2} m(\varphi) \nabla \mu \right ) & = \left ( 1 - \frac{\overline{\rho}_{1}}{\overline{\rho}_{2}} \right ) \Gamma_{2}, \\
\div \bm{v} & = 0,\\
\pd_{t} \varphi + \div (\varphi \bm{v} - m(\varphi) \nabla \mu ) & = \frac{2}{\overline{\rho}_{2}} \Gamma_{2}. 
\end{align}
\end{subequations}
\begin{remark}
The interesting feature of \eqref{Scaled:ZeroExcess} is that the densities are allowed to be different, i.e., $\overline{\rho}_{1} \neq \overline{\rho}_{2}$ while we retain a solenoid velocity field.
\end{remark}

\subsubsection{No mass exchange for fluid 1}
Setting $$\Gamma_{1} = 0, \quad \Gamma := \frac{\Gamma_{2}}{\overline{\rho}_{2}}$$
leads to the model consisting of \eqref{mom}, \eqref{mu}, \eqref{sigma} and
\begin{subequations}\label{Zero:Healthy:Source}
\begin{align}
\pd_{t} \rho + \div \left ( \rho \bm{v} - \frac{\overline{\rho}_{2} - \overline{\rho}_{1}}{2} m(\varphi) \nabla \mu \right )  & = \overline{\rho}_{2} \Gamma, \\
\div \bm{v}  & = \Gamma,  \label{div:ZHS} \\
\pd_{t} \varphi + \div (\varphi \bm{v} - m(\varphi) \nabla \mu) & =  \Gamma. \label{varphi:ZHS} 
\end{align}
\end{subequations}
Neglecting $\sigma$, the resulting model is the Navier--Stokes analogue of \cite[(2.43)]{GLSS}.  Note that the source terms for the divergence equation \eqref{div:ZHS} and the equation \eqref{varphi:ZHS} are equal, and this feature has also appeared in the works of \cite{JWZ} (see \cite{LTZ,WW2012,WZ2013} for the case $\Gamma=0$).

\subsubsection{Recovering Darcy's law}\label{sec:Darcy}
A Darcy's law can be obtained from the momentum equation \eqref{mom} by performing a proper averaging procedure in a Hele--Shaw geometry.  We assume that the domain $\Omega$ is given as $\Omega = U \times [0,\delta]$ for some positive constant $\delta \ll 1$ and bounded domain $U \subset \R^{d-1}$.  Rescaling the pressure $p$ and the constant $B$ appropriately with $\delta$, and employing formal asymptotic expansions of the variables $(\rho, \bm{v}, \varphi, \mu, \sigma, \Gamma_{i})$ in $\delta$, by examining the equations order by order in $\delta$, one can derive the Cahn--Hilliard--Darcy model of \cite{GLSS} from the Navier--Stokes--Cahn--Hilliard system \eqref{NSCH}.  For further details, we refer to the recent work \cite{DGL}, in which a Hele--Shaw--Cahn--Hilliard model is derived from the Navier--Stokes--Cahn--Hilliard model of \cite{AGG}, see also \cite[Chap. 4]{Ockendon}.

\subsection{Some analytical issues}\label{sec:analy:issue}

\subsubsection{Boundary conditions for the velocity}
We now examine the boundary term involving the velocity $\bm{v}$ in the energy identity \eqref{Energy:Id}, which reads as
\begin{equation}\label{bdyterm}
\begin{aligned}
 \int_{\pd \Omega} \left (\frac{\rho}{2} \abs{\bm{v}}^{2}  + p + A \Psi(\varphi) + \frac{B}{2} \abs{\nabla \varphi}^{2} + N(\varphi, \sigma) \right ) \bm{v} \cdot \bm{\nu} - 2 \eta(\varphi) (\der \bm{v}) \bm{v} \cdot \bm{\nu} \dHaus.
\end{aligned}
\end{equation}
Natural boundary conditions for the velocity are
\begin{equation*}
\begin{alignedat}{5}
&\text{no-slip b.c.} &&\quad && \bm{v} = \bm{0} && \quad && \text{ on } \pd \Omega,\\
&\text{or free-slip b.c.} && \quad  && \bm{v} \cdot \bm{\nu} = 0\quad \text{and}\quad \nabla \times(\bm{v}\times \bm{\nu})=0 && \quad  &&\text{ on } \pd \Omega,
\end{alignedat}
\end{equation*}
so that \eqref{bdyterm} would vanish in both cases.  However, integrating \eqref{div} and applying the divergence theorem yields that
\begin{align}\label{compatibility}
0 = \int_{\pd \Omega} \bm{v} \cdot \bm{\nu} \dHaus = \int_{\Omega} \div \bm{v} \dx = \int_{\Omega} U_{v} \dx,
\end{align}
which implies that $U_{v}$ necessarily has to have zero mean over $\Omega$.  This is not an issue for the models \eqref{ZeroExcess:equal:densities} and \eqref{Scaled:ZeroExcess}, but for non-zero $U_{v}$, analysis of the corresponding Cahn--Hilliard--Darcy model has utilized this compatibility condition, see \cite{GLDarcy,JWZ}.  In the case that the source term $U_{v}$ is a function of the variables $\varphi$ and $\sigma$, the condition \eqref{compatibility} may not hold in general, and thus alternative boundary conditions for the velocity have to be considered, see for example \cite{GLRome} for the analysis of a Cahn--Hilliard--Darcy model with a source term $U_{v}$ depending on $\varphi$ and $\sigma$, along with Dirichlet and Robin boundary conditions prescribed for the pressure.

\subsubsection{Remarks on the source terms}
To obtain useful a priori estimates from the energy identity \eqref{Energy:Id}, which are essential for the proof of existence of solutions, it suffices to control the source terms (neglecting the boundary term for the moment)
\begin{equation}\label{Source:Term}
\begin{aligned}
& \int_{\Omega} U_{v} \left ( \varphi \mu + \sigma N_{\sigma} + \frac{\overline{\rho}_{2} + \overline{\rho}_{1}}{2} \frac{\abs{\bm{v}}^{2}}{2} - p - A \Psi(\varphi) - \frac{B}{2} \abs{\nabla \varphi}^{2} - N(\varphi, \sigma) \right ) \dx \\
& \quad - \int_{\Omega} \left ( \mu - \frac{\overline{\rho}_{2} - \overline{\rho}_{1}}{2} \frac{\abs{\bm{v}}^{2}}{2} \right ) U_{\varphi} + N_{\sigma} S \dx
\end{aligned}
\end{equation}
with the ``good" part
\begin{align*}
& \frac{\dd}{\dt} \int_{\Omega} \frac{\rho}{2} \abs{\bm{v}}^{2} + A \Psi(\varphi) + \frac{B}{2} \abs{\nabla \varphi}^{2} + N(\varphi, \sigma) \dx \\
& \quad + \int_{\Omega} m(\varphi) \abs{\nabla \mu}^{2} + n(\varphi) \abs{\nabla N_{\sigma}}^{2} + 2 \eta(\varphi) \abs{\der \bm{v}}^{2} \dx.
\end{align*}
The fundamental difference between the source terms encountered in the prior work \cite{CGH,FGR,GLDirichlet,GLDarcy,GLNeumann,JWZ} for Cahn--Hilliard/Cahn--Hilliard--Darcy systems and \eqref{Source:Term} is the appearance of $\abs{\bm{v}}^{2}$ multiplied with $U_{v}$ and $U_{\varphi}$.  Suitable assumptions on the integrability of $U_{v}$ and $U_{\varphi}$ are required to control their products with the kinetic energy.  Moreover, estimating \eqref{Source:Term} is further complicated by the presence of the term $U_{v} p$ involving the pressure.  In contrast to Cahn--Hilliard--Darcy models, for which an estimate of the $L^{2}$-norm of the pressure can be obtained by studying the Darcy's equation as a second order elliptic equation for the pressure (see \cite{GLDarcy,GLRome,JWZ}), in the present setting involving the Navier--Stokes equations, such estimates seem not available at first glance.  However, we mention the work of Abels \cite{AbelsPressure} for some results in this direction.

We also point out that the pressure reformulations in Section \ref{sec:Reformulation} do not entirely remove the above technical issues, as from \eqref{alt:cv:q}, \eqref{alt:cv:r}, \eqref{alt:cv:s}, the source term $U_{v}$ will always be multiplied with a nonlinear term containing $\abs{\bm{v}}^{2}$ and the pressure.

\subsection{Sharp interface limit}\label{sec:SIM}
Let us consider the general model \eqref{NSCH} with the following choices
\begin{align}
A = \frac{\beta}{\eps}, \quad B = \beta \eps, \quad N(\varphi, \sigma) = \frac{1}{2} \abs{\sigma}^{2} + \chi \sigma (1-\varphi),
\end{align}
where $\beta > 0$ denotes the surface tension, $\eps > 0$ is related to the thickness of the diffuse interface, $\chi \geq 0$ is the chemotaxis parameter.  Furthermore, we allow the source terms $S = S(\varphi, \mu, \sigma)$, $\Gamma_{1} = \Gamma_{1}(\varphi, \mu, \sigma)$, $\Gamma_{2} = \Gamma_{2}(\varphi, \mu, \sigma)$ to depend on $\varphi, \mu, \sigma$ but not on any derivatives.  The motivation for considering a dependence on $\mu$ for the source terms comes from the work of \cite{Hawkins,Kampmann}.  We take a non-degenerate mobility $0 < n_{0} \leq n(s)$ for all $s \in \R$.

From the formally matched asymptotic analysis performed in \cite{AGG}, different sharp interface models in the limit $\eps \to 0$ can be derived from the Navier--Stokes--Cahn--Hilliard system for different choices of the mobility function $m(\varphi)$.  To compare with the results in \cite{GLSS} we confine ourselves to the constant mobility $m(\varphi) = 1$.  Furthermore, we rescale the potential $\Psi$ such that
\begin{align*}
\int_{-1}^{1} \sqrt{2 \Psi(s)} \ds = 1,
\end{align*}
so that the constant $\gamma$ in \cite[(3.17)]{GLSS} becomes $1$.

Let $\Omega_{1}(t)$ denote the region for fluid 1 and $\Omega_{2}(t)$ denote the region for fluid 2, which are separated by a hypersurface $\Sigma(t)$.  We denote the variables defined over $\Omega_{1}$ with the subscript $1$, and similarly for variables defined over $\Omega_{2}$ with the subscript $2$.  For the sake of convenience, we use the following notations
\begin{align*}
\rho = \begin{cases}
\overline{\rho}_{1} & \text{ in } \Omega_{1}, \\
\overline{\rho}_{2} & \text{ in } \Omega_{2},
\end{cases} \quad
\eta = \begin{cases}
\eta(-1) & \text{ in } \Omega_{1}, \\
\eta(+1) & \text{ in } \Omega_{2},
\end{cases} \quad
n =  \begin{cases}
n(-1) & \text{ in } \Omega_{1}, \\
n(+1) & \text{ in } \Omega_{2},
\end{cases} \\
\Gamma_{i} =  \begin{cases}
\Gamma_{i}(-1, \mu_{1}, \sigma_{1}) & \text{ in } \Omega_{1}, \\
\Gamma_{i}(1, \mu_{2}, \sigma_{2}) & \text{ in } \Omega_{2},
\end{cases} \quad S =  \begin{cases}
S(-1, \mu_{1}, \sigma_{1}) & \text{ in } \Omega_{1}, \\
S(1, \mu_{2}, \sigma_{2}) & \text{ in } \Omega_{2},
\end{cases}
\end{align*}
when writing equations over the union $\Omega_{1} \cup \Omega_{2}$.
Furthermore, let $\bm{\nu}$ and $\velo$ denote the unit normal (pointing into $\Omega_{2}$) and the normal velocity of $\Sigma(t)$, respectively.  We define the jump of a quantity $f$ along $\Sigma$ as
\begin{align*}
\jump{f}(x) := \lim_{\delta \searrow 0} \left ( f_{2}( \bm{x} + \delta \bm{\nu}(\bm{x})) - f_{1}(\bm{x} - \delta \bm{\nu}(\bm{x})) \right )
\end{align*}
for a point $\bm{x} \in \Sigma$ such that $\bm{x} + \delta \bm{\nu}(\bm{x}) \in \Omega_{2}$ and $\bm{x} - \delta \bm{\nu}(\bm{x}) \in \Omega_{1}$.  Then, following the procedure outlined in \cite[\S 4, Case I]{AGG} and \cite[\S 3]{GLSS}, while considering a similar treatment for the pressure variable as in \cite[\S 3]{DGL}, i.e., the pressure variable has a term scaling with $\frac{1}{\eps}$ in its asymptotic expansion, we can obtain the following sharp interface model:
\begin{subequations}\label{SIM}
\begin{alignat}{3}
\div \bm{v} & = \overline{\rho}_{1} \Gamma_{1} + \overline{\rho}_{1} \Gamma_{2} && \text{ in } \Omega_{1} \cup \Omega_{2}, \label{SIM:div}\\
\rho (\pd_{t} \bm{v} + \div (\bm{v} \otimes \bm{v})) -  \tfrac{\overline{\rho}_{2}-\overline{\rho}_{1}}{2} \div (\bm{v} \otimes \nabla \mu ) & = - \nabla p + 2 \eta \, \div (\der \bm{v})  && \text{ in } \Omega_{1} \cup \Omega_{2}, \label{SIM:NS} \\
\pd_{t} \sigma + \div (\sigma \bm{v}) & = n \Laplace \sigma + S && \text{ in } \Omega_{1} \cup \Omega_{2}, \label{SIM:sigma} \\
-\Laplace \mu_{2} & = - 2 \overline{\rho}_{1}^{-1} \Gamma_{1} && \text{ in } \Omega_{2}, \label{SIM:muT} \\
-\Laplace \mu_{1} & = 2 \overline{\rho}_{2}^{-1} \Gamma_{2} && \text{ in } \Omega_{1}, \label{SIM:muH}
\end{alignat}
\end{subequations}
along with the free boundary conditions:
\begin{subequations}\label{SIM:FB}
\begin{alignat}{3}
\jump{\mu} = 0, \quad \jump{\sigma} = 2 \chi & \text{\quad on } \Sigma, \label{SIM:FB:mu} \\
2 \mu = \beta \kappa - \frac{1}{2} \jump{\abs{\sigma}^{2}}, \quad 2 \chi (\bm{v} \cdot \bm{\nu} - \velo) = \jump{n \pdnu \sigma} & \text{\quad on } \Sigma, \label{SIM:FB:sigma} \\
\jump{\bm{v}} = \bm{0}, \quad \bm{v} \cdot \bm{\nu} -\velo = \tfrac{1}{2} \jump{\pdnu \mu}, \quad \jump{p \Id - 2 \eta \der \bm{v}} \bm{\nu} = \kappa \beta \bm{\nu} & \text{\quad on } \Sigma, \label{SIM:FB:v}
\end{alignat}
\end{subequations}
where $\kappa$ is the mean curvature of the  hypersurface $\Sigma$.  We remark that the equations \eqref{SIM:div}, \eqref{SIM:sigma}--\eqref{SIM:FB:sigma} also appear in the sharp interface model \cite[(3.49)--(3.50)]{GLSS}, while \eqref{SIM:NS}, \eqref{SIM:FB:v} are present in the sharp interface limit of \cite[\S 4, Case I]{AGG}.

\begin{remark}
We point out that it is also possible to consider mobilities scaling with $\eps$, i.e., $m_{1}(\varphi) = \eps$, or a two-sided degenerate mobility $m_{2}(\varphi) = (1-\varphi^{2})_{+} := \max(1-\varphi^{2},0)$.  In these two cases, the analysis of \cite[Cases II and IV]{AGG} will yield a sharp interface model consisting of \eqref{SIM:div}, \eqref{SIM:sigma} with
\begin{align*}
\rho (\pd_{t} \bm{v} + \div (\bm{v} \otimes \bm{v})) = - \nabla p + 2 \eta \, \div (\der \bm{v}) & \text{\quad in } \Omega_{1} \cup \Omega_{2},  \\
2 \overline{\rho}_{1}^{-1} \Gamma_{1}(1, \mu, \sigma_{2}) = 0 & \text{\quad in } \Omega_{2}, \\
2 \overline{\rho}_{2}^{-1} \Gamma_{2}(-1, \mu, \sigma_{1}) = 0 & \text{\quad in } \Omega_{1}, \\
2 \mu = \beta \kappa - \frac{1}{2} \jump{\abs{\sigma}^{2}}, \quad \jump{n \pdnu \sigma} = 0, \quad \jump{\sigma} = 2 \chi & \text{\quad on } \Sigma, \\
\jump{\bm{v}}  = \bm{0}, \quad  \jump{p \Id - 2 \eta \der \bm{v}} \bm{\nu} = \kappa \beta \bm{\nu}, \quad \velo = \bm{v} \cdot \bm{\nu} & \text{\quad on } \Sigma.
\end{align*}
In particular, the choices of $m_{1}$ and $m_{2}$ lead to the relation $\velo = \bm{v} \cdot \bm{\nu}$, which means that the interface is passively transported by the fluid velocity.  Note that there are no differential operators acting on the variable $\mu$, and the source terms $\Gamma_{1}$, $\Gamma_{2}$ must be chosen such that $\Gamma_{1} = 0$ in $\Omega_{2}$ and $\Gamma_{2} = 0$ in $\Omega_{1}$ in order for the above sharp interface model to be well-defined.
\end{remark}

\section{Analysis of the simplified model with zero excess of total mass and equal densities}\label{sec:Analysis}

\subsection{Problem setting}

For the rest of the paper, we analyze the model \eqref{ZeroExcess:equal:densities} with zero excess of total mass and equal densities, and consider the following typical form of the chemical free energy $N$:
\begin{align}
N(\varphi, \sigma) = \frac{1}{2} \abs{\sigma}^{2} + \chi \sigma (1-\varphi)\quad \text{ for } \chi \geq 0,
\end{align}
which has also been considered in \cite{GLDirichlet,GLDarcy,GLRome,GLNeumann}.  Then we have
\begin{align}
N_{\varphi} = - \chi \sigma, \quad N_{\sigma} = \sigma + \chi (1-\varphi).
\end{align}
For any given $T \in (0,+\infty)$, we set $\rho_{*} = 1$ and absorb the prefactor $2$ into the source term $\Gamma$, leading to the system
\begin{subequations}\label{ZEED}
\begin{alignat}{3}
 \pd_{t} \bm{v} + (\bm{v} \cdot \nabla ) \bm{v}  & = - \nabla q + \div (2 \eta(\varphi) \der \bm{v}) + (\mu + \chi \sigma) \nabla \varphi && \text{\quad in } \Omega \times (0,T), \label{ZEED:mom} \\
 \div \bm{v} & = 0 && \text{\quad in } \Omega \times (0,T), \\
\pd_{t} \varphi + \bm{v} \cdot \nabla \varphi & = \div (m(\varphi) \nabla \mu) +  \Gamma && \text{\quad in } \Omega \times (0,T), \label{ZEED:varphi} \\
\mu & = A \Psi'(\varphi) - B \Laplace \varphi - \chi \sigma && \text{\quad in } \Omega \times (0,T), \label{ZEED:mu} \\
\pd_{t} \sigma + \bm{v} \cdot \nabla \sigma & = \div (n(\varphi) \nabla (\sigma + \chi (1-\varphi))) + S && \text{\quad in } \Omega \times (0,T), \label{ZEED:sigma}
\end{alignat}
\end{subequations}
where we have also employed the reformulation \eqref{alt:mom:q} for the momentum equation with the modified pressure variable $q$.
We now prescribe the following initial-boundary conditions
\begin{subequations}\label{ZEED:bdy}
\begin{alignat}{3}
 \bm{v}(0) = \bm{v}_{0},\quad \varphi(0) = \varphi_{0}, \quad \sigma(0) = \sigma_{0} \quad & \text{ in } \Omega, \\
\bm{v} = \bm{0},\quad  m(\varphi)\pdnu \mu = \pdnu \varphi= n(\varphi) \pdnu N_{\sigma}  = 0 \quad & \text{ on } \pd \Omega \times (0,T),
\end{alignat}
\end{subequations}
so that the energy identity \eqref{Energy:Id} reduces to
\begin{equation}\label{Energy:ZEED}
\begin{aligned}
& \frac{\dd}{\dt} \int_{\Omega} \left[\frac{1}{2} \abs{\bm{v}}^{2} + A \Psi(\varphi) + \frac{B}{2} \abs{\nabla \varphi}^{2} + \frac{1}{2} \abs{\sigma}^{2} + \chi \sigma (1-\varphi)\right] \dx \\
& \quad + \int_{\Omega} \left[m(\varphi) \abs{\nabla \mu}^{2} + n(\varphi) \abs{\nabla (\sigma + \chi (1- \varphi))}^{2} + 2 \eta(\varphi) \abs{\der \bm{v}}^{2}\right] \dx  \\
& \quad - \int_{\Omega} \left[ (\sigma + \chi (1-\varphi)) S + \mu \Gamma\right] \dx  = 0.
\end{aligned}
\end{equation}

\subsection{Preliminaries}
Assume that $\Omega \subset \R^{d}$, $d = 2,3$, is a bounded domain with smooth boundary $\partial\Omega$.
We denote by $Q = \Omega \times (0,T)$ the space-time cylinder with $\Sigma := \pd \Omega \times (0,T)$.  For the standard Lebesgue and Sobolev spaces, we use the notations $L^{p} := L^{p}(\Omega)$ and $W^{k,p} := W^{k,p}(\Omega)$ for any $p \in [1,+\infty]$, $k > 0$ equipped with the norms $\norm{\cdot}_{L^{p}}$ and $\norm{\cdot}_{W^{k,p}}$.  In the case $p = 2$ we use $H^{k} := W^{k,2}$ and the norm $\norm{\cdot}_{H^{k}}$.  The $L^{2}$-scalar product between two functions $f$ and $g$ is denoted by $(f,g)$.  The dual space of a Banach space $X$ is denoted by $X'$, and  the duality pairing between $X$ and its dual will be denoted by $\inner{\cdot}{\cdot}_{X}$.

We express $\R^{d}$-valued functions and function spaces consisting of vector-valued/tensor-valued functions in boldface.  The matrix of second derivatives of a scalar function $f$ is denoted by $\nabla^{2} f$, and similarly the tensor of second derivatives of a vector $\bm{f}$ is denoted by $\nabla^{2} \bm{f}$.  For Bochner spaces, we will often use the isometric isomorphism between $L^{p}(0,T;L^{p}(\Omega))$ and $L^{p}(Q)$ for $1 \leq p < \infty$.  Similarly, we reuse the notation $L^{p}(Q)$ for $L^{p}(0,T;\bm{L}^{p}(\Omega))$ for $1 \leq p < \infty$.
Since the pressure is determined up to a time-dependent constant, we ask that the pressure $q$ belongs to the space $L^{2}_{0} := \{ f \in L^{2}(\Omega) : \mean{f} = 0 \}$, where $\mean{f} := \frac{1}{\abs{\Omega}} \int_{\Omega} f \dx$ denotes the mean of $f$.  Furthermore, we set
\begin{align*}
H^{2}_{N} := \{ f \in H^{2} : \pdnu f = 0 \text{ on } \pd \Omega \}.
\end{align*}

We now introduce the classical function spaces for the Navier--Stokes equations \cite{Sohr}.  For a vector-valued/tensor-valued Banach space $\bm{X}$, we define $\bm{X}_{0,\sigma}$ as the completion of $C^{\infty}_{0,\sigma} := \{ \bm{f} \in (C^{\infty}_{0}(\Omega))^{d} : \div \bm{f} = 0 \}$ with respect to the $\bm{X}$-norm.  In the case $\bm{X} = \bm{L}^{2}$, we use the notation $\bm{L}^{2}_{\sigma} := \bm{L}^{2}_{0,\sigma}$.  The space $\bm{H}^{1}_{0,\sigma}$ is endowed with the scalar product
\begin{align*}
(\bm{u}, \bm{v})_{\bm{H}^{1}_{0,\sigma}} := (\nabla \bm{u}, \nabla \bm{v}) \quad \forall \bm{u}, \bm{v} \in\bm{H}^{1}_{0,\sigma},
\end{align*}
and we denote its topological dual as $\bm{H}^{-1}$.  It is well-known that $\bm{L}^{2}$ can be decomposed into the sum $\bm{L}^{2}_{\sigma} \oplus \bm{G}(\Omega)$, where $\bm{G}(\Omega) := \{ \bm{f} \in \bm{L}^{2} : \exists z \in L^{2} \text{ with } \bm{f} = \nabla z \}$, i.e., the space $\bm{G}(\Omega)$ is the orthogonal complement of $\bm{L}^{2}_{\sigma}$.  Then, for every $\bm{f} \in \bm{L}^{2}$, we have the unique decomposition (up to an additive constant for $z$)
\begin{align*}
\bm{f} = \bm{f}_{0} + \nabla z, \text{ where }  \bm{f}_{0} \in \bm{L}^{2}_{\sigma}, \quad \nabla z \in \bm{G}(\Omega),
\end{align*}
and as a consequence, we obtain a bounded linear operator $P_{H} : \bm{L}^{2} \to \bm{L}^{2}_{\sigma}$ defined by $P_{H}(\bm{f}) = \bm{f}_{0}$, which is more commonly known as the Leray projection onto the space of divergence-free functions.  We recall the Stokes operator $\Stokes : D(\Stokes)\to \bm{L}^{2}_{\sigma}$ defined as
\begin{align*}
(\Stokes \bm{u}, \bm{\zeta}) = (\nabla \bm{u}, \nabla \bm{\zeta}) \quad \forall\, \bm{\zeta} \in \bm{H}^{1}_{0,\sigma}
\end{align*}
with domain $D(\Stokes) = \bm{H}^{2} \cap \bm{H}^{1}_{0,\sigma}$. Then the following estimates hold (see \cite[Lemma 3.4]{Lorca} and \cite{Sohr}):
\begin{lemma}\label{stoo}
For any $\bm{u} \in D(\Stokes)$, consider the Helmholtz decomposition $\Stokes \bm{u} = -\Laplace \bm{u} + \nabla \pi$ where the pressure-like function $\pi$ is taken such that $\int_{\Omega} \pi \dx = 0$.
Then for any $\nu > 0$, there exists a positive constant $C_{\nu}$ independent of $\bm{u}$ such that
\begin{align}
\norm{\pi}_{L^{2} \backslash {\mathbb{R}}} \leq \nu\norm{\Stokes \bm{u}}_{\bm{L}^2}+C_\nu\norm{\nabla\bm{u}}_{\bm{L}^{2}}. \label{Stokes I}
\end{align}
Moreover, there exists a positive constant $c=c(d, \Omega)$ such that
\begin{align}
\norm{\bm{u}}_{\bm{H}^{2}}+\norm{\pi}_{H^{1} \backslash {\mathbb{R}}} \leq c\norm{\Stokes\bm{u}}_{\bm{L}^{2}}. \label{Stokes II}
 \end{align}
\end{lemma}
\noindent Furthermore, we state some inequalities that will be useful in the subsequent analysis:
\begin{itemize}
\item \emph{Poincar\'{e}'s inequality}: There exists a positive constant $C$ depending only on $\Omega$ such that, for all $f \in W^{1,p}$, $p \in [1,\infty]$,
\begin{align}
\bignorm{f - \overline{f}}_{L^{p}} & \leq C \norm{\nabla f}_{L^{p}}. \label{regular:Poincare}
\end{align}
\item \emph{The Gagliardo--Nirenberg interpolation inequality} in dimension $d$:  Let $\Omega$ be a bounded domain with Lipschitz boundary, and $f \in W^{m,r} \cap L^{q}$, $1 \leq q,r \leq \infty$.  For any integer $j$, $0 \leq j < m$, suppose there is $\alpha \in \R$ such that
\begin{align}\label{GN:constraints}
\frac{1}{p} = \frac{j}{d} + \left ( \frac{1}{r} - \frac{m}{d} \right ) \alpha + \frac{1-\alpha}{q}, \quad \frac{j}{m} \leq \alpha \leq 1.
\end{align}
If $r \in (1,\infty)$ and $m - j - \frac{d}{r}$ is a non-negative integer, we additionally assume $\alpha \neq 1$.  Under these assumptions, there exists a positive constant $C$ depending only on $\Omega$, $m$, $j$, $q$, $r$, and $\alpha$ such that
\begin{align}
\label{GagNirenIneq}
\norm{D^{j} f}_{L^{p}} \leq C \norm{f}_{W^{m,r}}^{\alpha} \norm{f}_{L^{q}}^{1-\alpha},
\end{align}
where $D^{j} f$ denotes the $j$-th weak partial derivatives of $f$.  An alternate variant of the Gagliardo--Nirenbeg inequality we will use is
\begin{align*}
\norm{D^{j}f}_{L^{p}} \leq C \norm{D^{m}f}_{L^{r}}^{\alpha} \norm{f}_{L^{q}}^{1-\alpha} + C \norm{f}_{L^{s}},
\end{align*}
where $p$, $j$, $m$, $r$, $q$ and $\alpha$ satisfy \eqref{GN:constraints}, and $s > 0$ is arbitrary.
\item \emph{Korn's inequality}:  For any $\bm{u} \in \bm{H}^{1}_{0}$, there exists a positive constant $C$, depending only on $\Omega$ such that
\begin{align}\label{Korn}
\norm{\nabla \bm{u}}_{\bm{L}^{2}} \leq C \norm{\der \bm{u}}_{\bm{L}^{2}}, \text{ where } \der \bm{u} := \frac{1}{2} (\nabla \bm{u} + (\nabla \bm{u})^{\top}).
\end{align}
\item \emph{The Br\'{e}zis--Gallouet interpolation inequality} for $d=2$ (cf. \cite{BG,Engler}):  There exists a positive constant $C$ depending only on $\Omega$ such that
\begin{align}\label{Bre}
 \norm{g}_{L^{\infty}} \leq C \norm{g}_{H^{1}}\sqrt{\ln (1+ \norm{g}_{H^{2}})} + C \norm{g}_{H^{1}}, \quad \forall \, g \in H^{2}(\Omega).
 \end{align}
\item Elliptic estimates: If $f \in H^{2}(\Omega)$ satisfies $\pdnu f = 0$ on $\pd \Omega$, then there exists a positive constant $C$ depending only on $\Omega$ such that
\begin{align}\label{H2EllEst}
\norm{f}_{H^{2}} \leq C \left ( \norm{\Laplace f}_{L^{2}} + \norm{f}_{L^{2}} \right ).
\end{align}
If $f \in H^{4}(\Omega)$ satisfies $\pdnu f = \pdnu (\Laplace f) = 0$ on $\pd \Omega$, then there exists a positive constant $C$ depending only on $\Omega$ such that
\begin{align}\label{H4EllEst}
\norm{f}_{H^{4}} \leq C \left ( \norm{\Laplace^{2} f}_{L^{2}} + \norm{f}_{L^{2}} \right ).
\end{align}
\end{itemize}

\subsection{Main results}\label{sec:ZEED:weak}

Our first result concerns the existence of global weak solutions to the system \eqref{ZEED}--\eqref{ZEED:bdy} in both two and three dimensions.

\begin{thm}[Global weak solutions in 2D and 3D]\label{thm:ZEED:weaksoln}
We assume that
\begin{enumerate}
\item $m, n, \eta \in C^{0}(\R)$ and satisfy
\begin{align*}
m_{0} \leq m(s) \leq m_{1}, \quad n_{0} \leq n(s) \leq n_{1}, \quad \eta_{0} \leq \eta(s) \leq \eta_{1} \quad \forall s \in \R,
\end{align*}
where $m_0$, $m_1$, $n_0$, $n_1$, $\eta_0$, $\eta_1$ are given positive constants.
\item The external source terms $S \in L^{2}(Q)$ and $\Gamma \in L^{2}(0,T;L^{\infty}(\Omega))$ are prescribed functions.
\item The potential $\Psi \in C^{2}(\R)$ is non-negative and satisfies
\begin{equation}\label{Psi:assump}
\begin{aligned}
\abs{\Psi''(s)} \leq C_{0} (1 + \abs{s}^{r}), \quad  \abs{\Psi'(s)} \leq C_{1} \Psi(s) + C_{2},  \quad \Psi(s) \geq C_{3} \abs{s}^{2} - C_{4},
\end{aligned}
\end{equation}
for all $s \in \R$ with positive constants $C_{0}, C_{1}, C_{2}, C_{3}$ and $C_{4}$ that are independent of $s$. If $d=3$,  $ r \in [0,4)$; if $d=2$, $r\in [0,+\infty)$.
\item The coefficients $A > 0$, $B > 0$ and $\chi \geq 0$ are given constants that satisfy
\begin{align}\label{assump:A}
A > \frac{2 \chi^{2}}{C_{3}}.
\end{align}
\end{enumerate}
Let $T\in(0,+\infty)$ be an arbitrary but fixed terminal time.
Then for any initial data $\varphi_{0} \in H^{1}(\Omega)$, $\sigma_{0} \in L^{2}(\Omega)$ and $\bm{v}_{0} \in \bm{L}^{2}(\Omega)$, there exists at least one global weak solution  $(\varphi, \mu, \sigma, \bm{v})$ to problem \eqref{ZEED}--\eqref{ZEED:bdy} on $[0,T]$ such that
\begin{align*}
\varphi & \in L^{\infty}(0,T;H^{1}) \cap L^{2}(0,T;H^{3}) \cap H^{1}(0,T;(H^{1})'), \quad \mu \in L^{2}(0,T;H^{1}), \\
\sigma & \in L^{\infty}(0,T;L^{2}) \cap L^{2}(0,T;H^{1}) \cap  W^{1,y}(0,T;(H^{1})'), \\
\bm{v} & \in L^{\infty}(0,T;\bm{L}^{2}_{\sigma}) \cap L^{2}(0,T;\bm{H}^{1}_{0,\sigma}) \cap  W^{1,y}(0,T;\bm{H}^{-1}),
\end{align*}
where $y = 2$ if $d=2$, $y = \frac{4}{3}$ if $d=3$. Moreover, the quadruple $(\varphi, \mu, \sigma, \bm{v})$ satisfies for a.e. $t \in (0,T)$
\begin{subequations}\label{ZEED:weakform}
\begin{alignat}{3}
0 & = \inner{\pd_{t}\varphi}{\zeta}_{H^{1}} + (\bm{v} \cdot \nabla \varphi, \zeta) + (m(\varphi) \nabla \mu, \nabla \zeta) - (\Gamma, \zeta), \\
0 & = (\mu - A \Psi'(\varphi) + \chi \sigma, \zeta) - (B\nabla \varphi, \nabla \zeta), \\
0 & = \inner{\pd_{t} \sigma + \bm{v} \cdot \nabla \sigma}{\zeta}_{H^{1}} + (n(\varphi) \nabla(\sigma + \chi (1-\varphi)), \nabla \zeta) - (S, \zeta), \\
0 & = \inner{\pd_{t} \bm{v} + (\bm{v} \cdot \nabla) \bm{v}}{\bm{\xi}}_{\bm{H}^{1}} + (2 \eta(\varphi) \der \bm{v}, \der \bm{\xi}) - ((\mu + \chi \sigma) \nabla \varphi, \bm{\xi}), \label{ZEED:weakform:velo}
\end{alignat}
\end{subequations}
 and
\begin{align*}
\varphi(0) = \varphi_{0}, \quad \inner{\sigma(0)}{\zeta}_{H^{1}} = \inner{\sigma_{0}}{\zeta}_{H^{1}}, \quad \inner{\bm{v}(0)}{\bm{\xi}}_{\bm{H}^{1}} = \inner{\bm{v}_{0}}{\bm{\xi}}_{\bm{H}^{1}},
\end{align*}
for all $\zeta \in H^{1}(\Omega)$ and $\bm{\xi} \in \bm{H}^{1}_{0,\sigma}$.
\end{thm}

\begin{remark}
\
$\mathrm{(1)}.$ The initial value for $\varphi$ is attained precisely due to the (compact) embedding $L^{\infty}(0,T;H^{1}) \cap H^{1}(0,T;(H^{1})') \subset \subset C^{0}([0,T];L^{2})$.  Meanwhile, for $\sigma$ and $\bm{v}$, in two dimensions we use the continuous embedding
$L^{2}(0,T;H^{1}) \cap H^{1}(0,T;(H^{1})') \subset  C^{0}([0,T];L^{2})$
to deduce that $\sigma(0) = \sigma_{0}$ and $\bm{v}(0) = \bm{v}_{0}$.
In three dimensions, we use the continuous embedding $W^{1,\frac{4}{3}}(0,T;(H^{1})') \subset C^{0}([0,T];(H^{1})')$
to deduce that the initial conditions for $\sigma$ and $\bm{v}$ are attained.

$\mathrm{(2)}.$ To recover the pressure, we argue as follows.  The distribution $\bm{F}_{\bm{v}}$ defined as
\begin{align*}
\bm{F}_{\bm{v}} := \pd_{t} \bm{v} + (\bm{v} \cdot \nabla) \bm{v} - \div (2 \eta(\varphi) \der \bm{v}) - (\mu + \chi \sigma) \nabla \varphi
\end{align*}
belongs to $L^{y}(0,T;\bm{H}^{-1})$ and vanishes on the subspace $C^{\infty}_{0}([0,T];C^{\infty}_{0,\sigma})$ by \eqref{ZEED:weakform:velo}.  Applying \cite[Lemma IV.1.4.1]{Sohr} allows us to deduce that there exists a function $q \in L^{y}(0,T;L^{2}_{0})$ such that
$\bm{F}_{\bm{v}} = -\nabla q$ holds as an equality in the sense of distributions for a.e. $t \in (0,T)$.
\end{remark}\medskip

Our next result yields the existence of global strong solutions to problem \eqref{ZEED}--\eqref{ZEED:bdy} when the spatial dimension is two.

\begin{thm}[Global strong solutions in 2D]\label{thm:2D:strong}
In addition to the assumptions of Theorem \ref{thm:ZEED:weaksoln}, we assume that
\begin{enumerate}
\item $m \in C^{2}_{b}(\R)$, $n \in C^{1}_{b}(\R)$, $\eta \in C^{1}_{b}(\R)$ with bounded derivatives.
\item The potential $\Psi \in C^{3}(\R)$ satisfies
\begin{align}
\abs{\Psi'''(s)} & \leq C_{0} (1 + \abs{s}^{r-1}) \quad \text{ for } r \in [1,+\infty),
\end{align}
instead of the first assumption in \eqref{Psi:assump}.
\item The source term $\Gamma$ belongs to $L^{2}(0,T;H^{1}(\Omega) \cap L^{\infty}(\Omega))$.
\end{enumerate}
Then for any $T>0$ and arbitrary initial data $\bm{v}_{0} \in \bm{H}^{1}(\Omega)$, $\sigma_{0} \in H^{1}(\Omega)$ and $\varphi_{0} \in H^{3}(\Omega) \cap H^{2}_{N}(\Omega)$,
there exists at least one global strong solution $(\varphi, \mu, \sigma, \bm{v})$ to problem \eqref{ZEED}--\eqref{ZEED:bdy} such that
\begin{align*}
\varphi & \in L^{\infty}(0,T;H^{3} \cap H^{2}_{N}) \cap L^{2}(0,T;H^{4}) \cap H^{1}(0,T;H^{1}) \cap C^{0}([0,T];C^{1,\delta}(\overline{\Omega})), \\
\mu & \in L^{\infty}(0,T;H^{1})\cap L^{2}(0,T;H^{3} \cap H^{2}_{N}), \\
\sigma & \in L^{\infty}(0,T;H^{1}) \cap L^{2}(0,T;H^{2}_{N}) \cap  H^{1}(0,T;L^{2}),\\
\bm{v} & \in L^{\infty}(0,T;\bm{H}^{1}_{0}) \cap L^{2}(0,T;\bm{H}^{2}) \cap H^{1}(0,T;\bm{L}^{2}_{\sigma}) , \\
q & \in L^{2}(0,T;H^{1} \cap L^{2}_{0}),
\end{align*}
for some $0 < \delta < 1$.
For a.e. $(x,t) \in Q$, the quadruple $(\varphi, \mu, \sigma, \bm{v})$ satisfies
\begin{subequations}\label{Str:Soln:ZEED}
\begin{align}
0 & = \pd_{t}\varphi + \bm{v} \cdot \nabla \varphi - \div (m(\varphi) \nabla \mu) - \Gamma, \\
0 & = \mu - A \Psi'(\varphi) + B \Laplace \varphi + \chi \sigma, \\
0 & = \pd_{t}\sigma + \bm{v} \cdot \nabla \sigma - \div (n(\varphi) \nabla (\sigma + \chi (1-\varphi))) - S, \\
0 & = \pd_{t} \bm{v} + (\bm{v} \cdot \nabla) \bm{v} - \div (2 (\eta(\varphi) \der \bm{v}) + \nabla q - (\mu + \chi \sigma) \nabla \varphi.
\end{align}
\end{subequations}
Furthermore, it holds
\begin{align*}
\varphi(0) = \varphi_{0}, \quad \sigma(0) = \sigma_{0} \quad \bm{v}(0) = \bm{v}_{0} \quad  \text{ a.e. in } \Omega.
\end{align*}
\end{thm}

\begin{remark}
$\mathrm{(1)}.$ We used the compact embedding $L^{\infty}(0,T;H^{3}) \cap H^{1}(0,T;H^{1}) \subset \subset C^{0}([0,T];C^{1,\delta}(\overline{\Omega}))$ to deduce the H\"{o}lder spatial regularity for $\varphi$.

$\mathrm{(2)}.$ The assumption $\varphi_{0} \in H^{2}_{N}(\Omega)$ implies that $\mu_{0} := A \Psi'(\varphi_{0}) - B \Laplace \varphi_{0} - \chi \sigma_{0}  \in L^{2}(\Omega)$, and if $\varphi_{0} \in H^{3}(\Omega)$, then $\mu_{0} \in H^{1}(\Omega)$.

$\mathrm{(3)}.$ We point out that one can weaken the assumptions to $m \in C^{1}_{b}(\R)$ with bounded derivatives and $\varphi_{0} \in H^{2}_{N}(\Omega)$ to derive a strong solution with less regularities:
\begin{align*}
\varphi & \in L^{\infty}(0,T;H^{3} \cap H^{2}_{N}) \cap L^{2}(0,T;H^{4}) \cap H^{1}(0,T;L^{2}), \\
\mu & \in L^{\infty}(0,T;L^{2}) \cap L^{2}(0,T;H^{2}_{N}).
\end{align*}
In essence, we neglect the fourth a priori estimate in Section \ref{sec:ZEED:2D:strong}.
\end{remark}\medskip

The last theorem gives a continuous dependence on the initial data and source terms for global strong solutions in two dimensions.

\begin{thm}[Continuous dependence in 2D]\label{thm:2D:ctsdep}
Let $\{\varphi_{i}, \mu_{i}, \sigma_{i}, \bm{v}_{i}, q_{i}\}_{i=1,2}$ denote two global strong solutions to problem \eqref{ZEED}--\eqref{ZEED:bdy} corresponding to initial data $\{\varphi_{0,i}, \sigma_{0,i}, \bm{v}_{0,i}\}_{i=1,2}$ and source terms $\{\Gamma_{i}, S_{i}\}_{i=1,2}$ obtained in Theorem \ref{thm:2D:strong}.
Besides, we assume that $m(\cdot)$, $\eta(\cdot)$ and $n(\cdot)$ are Lipschitz continuous, and $\Psi$ satisfies
\begin{align}\label{assump:Psi:ctsdep}
\abs{\Psi'(s_1) - \Psi'(s_2)} \leq C_{5} \left ( 1 + \abs{s_1}^{r} + \abs{s_2}^{r} \right ) \abs{s_1-s_2} \quad \forall\, s_1, s_2 \in \R,
\end{align}
for some $r > 0$ and constant $C_{5} > 0$. Then there exists a positive constant $C$, depending on $T$, $\Omega$, $\chi$, $\eta_{0}$, $n_{0}$, $m_{0}$, $\norm{\eta'}_{L^{\infty}(\R)}$, $\norm{n'}_{L^{\infty}(\R)}$, $\norm{m'}_{L^{\infty}(\R)}$,  $A$, $B$, $C_{i}$ ($i=1,...,5$), $\norm{\nabla \sigma_{i}}_{L^{2}(0,T;\bm{L}^{4})}$, $\norm{\varphi_{i}}_{L^{2r}(0,T;L^{\infty})}$, $\norm{\mu_{i} + \chi \sigma_{i}}_{L^{2}(0,T;L^{\infty})}$, $\norm{\nabla \varphi_{i}}_{L^{\infty}(0,T;\bm{L}^{\infty})}$, $\norm{\mu_{i}}_{L^{2}(0,T;H^{3})}$, $\norm{\nabla (\sigma_{i} - \chi \varphi_{i})}_{L^{4}(Q)}$ and $\norm{\der \bm{v}_{i}}_{L^{4}(Q)}$ such that
\begin{equation}\label{2D:CtsDep:Result}
\begin{aligned}
\sup_{t \in (0,T]} & \left ( \norm{(\varphi_{1} - \varphi_{2})(t)}_{L^{2}}^{2} + \norm{(\bm{v}_{1} - \bm{v}_{2})(t)}_{\bm{L}^{2}}^{2} + \norm{(\sigma_{1} - \sigma_{2})(t)}_{L^{2}}^{2} \right ) +  \norm{\mu_{1} - \mu_{2}}_{L^{2}(Q)}^{2} \\
& \qquad  + \norm{\nabla (\sigma_{1} - \sigma_{2})}_{L^{2}(Q)}^{2} + \norm{\nabla (\bm{v}_{1} - \bm{v}_{2})}_{L^{2}(Q)}^{2} + \norm{\varphi_{1} - \varphi_{2}}_{L^{2}(0,T;H^{2})}^{2} \\
& \quad \leq C \left (\norm{\sigma_{0,1} - \sigma_{0,2}}_{L^{2}}^{2} + \norm{\bm{v}_{0,1} - \bm{v}_{0,2}}_{\bm{L}^{2}}^{2} + \norm{\varphi_{0,1} - \varphi_{0,2}}_{L^{2}}^{2} \right ) \\
& \qquad + C \left ( \norm{\Gamma_{1} - \Gamma_{2}}_{L^{2}(Q)}^{2} + \norm{S_{1} - S_{2}}_{L^{2}(Q)}^{2} \right ).
\end{aligned}
\end{equation}
In particular, the global strong solution to problem \eqref{ZEED}--\eqref{ZEED:bdy} is unique under the above assumptions.
\end{thm}


\begin{remark}
It is worth mentioning that above results concerning global weak existence and strong well-posedness also hold when the mass transfer terms $\Gamma$ and $S$ are functions of $\varphi$ and $\sigma$.  For example, one may consider the typical choice
\begin{align*}
\Gamma = h(\varphi)P(\sigma), \quad S = - h(\varphi) \mathcal{C} \sigma,
\end{align*}
where $h(s)$ is a non-negative, bounded, Lipschitz function such that $h(-1) = 0$ and $h(1) = 1$, $\mathcal{C} \geq 0$ can be seen as a constant consumption rate, and $P$ is a non-negative, bounded, Lipschitz function modeling growth by consuming the chemical species.
Similar types of mass transfer terms have been studied in \cite{GLDirichlet,GLNeumann} for the Cahn--Hilliard--Darcy system modeling tumor growth.
\end{remark}

\section{Proof of Theorem \ref{thm:ZEED:weaksoln}: global weak solutions}\label{sec:weak}

Theorem \ref{thm:ZEED:weaksoln} can be proved by using a suitable Galerkin approximation scheme based on the energy identity \eqref{Energy:ZEED}.

\subsection{Galerkin approximation}
We recall that the eigenfunctions $\{w_{j}\}_{j \in \N}$ of the Neumann--Laplacian which forms an orthonormal basis of $L^{2}$ and is also a basis of $H^{2}_{N}$, see for example \cite{GLDarcy}.
Furthermore, we can choose $w_{1} = 1$. Meanwhile, by the Lax--Milgram theorem, for every $\bm{f} \in \bm{L}^{2}_{\sigma}$ there exists a unique $\bm{v} \in  \bm{H}^{1}_{0,\sigma}$ satisfying
\begin{align}\label{Weak:Stokes}
(\nabla \bm{v}, \nabla \bm{\zeta}) = (\bm{f}, \bm{\zeta}) \quad \forall\, \bm{\zeta} \in \bm{H}^{1}_{0,\sigma}.
\end{align}
Defining the operator $\AInv : \bm{L}^{2}_{\sigma} \to \bm{H}^{1}_{0,\sigma}$ by $\AInv (\bm{f}) = \bm{v}$, where $\bm{v}$ solves \eqref{Weak:Stokes} with data $\bm{f}$.  The operator $\AInv$ can also be identified as the inverse to the Stokes operator $\Stokes$.  It is clear that $\AInv$ is a linear, continuous and self-adjoint operator.  Furthermore, by the compact embedding $\bm{H}^{1}_{0,\sigma} \subset \subset \bm{L}^{2}_{\sigma}$ it holds that $\AInv$ is a compact operator, and thus by spectral theory there exists a increasing sequence $\{\lambda_{j}\}_{j \in \N}$ of positive real numbers and a family $\{\bm{y}_{j}\}_{j \in \N} \subset \bm{H}^{2} \cap \bm{H}^{1}_{0,\sigma}$ of eigenfunctions that form an orthonormal basis of $\bm{L}^{2}_{\sigma}$ and an orthogonal basis of $\bm{H}^{1}_{0,\sigma}$.

We define the finite-dimensional subspaces
\begin{align*}
W_{n} := \mathrm{span}\{w_{1}, \dots, w_{n}\}  \text{ and } \bm{Y}_{n} := \mathrm{span}\{ \bm{y}_{1}, \dots, \bm{y}_{n}\}
\end{align*}
with corresponding orthogonal projections $\P_{W_{n}}$ and $\P_{\bm{Y}_{n}}$, and look for functions of the form
\begin{align*}
\varphi_{n}(x,t) & := \sum_{i=1}^{n} \alpha_{n,i}(t) w_{i}(x),  & \mu_{n}(x,t) & := \sum_{i=1}^{n} \beta_{n,i}(t) w_{i}(x), \\
\sigma_{n}(x,t) & := \sum_{i=1}^{n} \gamma_{n,i}(t) w_{i}(x),  & \bm{v}_{n}(x,t) & := \sum_{i=1}^{n} \delta_{n,i}(t) \bm{y}_{i}(x)
\end{align*}
that solve the following approximating problem:
\begin{subequations}\label{Galerkin:Approx}
\begin{alignat}{3}
0 & = \pd_{t} \varphi_{n} + \P_{W_{n}}(\bm{v}_{n} \cdot \nabla \varphi_{n}) - \P_{W_{n}}(\div (m(\varphi_{n}) \nabla \mu_{n})) - \P_{W_{n}}(\Gamma),  \label{Galerkin:varphi}  \\
0 & = \pd_{t} \sigma_{n} + \P_{W_{n}}(\bm{v}_{n} \cdot \nabla \sigma_{n}) - \P_{W_{n}}(\div(n(\varphi_{n}) \nabla (\sigma_{n} + \chi (1-\varphi_{n})))) - \P_{W_{n}}(S), \label{Galerkin:sigma} \\
0 & = \pd_{t} \bm{v}_{n} + \P_{\bm{Y}_{n}}((\bm{v}_{n} \cdot \nabla) \bm{v}_{n}) - \P_{\bm{Y}_{n}}(\div (2 \eta(\varphi_{n}) \der \bm{v}_{n})) - \P_{\bm{Y}_{n}}((\mu_{n} + \chi \sigma_{n}) \nabla \varphi_{n}), \label{Galerkin:velo} \\
0 & = \mu_{n} - A \P_{W_{n}}(\Psi'(\varphi_{n})) + B \Laplace \varphi_{n} + \chi \sigma_{n}, \label{Galerkin:mu}
\end{alignat}
\end{subequations}
subject to the initial data
\begin{align*}
\varphi_{n}(0) = \P_{W_{n}}(\varphi_{0}), \quad \sigma_{n}(0) = \P_{W_{n}}(\sigma_{0}), \quad \bm{v}_{n}(0) = \P_{\bm{Y}_{n}}(\bm{v}_{0}).
\end{align*}
It is easy to see that \eqref{Galerkin:Approx} is equivalent to
\begin{subequations}
\begin{alignat}{3}
0 & = (\pd_{t} \varphi_{n}, w_{i}) + (\bm{v}_{n} \cdot \nabla \varphi_{n}, w_{i}) + (m(\varphi_{n}) \nabla \mu_{n}, \nabla w_{i}) - (\Gamma, w_{i}), \label{Galerkin:varphi:var} \\
0 & = (\pd_{t} \sigma_{n}, w_{i}) + (\bm{v}_{n} \cdot \nabla \sigma_{n}, w_{i}) + (n(\varphi_{n}) \nabla (\sigma_{n} + \chi(1-\varphi_{n})), \nabla w_{i}) - (S, w_{i}), \label{Galerkin:sigma:var} \\
0 & = (\pd_{t} \bm{v}_{n}, \bm{y}_{i}) + ((\bm{v}_{n} \cdot \nabla)\bm{v}_{n}, \bm{y}_{i}) + (2 \eta(\varphi_{n}) \der \bm{v}_{n}, \der \bm{y}_{i}) - ((\mu_{n} + \chi \sigma_{n}) \nabla \varphi_{n}, \bm{y}_{i}), \label{Galerkin:velo:var} \\
0 & = (\mu_{n}, w_{i}) - (A \Psi'(\varphi_{n}), w_{i}) - (B \nabla \varphi_{n}, \nabla w_{i}) + (\chi \sigma_{n}, w_{i}), \label{Galerkin:mu:var}
\end{alignat}
\end{subequations}
for all $i \in \N$.  The approximating problem \eqref{Galerkin:Approx} is equivalent to solving a system of nonlinear ODEs in the $4n$ unknowns $\{\alpha_{n,i}, \beta_{n,i}, \gamma_{n,i}, \delta_{n,i}\}_{i = 1}^{n}$.  The continuity of $\Psi'(\cdot)$, $m(\cdot)$, $n(\cdot)$, $\eta(\cdot)$ and the Cauchy--Peano theorem allow us to deduce the existence of a local solution $(\bm{\alpha}, \bm{\beta}, \bm{\gamma}, \bm{\delta}) \in C^{1}([0,t_{n}), \R^{4n})$ for some $t_{n} \in (0,T]$.

\subsection{A priori estimates and passage to the limit as $n\to+\infty$}\label{sec:apriori:weak}
We now derive a priori estimates to show that $t_{n} = T$ for all $n \in \N$ and the approximate solutions $(\varphi_{n}, \mu_{n}, \sigma_{n}, \bm{v}_{n})$ are uniformly bounded with respect to $n$ in suitable function spaces.
Below the symbol $C$ denotes a generic positive constant that is independent of the parameter $n$ and may vary from line to line.

\paragraph{First estimate.} We obtain from testing \eqref{Galerkin:varphi} with $\mu_{n}$, testing \eqref{Galerkin:sigma} with $\sigma_{n} + \chi (1-\varphi_{n})$, testing \eqref{Galerkin:velo} with $\bm{v}_{n}$ and testing \eqref{Galerkin:mu} with $\pd_{t} \varphi_{n}$ an analogous energy identity to \eqref{Energy:ZEED}:
\begin{equation}\label{Energy:ZEED:Galerkin}
\begin{aligned}
& \frac{\dd}{\dt} \int_{\Omega} \frac{1}{2} \abs{\bm{v}_{n}}^{2} + A \Psi(\varphi_{n}) + \frac{B}{2} \abs{\nabla \varphi_{n}}^{2} + \frac{1}{2} \abs{\sigma_{n}}^{2} + \chi \sigma_{n} (1-\varphi_{n}) \dx \\
& \qquad + \int_{\Omega} m(\varphi_{n}) \abs{\nabla \mu_{n}}^{2} + n(\varphi_{n}) \abs{\nabla (\sigma_{n} + \chi (1- \varphi_{n}))}^{2} + 2 \eta(\varphi_{n}) \abs{\der \bm{v}_{n}}^{2} \dx  \\
& \quad = \int_{\Omega}  (\sigma_{n} + \chi (1-\varphi_{n})) S + \mu_{n} \Gamma \dx.
\end{aligned}
\end{equation}
Here, we used the fact $\pd_{t}\varphi_{n} \in W_{n}$ so that
\begin{align*}
(\P_{W_{n}}(\Psi'(\varphi_{n})), \pd_{t}\varphi_{n}) = (\Psi'(\varphi_{n}), \P_{W_{n}}(\pd_{t}\varphi_{n})) = (\Psi'(\varphi_{n}),\pd_{t}\varphi_{n}).
\end{align*}
Similarly, during the testing procedure, we shift the projection operators from the nonlinearities onto $\mu_{n}, \sigma_{n} + \chi (1-\varphi_{n})$ and $\bm{v}_{n}$, and this leads to \eqref{Energy:ZEED:Galerkin}.

Besides, by the Sobolev embedding $H^{1} \subset L^{s}$ for $s \in[1,+\infty)$ if $d = 2$ and for $s \in[1, 6]$ if $d = 3$, and the assumption \eqref{Psi:assump}, we see that $\Psi(\varphi_{0}) \in L^{1}$.
Thus, the initial energy satisfies
\begin{align}\label{InitialEnergy}
\int_{\Omega} \frac{1}{2} \abs{\bm{v}_0}^{2}  + A \Psi(\varphi_{0}) + \frac{B}{2} \abs{\nabla \varphi_{0}}^{2} + \frac{1}{2} \abs{\sigma_{0}}^{2} + \chi \sigma_{0} (1-\varphi_{0}) \dx  < +\infty.
\end{align}
Then the initial approximating energy
\begin{align*}
\int_{\Omega} \frac{1}{2} \abs{\bm{v}_{n}(0)}^{2} +  A \Psi(\varphi_{n}(0)) + \frac{B}{2} \abs{\nabla \varphi_{n}(0)}^{2} + \frac{1}{2} \abs{\sigma_{n}(0)}^{2} + \chi \sigma_{n}(0)(1-\varphi_{n}(0)) \dx
\end{align*}
is uniformly bounded independently of $n$ due to the fact that
\begin{align*}
\norm{\P_{W_{n}}(f)}_{X} \leq C \norm{f}_{X} \text{ for any } X \supseteq H^{2}_{N}, \quad \norm{\P_{\bm{Y}_{n}}(\bm{f})}_{\bm{Z}} \leq C \norm{\bm{f}}_{\bm{Z}} \text{ for } \bm{Z} = \bm{L}^{2}_{\sigma} \text{ or } \bm{H}^{1}_{0,\sigma}.
\end{align*}

By the assumptions on $\Gamma$ and $S$, we see that the right-hand side of \eqref{Energy:ZEED:Galerkin} can be estimated as
\begin{equation}\label{Source:est:1}
\begin{aligned}
 \abs{\mathrm{RHS}}
 &\leq \norm{S}_{L^{2}} \left (\norm{\sigma_{n}}_{L^{2}} + \chi \norm{\varphi_{n}}_{L^{2}} + C \right ) + \norm{\mu_{n}}_{L^{1}} \norm{\Gamma}_{L^{\infty}} \\
& \leq C \left ( 1 + \norm{\varphi_{n}}_{L^{2}} + \norm{\sigma_{n}}_{L^{2}} \right ) \norm{S}_{L^{2}} + \norm{\Gamma}_{L^{\infty}} \norm{\mu_{n} - \mean{\mu_{n}}}_{L^{1}} + \abs{\mean{\mu_{n}}} \norm{\Gamma}_{L^{\infty}}.
\end{aligned}
\end{equation}
Note that by the assumption \eqref{Psi:assump} it holds that
\begin{equation}\label{mean:mu}
\begin{aligned}
\abs{\mean{\mu_{n}}} \leq C \left ( \norm{\Psi'(\varphi_{n})}_{L^{1}} + \norm{\sigma_{n}}_{L^{1}} \right ) \leq C \left ( 1 + \norm{\Psi(\varphi_{n})}_{L^{1}} + \norm{\sigma_{n}}_{L^{1}} \right).
\end{aligned}
\end{equation}
Then, applying Poincar\'{e}'s inequality in $L^{1}$, H\"{o}lder's inequality and Young's inequality to \eqref{Source:est:1}, we have
\begin{equation}\label{Source:est:2}
\begin{aligned}
 \abs{\mathrm{RHS}} & \leq C \left ( 1 + \norm{\varphi_{n}}_{L^{2}}^{2} + \norm{\sigma_{n}}_{L^{2}}^{2} \right ) + \norm{S}_{L^{2}}^{2} + C \norm{\Gamma}_{L^{\infty}}^{2} + \frac{m_{0}}{2} \norm{\nabla \mu_{n}}_{\bm{L}^{2}}^{2} \\
 & \quad + C \norm{\Gamma}_{L^{\infty}} ( 1 + \norm{\Psi(\varphi_{n})}_{L^{1}} + \norm{\sigma_{n}}_{L^{1}} ) .
\end{aligned}
\end{equation}
Employing the lower bounds of $\eta$, $m$ and $n$, the inequality
\begin{align*}
\frac{1}{2} \norm{\nabla \sigma_{n}}_{\bm{L}^{2}}^{2} \leq \norm{\nabla (\sigma_{n} - \varphi_{n})}_{\bm{L}^{2}}^{2} + \norm{\nabla \varphi_{n}}_{\bm{L}^{2}}^{2}
\end{align*}
and the quadratic lower bound of $\Psi$ in \eqref{Psi:assump}, we obtain from \eqref{Energy:ZEED:Galerkin} the inequality
\begin{equation}\label{Energy:ZEED:2}
\begin{aligned}
& \frac{\dd}{\dt} \int_{\Omega} \frac{1}{2} \abs{\bm{v}_{n}}^{2} + A\Psi(\varphi_{n}) + \frac{B}{2} \abs{\nabla \varphi_{n}}^{2} + \frac{1}{2} \abs{\sigma_{n}}^{2} + \chi \sigma_{n}(1-\varphi_{n}) \dx \\
& \qquad + \int_{\Omega} \frac{m_{0}}{2} \abs{\nabla \mu_{n}}^{2} + \frac{n_{0}}{2} \abs{\nabla \sigma_{n}}^{2} + 2 \eta_{0} \abs{\der \bm{v}_{n}}^{2} \dx \\
& \quad \leq \norm{S}_{L^{2}}^{2} + C \left (1 + \norm{\Gamma}_{L^{\infty}}^{2} \right ) \left ( 1 + \norm{\Psi(\varphi_{n})}_{L^{1}} + \norm{\nabla \varphi_{n}}_{\bm{L}^{2}}^{2} + \norm{\sigma_{n}}_{L^{2}}^{2} \right ).
\end{aligned}
\end{equation}
 Applying Young's inequality and \eqref{Psi:assump}, we see that (cf. \cite[(4.37)]{GLDarcy} or \cite{GLNeumann})
\begin{equation}\label{Energy:lowerbound}
\begin{aligned}
& \int_{\Omega} A \Psi(\varphi_{n}) + \frac{1}{2} \abs{\sigma_{n}}^{2} + \chi \sigma_{n} (1-\varphi_{n}) \dx \\
& \quad  \geq \left ( A - \frac{2 \chi^{2}}{C_{3}} \right ) \norm{\Psi(\varphi_{n})}_{L^{1}} + \frac{1}{4} \norm{\sigma_{n}}_{L^{2}}^{2} - C.
\end{aligned}
\end{equation}
Then, integrating \eqref{Energy:ZEED:2} in time from $0$ to arbitrary $s \in (0,T]$, by the assumption \eqref{assump:A} and \eqref{Energy:lowerbound}, we find that
\begin{equation*}
\begin{aligned}
& \left ( \norm{\bm{v}_{n}(s)}_{\bm{L}^{2}}^{2} + \norm{\Psi(\varphi_{n}(s))}_{L^{1}} + \norm{\nabla \varphi_{n}(s)}_{\bm{L}^{2}}^{2} + \norm{\sigma_{n}(s)}_{L^{2}}^{2} \right ) \\
& \qquad + \norm{\nabla \mu_{n}}_{L^{2}(0,s;\bm{L}^{2})}^{2} + \norm{\nabla \sigma_{n}}_{L^{2}(0,s;\bm{L}^{2})}^{2} + \norm{\der \bm{v}_{n}}_{L^{2}(0,s;\bm{L}^{2})}^{2} \\
& \quad \leq \norm{S}_{L^{2}(0,T;L^{2})}^{2} + C \int_{0}^{s} \left (1 + \norm{\Gamma}_{L^{\infty}}^{2} \right ) \left ( 1 + \norm{\Psi(\varphi_{n})}_{L^{1}} + \norm{\nabla \varphi_{n}}_{\bm{L}^{2}}^{2} + \norm{\sigma_{n}}_{L^{2}}^{2} \right ) \dt + C.
\end{aligned}
\end{equation*}
Applying Gronwall's lemma in integral form (see e.g., \cite[Lem. 3.1]{GLNeumann}) leads to
\begin{equation*}
\begin{aligned}
& \sup_{t\in [0,T]} \left ( \norm{\bm{v}_{n}(t)}_{\bm{L}^{2}}^{2} + \norm{\Psi(\varphi_{n}(t))}_{L^{1}} + \norm{\nabla \varphi_{n}(t)}_{\bm{L}^{2}}^{2} + \norm{\sigma_{n}(t)}_{L^{2}}^{2} \right ) \\
& \quad + \norm{\nabla \mu_{n}}_{L^{2}(Q)}^{2} + \norm{\nabla \sigma_{n}}_{L^{2}(Q)}^{2} + \norm{\der \bm{v}_{n}}_{L^{2}(Q)}^{2} \leq C.
\end{aligned}
\end{equation*}
Using Korn's inequality, Poincar\'{e} inequality, \eqref{Psi:assump} and \eqref{mean:mu}, we see that
\begin{equation}\label{Apriori:1}
\begin{aligned}
& \norm{\bm{v}_{n}}_{L^{\infty}(0,T;\bm{L}^{2}) \cap L^{2}(0,T;\bm{H}^{1})} + \norm{\varphi_{n}}_{L^{\infty}(0,T;H^{1})}+ \norm{\Psi(\varphi_{n})}_{L^{\infty}(0,T;L^{1})} \\
 & \quad  + \norm{\sigma_{n}}_{L^{\infty}(0,T;L^{2}) \cap L^{2}(0,T;H^{1})} + \norm{\mu_{n}}_{L^{2}(0,T;H^{1})} \leq C.
\end{aligned}
\end{equation}
The uniform estimate \eqref{Apriori:1} ensures that we can extend $(\varphi_{n}, \mu_{n}, \sigma_{n}, \bm{v}_{n})$ to the full time interval $[0,T]$, and so $t_{n} = T$ for all $n \in \N$.

\paragraph{Second estimate.} In order to obtain higher order estimate for $\varphi_{n}$, we view \eqref{Galerkin:mu} as an elliptic equation for $\varphi_{n}$ subject to homogeneous Neumann boundary condition.  For $d = 2$, we see that
\begin{align*}
 \norm{\P_{W_{n}}(\Psi'(\varphi_{n}))}_{H^{1}}
 &\leq  C \left (\norm{\nabla(\Psi'(\varphi_{n}))}_{\bm{L}^{2}} + \norm{\Psi'(\varphi_{n})}_{L^{2}} \right ) \\
&  \leq C \norm{\Psi''(\varphi_{n})}_{L^{4}} \norm{\nabla \varphi_{n}}_{\bm{L}^{4}} + C \norm{\Psi'(\varphi_{n})}_{L^{2}} \\
& \leq C \left (1+\norm{\varphi_{n}}_{L^{4r}}^{r} \right ) \norm{\varphi_{n}}_{H^{2}}^{\frac{1}{2}} \norm{\nabla \varphi_{n}}_{\bm{L}^{2}}^{\frac{1}{2}} + C \left (1+\norm{\varphi_{n}}_{L^{2(r+1)}}^{r+1} \right ) \nonumber\\
& \leq C\norm{\varphi_{n}}_{H^{2}}^{\frac{1}{2}} + C,
\end{align*}
where we have used \eqref{Psi:assump}, the Gagliardo--Nirenberg inequality $\norm{f}_{L^{4}}\leq C\norm{f}_{H^{1}}^{\frac{1}{2}} \norm{f}_{L^{2}}^{\frac{1}{2}}$ for two dimensions, and the fact that \eqref{Apriori:1} implies that $\varphi_{n} \in L^{\infty}(0,T;L^{s})$ for all $s \in[1,+\infty)$.  Then by the elliptic estimate for \eqref{Galerkin:mu} and Young's inequality, we get
\begin{align}
\norm{\varphi_{n}}_{H^3}
&\leq C\left(\norm{\mu_{n}}_{H^{1}}+ A \norm{\P_{W_{n}}(\Psi'(\varphi_{n}))}_{H^{1}} + \chi \norm{\sigma_{n}}_{H^{1}}+\norm{\varphi_{n}}_{L^{2}}\right)\nonumber\\
&\leq {\frac{1}{2}}\norm{\varphi_{n}}_{H^{2}}+C\left(\norm{\mu_{n}}_{H^{1}}+ \norm{\sigma_{n}}_{H^{1}}+1\right),
\end{align}
which together with \eqref{Apriori:1} yields the higher order estimate
\begin{align}\label{Apriori:2}
 \norm{\varphi_{n}}_{L^{2}(0,T;H^{3})} + \norm{\Psi'(\varphi_{n})}_{L^{2}(0,T;H^{1})} \leq C.
\end{align}
In the case $d = 3$, as $\Psi$ has at most polynomial growth of order $6$, we can appeal to the bootstrapping argument used in \cite[\S 4.2]{GLDarcy} to obtain the higher order estimate \eqref{Apriori:2}.

\paragraph{Third estimate.} Using that $\nabla \varphi_{n} \in L^{2}(0,T;\bm{H}^{2}) \subset L^{2}(0,T;\bm{L}^{\infty})$ and $\bm{v}_{n} \in L^{\infty}(0,T;\bm{L}^{2})$, we have
\begin{align}\label{Apriori:3:a}
\norm{\bm{v}_{n} \cdot \nabla \varphi_{n}}_{L^{2}(Q)} \leq C.
\end{align}
Then, by testing \eqref{Galerkin:varphi} with an arbitrary test function $\zeta \in L^{2}(0,T;H^{1})$, we easily see from \eqref{Apriori:1} and \eqref{Apriori:3:a} that
\begin{align}\label{Apriori:3:b}
\norm{\pd_{t}\varphi_{n}}_{L^{2}(0,T;(H^{1})')} \leq C.
\end{align}

\paragraph{Fourth estimate.}
We now estimate the convection term for $\sigma$.  In two dimensions, we have the following Gagliardo--Nirenberg inequalities
\begin{align*}
\norm{f}_{L^{2s}} \leq C \norm{f}_{H^{1}}^{\frac{s-1}{s}} \norm{f}_{L^{2}}^{\frac{1}{s}}, \quad \norm{f}_{L^{\frac{2s}{s-1}}} \leq C \norm{f}_{H^{1}}^{\frac{1}{s}} \norm{f}_{L^{2}}^{\frac{s-1}{s}},
\end{align*}
for any $s > 1$.  Using boundedness of $\bm{v}_{n}$ and $\sigma_{n}$ in $L^{\infty}(0,T;L^{2}) \cap L^{2}(0,T;H^{1})$, for arbitrary $\zeta \in L^{2}(0,T;H^{1})$, we see that
\begin{align*}
& \abs{\int_{0}^{T} \int_{\Omega} \sigma_{n} \bm{v}_{n} \cdot \nabla \zeta \dx \dt} \leq \int_{0}^{T} \norm{\sigma_{n}}_{L^{2s}} \norm{\bm{v}_{n}}_{\bm{L}^{\frac{2s}{s-1}}} \norm{\nabla \zeta}_{\bm{L}^{2}} \dt \\
& \quad \leq C \norm{\sigma_{n}}_{L^{\infty}(0,T;L^{2})}^{\frac{1}{s}} \norm{\bm{v}_{n}}_{L^{\infty}(0,T;\bm{L}^{2})}^{\frac{s-1}{s}} \norm{\sigma_{n}}_{L^{2}(0,T;H^{1})}^{\frac{s-1}{s}} \norm{\bm{v}_{n}}_{L^{2}(0,T;\bm{H}^{1})}^{\frac{1}{s}} \norm{\nabla \zeta}_{L^{2}(Q)}
\end{align*}
for any $s > 1$.  This implies that $\bm{v}_{n} \cdot \nabla \sigma_{n}$ is uniformly bounded in $L^{2}(0,T;(H^{1})')$.
Then from the equation \eqref{Galerkin:sigma} we find that
\begin{align}\label{Apriori:4:2D}
\norm{\pd_{t} \sigma_{n}}_{L^{2}(0,T;(H^{1})')} + \norm{\bm{v}_{n} \cdot \nabla \sigma_{n}}_{L^{2}(0,T;(H^{1})')} \leq C.
\end{align}
In three dimensions, the Gagliardo--Nirenberg inequality $\norm{f}_{L^{3}} \leq C \norm{f}_{H^{1}}^{\frac{1}{2}} \norm{f}_{L^{2}}^{\frac{1}{2}}$ yields that for arbitrary $\zeta \in L^{4}(0,T;H^{1})$,
\begin{align*}
& \abs{\int_{0}^{T} \int_{\Omega} (\bm{v}_{n} \cdot \nabla \sigma_{n}) \zeta \dx \dt}
\leq \int_{0}^{T}  \norm{\bm{v}_{n}}_{\bm{L}^{3}} \norm{\nabla \sigma_{n}}_{\bm{L}^{2}} \norm{\zeta}_{L^{6}} \dt \\
& \quad \leq C \norm{\bm{v}_{n}}_{L^{\infty}(0,T;\bm{L}^{2})}^{\frac{1}{2}} \norm{\sigma_{n}}_{L^{2}(0,T;H^{1})}  \norm{\bm{v}_{n}}_{L^{2}(0,T;\bm{H}^{1})}^{\frac{1}{2}} \norm{\zeta}_{L^{4}(0,T;H^{1})},
\end{align*}
which implies
\begin{align}\label{Apriori:4:3D}
\norm{\pd_{t} \sigma_{n}}_{L^{\frac{4}{3}}(0,T;(H^{1})')} + \norm{\bm{v}_{n} \cdot \nabla \sigma_{n}}_{L^{\frac{4}{3}}(0,T;(H^{1})')} \leq C.
\end{align}

\paragraph{Fifth estimate.} We proceed to estimate the convection term in the momentum equation.
Note that, by the boundedness of $\mu_{n} + \chi \sigma_{n}$ in $L^{2}(0,T;H^{1})$ and $\nabla \varphi_{n}$ in $L^{\infty}(0,T;\bm{L}^{2})$, it holds
\begin{align*}
\norm{(\mu_{n} + \chi \sigma_{n})\nabla \varphi_{n}}_{L^{2}(0,T;\bm{L}^{\frac{3}{2}})} \leq \norm{\nabla \varphi_{n}}_{L^{\infty}(0,T;\bm{L}^{2})} \norm{\mu_{n} + \chi \sigma_{n}}_{L^{2}(0,T;L^{6})},
\end{align*}
and so we find that $(\mu_{n} + \chi \sigma_{n})\nabla \varphi_{n}$ is bounded in $L^{2}(0,T;\bm{L}^{\frac{3}{2}}) \subset L^{2}(0,T;\bm{H}^{-1})$.
  In two dimensions, thanks to the boundedness of $\bm{v}_{n}$ in $L^{\infty}(0,T;\bm{L}^{2}) \cap L^{2}(0,T;\bm{H}^{1})$
  and the Gagliardo--Nirenberg inequality $\norm{f}_{L^{4}} \leq C \norm{f}_{L^{2}}^{\frac{1}{2}} \norm{f}_{H^{1}}^{\frac{1}{2}}$,
  we see that for arbitrary $\bm{\zeta} \in L^{2}(0,T;\bm{H}^{1})$,
\begin{align*}
& \abs{ \int_{0}^{T} \int_{\Omega} (\bm{v}_{n} \otimes \bm{v}_{n}): \nabla \bm{\zeta} \dx \dt}  \leq \int_{0}^{T} \norm{\bm{v}_{n}}_{\bm{L}^{4}}^{2} \norm{\nabla \bm{\zeta}}_{\bm{L}^{2}}  \dt \\
& \quad \leq  C \norm{\bm{v}_{n}}_{L^{\infty}(0,T;\bm{L}^{2})} \norm{\bm{v}_{n}}_{L^{2}(0,T;\bm{H}^{1})} \norm{\bm{\zeta}}_{L^{2}(0,T;\bm{H}^{1})}.
\end{align*}
  This implies that
\begin{align}\label{Apriori:5:2D}
\norm{\pd_{t} \bm{v}_{n}}_{L^{2}(0,T;\bm{H}^{-1})} + \norm{(\bm{v}_{n} \cdot \nabla) \bm{v}_{n}}_{L^{2}(0,T;\bm{H}^{-1})} \leq C.
\end{align}
In three dimensions, the Gagliardo--Nirenberg inequality $\norm{f}_{L^{4}} \leq C \norm{f}_{L^{2}}^{\frac{1}{4}} \norm{f}_{H^{1}}^{\frac{3}{4}}$ yields for arbitrary $\bm{\zeta} \in L^{4}(0,T;\bm{H}^{1})$,
\begin{align*}
& \abs{ \int_{0}^{T} \int_{\Omega} (\bm{v}_{n} \otimes \bm{v}_{n}): \nabla \bm{\zeta} \dx \dt}  \leq \int_{0}^{T} \norm{\bm{v}_{n}}_{\bm{L}^{4}}^{2} \norm{\nabla \bm{\zeta}}_{\bm{L}^{2}} \dx \\
& \quad \leq C \norm{\bm{v}_{n}}_{L^{\infty}(0,T;\bm{L}^{2})}^{\frac{1}{2}} \norm{\bm{v}_{n}}_{L^{2}(0,T;\bm{H}^{1})}^{\frac{3}{2}} \norm{\bm{\zeta}}_{L^{4}(0,T;\bm{H}^{1})}.
\end{align*}
Hence, we obtain
\begin{align}\label{Apriori:5:3D}
\norm{\pd_{t} \bm{v}_{n}}_{L^{\frac{4}{3}}(0,T;\bm{H}^{-1})} + \norm{(\bm{v}_{n} \cdot \nabla) \bm{v}_{n}}_{L^{\frac{4}{3}}(0,T;\bm{H}^{-1})} \leq C.
\end{align}

\paragraph{Passage to the limit as $n\to+\infty$.} The above a priori estimates \eqref{Apriori:1}--\eqref{Apriori:5:3D} are sufficient to obtain compactness results that allow us to obtain a convergent subsequence of the approximate solutions $(\varphi_{n}, \mu_{n}, \sigma_{n}, \bm{v}_{n})$, whose limit denoted by $(\varphi, \mu, \sigma, \bm{v})$ satisfies the weak formulation \eqref{ZEED:weakform} as well as the initial conditions.  Since the procedure is standard (see, e.g., \cite{Boyer} for the Navier--Stokes--Cahn--Hilliard system), we omit the details.

\section{Proof of Theorem \ref{thm:2D:strong}: global strong solutions in two dimensions}\label{sec:ZEED:2D:strong}

Keeping the a priori estimates \eqref{Apriori:1}--\eqref{Apriori:3:b} in mind, we proceed to derive further higher order estimates for the approximate solutions $(\varphi_{n}, \mu_{n}, \sigma_{n}, \bm{v}_{n})$ obtained via the Galerkin scheme.

\paragraph{First estimate.} Testing \eqref{Galerkin:varphi} by $\Laplace^{2} \varphi_{n} \in W_{n}$ and using \eqref{Galerkin:mu}, we get
\begin{equation}\label{StrSol:2D:Est:phi:1}
\begin{aligned}
& \frac{\dd}{\dt} \int_\Omega \frac{1}{2} \abs{\Laplace \varphi_{n}}^{2} \dx + B \int_{\Omega} m(\varphi_{n}) \abs{\Laplace^{2} \varphi_{n}}^{2} \dx \\
&\quad = - \int_{\Omega} ( \bm{v}_{n} \cdot \nabla \varphi_{n}) \Laplace^{2} \varphi_{n} \dx +  A\int_{\Omega} m(\varphi_{n}) \Laplace( \P_{W_{n}}(\Psi'(\varphi_{n})) \Laplace^{2} \varphi_{n} \dx  \\
& \qquad  + \int_{\Omega} m'(\varphi_{n}) (\nabla \varphi_{n} \cdot \nabla \mu_{n}) \Laplace^{2} \varphi_{n} \dx
-\chi \int_{\Omega} m(\varphi_{n}) \Laplace \sigma_{n} \Laplace^{2} \varphi_{n} \dx \\
& \qquad + \int_{\Omega} \Gamma \Laplace^{2} \varphi_{n} \dx \\
&\quad =: R_{1} + R_{2} + R_{3} + R_{4} + R_{5}.
\end{aligned}
\end{equation}
The terms on the right-hand side of \eqref{StrSol:2D:Est:phi:1} can be estimated as follows.  Using H\"older's inequality, Young's inequality and the Gagliardo--Nirenberg inequality (for $d = 2$)
\begin{align}\label{vphiGN}
\norm{\nabla\varphi_{n}}_{\bm{L}^{\infty}}\leq C\norm{\nabla\varphi_{n}}_{\bm{H}^{3}}^{\frac{1}{3}} \norm{\nabla\varphi_{n}}_{\bm{L}^{2}}^{\frac{2}{3}},
\end{align}
the boundedness of $\bm{v}_{n}$ in $L^{\infty}(0,T;\bm{L}^{2})$ and $\varphi_{n}$ in $L^{\infty}(0,T;H^{1})$, and the elliptic estimate \eqref{H4EllEst}, we have
\begin{equation}\label{R1}
\begin{aligned}
R_{1} & \leq \norm{\bm{v}_{n}}_{\bm{L}^{2}} \norm{\nabla \varphi_{n}}_{\bm{L}^{\infty}} \norm{\Laplace^{2} \varphi_{n}}_{L^{2}}\nonumber\\
 &\leq C\norm{\nabla \varphi_{n}}_{\bm{H}^{3}}^{\frac{1}{3}} \norm{\nabla \varphi_{n}}_{\bm{L}^{2}}^{\frac{2}{3}} \norm{\Laplace^{2} \varphi_{n}}_{L^{2}} \\
&\leq C \left (\norm{\Laplace^{2} \varphi_{n}}_{L^{2}}+ \norm{\varphi_{n}}_{L^{2}} \right )^{\frac{1}{3}} \norm{\varphi_{n}}_{H^{1}}^{\frac{2}{3}} \norm{\Laplace^{2} \varphi_{n}}_{L^{2}} \\
& \leq \eps \norm{\Laplace^{2} \varphi_{n}}_{L^{2}}^{2} + C,
\end{aligned}
\end{equation}
for some positive constant $\eps$ yet to be determined.  Next, noticing that from integration by parts and using $\pdnu \varphi_{n} = \pdnu (\Laplace \varphi_{n}) = 0$ on $\pd \Omega$, it holds that
\begin{align*}
\norm{\Laplace  \varphi_{n}}_{L^{2}}^{2} = \int_{\Omega} \varphi_{n} \Laplace^{2} \varphi_{n} \dx \leq \norm{\Laplace^{2} \varphi_{n}}_{L^{2}} \norm{\varphi_{n}}_{L^{2}}.
\end{align*}
Then, together with the estimate \eqref{Apriori:1}, \eqref{H4EllEst}, Young's inequality, the Gagliardo--Nirenberg inequality $\norm{\nabla \varphi_{n}}_{\bm{L}^{4}} \leq C\norm{\varphi_{n}}_{H^{4}}^{\frac{1}{6}} \norm{\nabla \varphi_{n}}_{\bm{L}^{2}}^{\frac{5}{6}}$, we have
\begin{align*}
& \norm{\Laplace(\P_{W_{n}}(\Psi'(\varphi_{n})))}_{L^{2}}\nonumber\\
&\quad \leq C \left (\norm{\Laplace(\Psi'(\varphi_{n}))}_{L^{2}} + \norm{\Psi'(\varphi_{n})}_{L^{2}} \right ) \\
&\quad \leq C \norm{\Psi'''(\varphi_{n})}_{L^{\infty}} \norm{\nabla \varphi_{n}}_{\bm{L}^{4}}^{2} + C\norm{\Psi''(\varphi_{n})}_{L^{\infty}} \norm{\Laplace \varphi_{n}}_{L^{2}} + C \norm{\Psi'(\varphi_{n})}_{L^{2}} \\
&\quad \leq C \left (1+\norm{\varphi_{n}}_{L^{\infty}}^{r-1} \right ) \norm{\varphi_{n}}_{H^{4}}^{\frac{1}{3}} \norm{\nabla \varphi_{n}}_{\bm{L}^{2}}^{\frac{5}{3}} + C \left (1+\norm{\varphi_{n}}_{L^{\infty}}^{r} \right )\norm{\Laplace^{2} \varphi_{n}}_{L^{2}}^{\frac{1}{2}} \norm{\varphi}_{L^{2}}^{\frac{1}{2}} \\
&\qquad + C \left (1+\norm{\varphi_{n}}_{L^{\infty}}^{r+1} \right ) \\
&\quad \leq C \left (1 + \norm{\varphi_{n}}_{L^{\infty}}^{r-1} \right ) \norm{\Laplace^{2} \varphi_{n}}_{L^{2}}^{\frac{1}{3}} + C \left (1+\norm{\varphi_{n}}_{L^{\infty}}^{r} \right )\norm{\Laplace^{2} \varphi_{n}}_{L^{2}}^{\frac{1}{2}}+ C \left (1+\norm{\varphi_{n}}_{L^{\infty}}^{r+1} \right )\\
&\quad \leq C \left (1+\norm{\varphi_{n}}_{L^{\infty}}^{3r} \right ) \norm{\Laplace^{2} \varphi_{n}}_{L^{2}}^{\frac{1}{2}} + C \left (1 + \norm{\varphi_{n}}_{L^{\infty}}^{r+1} \right ).
\end{align*}
Applying the Br\'{e}zis--Gallouet inequality \eqref{Bre} for $\varphi_{n}$, using the estimates \eqref{H2EllEst}, \eqref{Apriori:1} and Young's inequality, we deduce that
\begin{equation}\label{Psi1}
\begin{aligned}
& \norm{\Laplace (\P_{W_{n}}(\Psi'(\varphi_{n})))}_{L^{2}} \\
& \quad \leq C \left [ 1+(\ln(1+\norm{\Laplace \varphi_{n}}_{L^{2}}))^{\frac{3r}{2}} \right ] \norm{\Laplace^{2} \varphi_{n}}_{L^{2}}^{\frac{1}{2}} + C \left [1+(\ln(1+\norm{\Laplace \varphi_{n}}_{L^{2}}))^{\frac{r+1}{2}} \right ]\\
& \quad \leq C \big(1+\norm{\Laplace \varphi_{n}}_{L^{2}}^{\frac{1}{2}}\big) \norm{\Laplace^{2} \varphi_{n}}_{L^{2}}^{\frac{1}{2}} + C \norm{\Laplace \varphi_{n}}_{L^{2}} + C.
\end{aligned}
\end{equation}
For the second inequality in \eqref{Psi1} we used the fact that for any $k \geq 1$,
\begin{align*}
z(x) := \frac{(\ln(1+x))^{k}}{x^{\frac{1}{2}}} \to 0\quad  \text{ as } x \to +\infty,
\end{align*}
as its derivative $z'(x) = \frac{(\ln(1+x))^{k-1}}{x^{\frac{1}{2}}} \left ( \frac{k}{1+x} - \frac{\ln(1+x)}{2x} \right )$ becomes negative for sufficiently large $x$.  As a consequence,
\begin{align}\label{R2}
R_{2} & \leq A m_{1} \norm{\Laplace (\P_{W_{n}}(\Psi'(\varphi_{n})))}_{L^{2}} \norm{\Laplace^{2} \varphi_{n}}_{L^{2}}\nonumber\\
 &\leq \eps \norm{\Laplace^{2} \varphi_{n}}_{L^{2}}^{2} + C \norm{\Laplace \varphi_{n}}_{L^{2}}^{2} + C.
\end{align}
Next, using the following estimate
\begin{align*}
\norm{\Laplace \mu_{n}}_{L^{2}} & \leq A\norm{\Laplace(\P_{W_{n}}(\Psi'(\varphi_{n})))}_{L^{2}} + B \norm{\Laplace^{2}\varphi_{n}}_{L^{2}}+\chi \norm{\Laplace \sigma_{n}}_{L^{2}} \\
&\leq B \norm{\Laplace^{2}\varphi_{n}}_{L^{2}} + \chi\norm{\Laplace \sigma_{n}}_{L^{2}} + C \norm{\Laplace \varphi_{n}}_{L^{2}}^{\frac{1}{2}} \norm{\Laplace^{2} \varphi_{n}}_{L^{2}}^{\frac{1}{2}} + C\norm{\Laplace \varphi_{n}}_{L^{2}}+ C \\
&\leq C \left (1 + \norm{\Laplace^{2} \varphi_{n}}_{L^{2}}+ \norm{\Laplace \varphi_{n}}_{L^{2}} + \norm{\Laplace \sigma_{n}}_{L^{2}} \right )
\end{align*}
 we obtain from \eqref{H2EllEst} and Young's inequality that for $R_{3}$,
\begin{equation}\label{R3}
\begin{aligned}
R_{3} & \leq \norm{m'(\varphi_{n})}_{L^{\infty}}\norm{\nabla \varphi_{n}}_{L^{4}}\norm{\nabla \mu_{n}}_{L^{4}}\norm{\Laplace^{2}\varphi_{n}}_{L^{2}} \\
&\leq C \norm{\varphi_{n}}_{H^{1}}^{\frac{1}{2}}\norm{\varphi_{n}}_{H^{2}}^{\frac{1}{2}}\norm{\mu_{n}}_{H^{1}}^{\frac{1}{2}}\norm{\mu_{n}}_{H^{2}}^{\frac{1}{2}}\norm{ \Laplace^{2}\varphi_{n} }_{L^{2}} \\
&\leq C \left (1+\norm{\Laplace \varphi_{n}}_{L^{2}}^{\frac{1}{2}} \right ) \norm{\mu_{n}}_{H^{1}}^{\frac{1}{2}} \left (\norm{\mu_{n}}_{H^{1}}^{\frac{1}{2}}+ \norm{\Laplace \mu_{n}}_{L^{2}}^{\frac{1}{2}} \right )\norm{\Laplace^{2} \varphi_{n}}_{L^{2}} \\
&\leq C \left ( 1+\norm{\Laplace \varphi_{n}}_{L^{2}}^{\frac{1}{2}} \right ) \norm{\mu_{n}}_{H^{1}}\norm{\Laplace^{2} \varphi_{n}}_{L^{2}} \\
&\quad + C\norm{\mu_{n}}_{H^{1}}^{\frac{1}{2}}\left ( 1 + \norm{\Laplace^{2} \varphi_{n}}_{L^{2}} + \norm{\Laplace \varphi_{n}}_{L^{2}} + \norm{\Laplace \sigma_{n}}_{L^{2}}  \right )^{\frac{1}{2}}\norm{\Laplace^{2}\varphi_{n}}_{L^{2}} \\
&\quad + C \norm{\mu_{n}}_{H^{1}}^{\frac{1}{2}}\norm{\Laplace \varphi_{n}}_{L^{2}}^{\frac{1}{2}} \left (1 + \norm{\Laplace^{2} \varphi_{n}}_{L^{2}}+ \norm{\Laplace \varphi_{n}}_{L^{2}} + \norm{\Laplace \sigma_{n}}_{L^{2}}  \right )^{\frac{1}{2}}\norm{\Laplace^{2} \varphi_{n}}_{L^{2}}\\
&\leq \eps \norm{\Laplace^{2} \varphi_{n}}_{L^{2}}^{2} + \eps \norm{\Laplace \sigma_{n}}_{L^{2}}^{2} +  C\left (1+\norm{\mu_{n}}_{H^{1}}^{2} \right )\norm{\Laplace \varphi_{n}}_{L^{2}}^{2} + C\norm{\mu_{n}}_{H^{1}}^{2} + C.
\end{aligned}
\end{equation}
By Young's inequality, we easily have
\begin{align}\label{R4}
R_{4} \leq \chi m_{1} \norm{\Laplace \sigma_{n}}_{L^{2}} \norm{\Laplace^{2} \varphi_{n}}_{L^{2}} \leq \eps \norm{\Laplace^{2} \varphi_{n}}_{L^{2}}^{2} + \frac{\chi^{2} m_{1}^{2}} {4\eps} \norm{\Laplace \sigma_{n}}_{L^{2}}^{2},
\end{align}
and
\begin{align}\label{R5}
R_{5} \leq \norm{\Gamma}_{L^{2}} \norm{\Laplace^{2} \varphi_{n}}_{L^{2}} \leq \eps \norm{\Laplace^{2} \varphi_{n}}_{L^{2}}^{2} + C \norm{\Gamma}_{L^{2}}^{2}.
\end{align}
Now testing \eqref{Galerkin:sigma} by $-\Laplace \sigma_{n} \in W_{n}$, we get
\begin{equation}\label{StrSol:2D:Est:sigm:1}
\begin{aligned}
& \frac{\dd}{\dt} \int_{\Omega} {\frac{1}{2}} \abs{\nabla \sigma_{n}}^{2} \dx +  \int_{\Omega} n(\varphi_{n}) \abs{\Laplace \sigma_{n}}^{2} \dx\\
& \quad =  \int_{\Omega} (\bm{v}_{n} \cdot \nabla \sigma_{n}) \Laplace \sigma_{n} \dx
+\chi \int_{\Omega} n(\varphi_{n}) \Laplace \varphi_{n} \Laplace \sigma_{n} \dx \\
&\qquad -\int_{\Omega} n'(\varphi_{n}) \left ( \nabla \varphi_{n} \cdot \nabla \left (\sigma_{n}-\chi \varphi_{n} \right ) \right ) \Laplace \sigma_{n} -\int_{\Omega} S \Laplace \sigma_{n} \dx \\
&\quad := R_{6} + R_{7} + R_{8} + R_{9}.
\end{aligned}
\end{equation}
The terms on the right-hand side of \eqref{StrSol:2D:Est:sigm:1} can be estimated as follows.  By the Gagliardo--Nirenberg inequality $\norm{f}_{L^{2}} \leq C \norm{f}_{H^{1}}^{\frac{1}{2}} \norm{f}_{L^{2}}^{\frac{1}{2}}$, \eqref{H2EllEst} and the estimate \eqref{Apriori:1}, it holds
\begin{equation}\label{R6}
\begin{aligned}
R_{6} & \leq C\norm{\bm{v}_{n}}_{\bm{L}^{4}} \norm{ \nabla \sigma_{n}}_{\bm{L}^{4}} \norm{\Laplace \sigma_{n}}_{L^{2}} \nonumber\\
&\leq C \norm{\bm{v}_{n}}_{\bm{L}^{2}}^{\frac{1}{2}} \norm{\nabla \bm{v}_{n}}_{\bm{L}^{2}}^{\frac{1}{2}} \norm{\nabla \sigma_{n} }_{\bm{L}^{2}}^{\frac{1}{2}}
\norm{ \sigma_{n} }_{H^{2}}^{\frac{1}{2}} \norm{ \Laplace \sigma_{n}}_{L^{2}}\\
&\leq C \norm{\nabla \bm{v}_{n}}_{\bm{L}^{2}}^{\frac{1}{2}} \norm{\nabla \sigma_{n} }_{\bm{L}^{2}}^{\frac{1}{2}} \left ( \norm{\sigma_{n}}_{L^{2}}^{\frac{1}{2}} +\norm{\Laplace \sigma_{n} }_{L^{2}}^{\frac{1}{2}} \right ) \norm{ \Laplace \sigma_{n}}_{L^{2}} \\
&\leq \eps \norm{ \Laplace \sigma_{n}}_{L^{2}}^{2} +C \left (1+\norm{\nabla \bm{v}_{n}}_{\bm{L}^{2}}^{2} \right ) \norm{\nabla \sigma_{n} }_{\bm{L}^{2}}^{2}.
\end{aligned}
\end{equation}
By Young's inequality, we easily have
\begin{align}\label{R7R9}
R_{7} \leq \eps \norm{ \Laplace \sigma_{n}}_{L^{2}}^{2} + C\norm{\Laplace \varphi_{n}}_{L^{2}}^{2}, \quad R_{9} \leq \eps \norm{ \Laplace \sigma_{n} }_{L^{2}}^{2} + C\norm{S}_{L^{2}}^{2}.
\end{align}
Meanwhile, for $R_{8}$ we have by applying the estimate \eqref{Apriori:1} and \eqref{H2EllEst},
\begin{equation}\label{R8}
\begin{aligned}
R_{8} &\leq \norm{ n'(\varphi_{n})}_{L^{\infty}} \norm{ \nabla \varphi_{n}}_{\bm{L}^{4}} \norm{ \nabla \sigma_{n}}_{\bm{L}^{4}} \norm{ \Laplace \sigma_{n} }_{L^{2}} +  \chi \norm{n'(\varphi_{n})}_{L^{\infty}} \norm{\nabla \varphi_{n}}_{\bm{L}^{4}}^{2} \norm{ \Laplace \sigma_{n}}_{L^{2}}\\
&\leq C \norm{\varphi_{n}}_{H^{2}}^{\frac{1}{2}} \norm{ \nabla \varphi_{n}}_{\bm{L}^{2}}^{\frac{1}{2}}        \norm{\sigma_{n}}_{H^{2}}^{\frac{1}{2}} \norm{ \nabla \sigma_{n}}_{\bm{L}^{2} }^{\frac{1}{2}}\norm{ \Laplace \sigma_{n} }_{L^{2}} + C \norm{\varphi_{n}}_{H^{2}}\norm{\nabla \varphi_{n}}_{\bm{L}^{2}}\norm{ \Laplace \sigma_{n}}_{L^{2}}\\
&\leq  C \left (\norm{\varphi_{n}}_{L^{2}}^{\frac{1}{2}} + \norm{\Laplace \varphi_{n}}_{L^{2}}^{\frac{1}{2}} \right ) \norm{ \nabla \varphi_{n}}_{\bm{L}^{2}}^{\frac{1}{2}}
\left (\norm{\sigma_{n}}_{L^{2}}^{\frac{1}{2}} + \norm{\Laplace \sigma_{n}}_{L^{2}}^{\frac{1}{2}} \right ) \norm{ \nabla \sigma_{n}}_{\bm{L}^{2}}^{\frac{1}{2}}\norm{ \Laplace \sigma_{n} }_{L^{2}}\\
&\quad + C \left (\norm{\varphi_{n}}_{L^{2}} + \norm{\Laplace \varphi_{n}}_{L^{2}} \right ) \norm{\nabla \varphi_{n}}_{\bm{L}^{2}}\norm{ \Laplace \sigma_{n}}_{L^{2}} \\
&\leq  C \left (1+\norm{\Laplace \varphi_{n}}_{L^{2}}^{\frac{1}{2}} \right ) \left (1+\norm{\Laplace \sigma_{n}}_{L^{2}}^{\frac{1}{2}} \right ) \norm{ \nabla \sigma_{n}}_{\bm{L}^{2}}^{\frac{1}{2}} \norm{ \Laplace \sigma_{n} }_{L^{2}}\\
&\quad + C \left (1 + \norm{\Laplace \varphi_{n}}_{L^{2}} \right ) \norm{ \Laplace \sigma_{n}}_{L^{2}}\\
&\leq \eps \norm{ \Laplace \sigma_{n} }_{L^{2}}^{2} +  C \left (1+\norm{ \nabla \sigma_{n}}_{\bm{L}^{2}}^{2} \right )\norm{\Laplace \varphi_{n}}_{L^{2}}^{2} + C\norm{ \nabla \sigma_{n}}_{\bm{L}^{2}}^{2} +C.
\end{aligned}
\end{equation}
Then, multiplying \eqref{StrSol:2D:Est:sigm:1} by $\eps^{-2}$, adding the resultant with \eqref{StrSol:2D:Est:phi:1}, we deduce from the above estimates \eqref{R1}--\eqref{R5}, \eqref{R6}--\eqref{R8} that
\begin{equation}\label{phisig1}
\begin{aligned}
& \frac{1}{2} \frac{\dd}{\dt}  \left( \norm{\Laplace \varphi_{n}}_{L^{2}}^{2} + \eps^{-2} \norm{\nabla \sigma_{n}}_{\bm{L}^{2}}^{2} \right) + (B m_{0} - 5\eps) \norm{\Laplace^{2} \varphi_{n}}_{L^{2}}^{2} \\
&\qquad + \left(\frac{n_{0}}{\eps^2} - \eps - \frac{4}{\eps}-\frac{\chi^{2} m_{1}^{2}}{4\eps} \right)  \norm{\Laplace \sigma_{n}}_{L^{2}}^{2} \\
&\quad \leq C \left (1+ \norm{\mu_{n}}_{H^{1}}^{2} +\norm{\nabla \bm{v}_{n}}_{\bm{L}^2}^{2} +\norm{\nabla \sigma_{n}}_{\bm{L}^{2}}^{2} \right ) \left (\norm{\Laplace \varphi_{n}}_{L^{2}}^{2} + \norm{\nabla \sigma_{n}}_{\bm{L}^{2}}^{2} \right )\\
&\qquad +C \left (1+ \norm{\mu_{n}}_{H^{1}}^{2}+\norm{\nabla \sigma_{n}}_{\bm{L}^{2}}^{2} + \norm{\Gamma}_{L^{2}}^{2} + \norm{S}_{L^{2}}^{2} \right ).
\end{aligned}
\end{equation}
We can choose
\begin{align*}
\eps =\min \left\{ 1, \ \frac{B m_{0}}{10},\ \frac{2 n_{0}}{20+\chi^{2} m_{1}^{2}}\right\},
\end{align*}
then, by Gronwall's lemma, \eqref{phisig1} and the lower order estimate \eqref{Apriori:1} we have
\begin{align}
\sup_{t\in [0,T]} \left(\norm{\Laplace \varphi_{n}(t)}_{L^{2}}^{2} +  \norm{\nabla \sigma_{n}(t)}_{\bm{L}^{2}}^2\right )+ \norm{\Laplace^{2} \varphi_{n}}_{L^{2}(Q)}^{2} + \norm{\Laplace \sigma_{n}}_{L^{2}(Q)}^{2} \leq C,\nonumber
\end{align}
which implies
\begin{align}\label{Apriori:3}
\norm{\varphi_{n}}_{L^{\infty}(0,T;H^{2})\cap L^{2}(0,T; H^{4})} +  \norm{\sigma_{n}}_{L^{\infty}(0,T; H^{1})\cap L^{2}(0,T; H^{2})} \leq C.
\end{align}
Besides, we infer from \eqref{Galerkin:mu}, \eqref{Psi1} and \eqref{Apriori:3} that
\begin{align}\label{munL2H2}
\norm{\mu_{n}}_{L^{2}(0,T; H^{2})} \leq C.
\end{align}

\paragraph{Second estimate.}  By definition of the Stokes operator $\Stokes$, we can write  $\Stokes \bm{v}_{n}=-\Laplace \bm{v}_{n}+ \nabla \pi_{n}$, which implies that (cf. \cite[(3.30)]{Boyer})
\begin{align*}
-  \int_{\Omega} \eta(\varphi_{n}) \Laplace \bm{v}_{n} \cdot \Stokes \bm{v}_{n} \dx = \int_{\Omega} \eta(\varphi_{n}) \abs{\Stokes \bm{v}_{n}}^{2} \dx
- \int_{\Omega} \eta(\varphi_{n}) \nabla \pi_{n} \cdot \Stokes \bm{v}_{n} \dx.
\end{align*}
Then, testing \eqref{Galerkin:velo} by $\Stokes \bm{v}_{n} \in \bm{Y}_{n}$,  using \eqref{Galerkin:mu} and the above fact, we get
\begin{equation}\label{vnH1}
\begin{aligned}
& \frac{\dd}{\dt} \int_{\Omega} {\frac{1}{2}} \abs{\nabla \bm{v}_{n}}^{2} \dx + \int_{\Omega} \eta(\varphi_{n}) \abs{\Stokes \bm{v}_{n}}^{2} \dx \\
& \quad = - \int_{\Omega} \left [(\bm{v}_{n} \cdot \nabla) \bm{v}_{n} \right ]\cdot \Stokes \bm{v}_{n} \dx  + \int_{\Omega} 2 \eta'(\varphi_{n}) \nabla \varphi_{n} \cdot (\der \bm{v}_{n} \Stokes \bm{v}_{n}) \dx\\
& \qquad + \int_{\Omega} \eta(\varphi_{n}) \nabla \pi_{n} \cdot \Stokes \bm{v}_{n} \dx  - B  \int_{\Omega} \Laplace \varphi_{n}  \nabla \varphi_{n} \cdot \Stokes \bm{v}_{n} \dx \\
&\qquad +  A \int_{\Omega} \P_{W_{n}}(\Psi'(\varphi_{n})) \nabla \varphi_{n} \cdot \Stokes \bm{v}_{n} \dx \\
& \quad =: R_{10} + R_{11} + R_{12} + R_{13} + R_{14}.
\end{aligned}
\end{equation}
Keeping in mind the estimate \eqref{Apriori:3}, the reminder terms on the right hand side of \eqref{vnH1} can be estimate as follows.
Using \eqref{Stokes II}, we have
\begin{align*}
R_{10} & \leq \norm{\bm{v}_{n}}_{\bm{L}^{4}} \norm{\nabla \bm{v}_{n}}_{\bm{L}^{4}} \norm{ \Stokes \bm{v}_{n}}_{\bm{L}^{2}}\nonumber\\
 &\leq C \norm{\bm{v}_{n}}_{\bm{L}^{2}}^{\frac{1}{2}} \norm{\nabla \bm{v}_{n}}_{\bm{L}^{2}} \norm{\nabla \bm{v}_{n}}_{\bm{H}^{1}}^{\frac{1}{2}} \norm{\Stokes \bm{v}_{n}}_{\bm{L}^{2}} \\
& \leq C \norm{\bm{v}_{n}}_{\bm{L}^{2}}^{\frac{1}{2}} \norm{\nabla \bm{v}_{n}}_{\bm{L}^{2}} \norm{ \Stokes \bm{v}_{n}}_{\bm{L}^{2}}^{\frac{3}{2}} \\
& \leq \eps \norm{ \Stokes \bm{v}_{n}}_{\bm{L}^{2}}^{2} + C \norm{\nabla \bm{v}_{n}}_{\bm{L}^{2}}^{4}.
\end{align*}
Using \eqref{vphiGN} and Young's inequality, it holds
\begin{align*}
R_{11} & \leq 2 \norm{\eta'(\varphi_{n})}_{L^{\infty}} \norm{\nabla \varphi_{n}}_{\bm{L}^{\infty}} \norm{\der \bm{v}_{n}}_{\bm{L}^{2}} \norm{ \Stokes \bm{v}_{n}}_{\bm{L}^{2}} \\
&\leq  C \norm{\nabla\varphi_{n}}_{\bm{H}^{3}}^{\frac{1}{3}} \norm{\nabla\varphi_{n}}_{\bm{L}^{2}}^{\frac{2}{3}} \norm{\nabla \bm{v}_{n}}_{\bm{L}^{2}} \norm{ \Stokes \bm{v}_{n}}_{\bm{L}^{2}}\\
&\leq  C \left (\norm{\Laplace^{2} \varphi_{n}}_{L^{2}} + \norm{\varphi_{n}}_{L^{2}} \right )^{\frac{1}{3}} \norm{\nabla \bm{v}_{n}}_{\bm{L}^{2}} \norm{ \Stokes \bm{v}_{n}}_{\bm{L}^{2}} \\
&\leq \eps \norm{ \Stokes \bm{v}_{n}}_{\bm{L}^{2}}^{2} + C \norm{\nabla \bm{v}_{n}}_{\bm{L}^{2}}^{3} + C \norm{\Laplace^{2} \varphi_{n}}_{L^{2}}^{2} + C.
\end{align*}
Applying integration by parts and the estimate for the Stokes problem (see Lemma \ref{stoo}), and the facts that $\div (\Stokes \bm{v}_{n})  = 0$ as well as $\Stokes \bm{v}_{n} \cdot \bm{\nu} = 0$ on $\pd \Omega$, we infer from the higher order estimate \eqref{Apriori:3} that
\begin{align*}
R_{12} &= -\int_{\Omega} \eta'(\varphi_{n}) \pi_{n} \nabla \varphi_{n}  \cdot \Stokes \bm{v}_{n} \dx \\
&\leq \norm{\eta'(\varphi_{n})}_{L^{\infty}} \norm{\pi_{n}}_{L^{4}} \norm{\nabla \varphi_{n}}_{\bm{L}^{4}}  \norm{ \Stokes \bm{v}_{n}}_{\bm{L}^{2}}\\
&\leq C \norm{\pi_{n}}_{L^{2}}^{\frac{1}{2}} \norm{\pi_{n}}_{H^{1}}^{\frac{1}{2}} \norm{\nabla \varphi_{n}}_{\bm{L}^{2}}^{\frac{1}{2}} \norm{\varphi_{n}}_{H^{2}}^{\frac{1}{2}}  \norm{ \Stokes \bm{v}_{n}}_{\bm{L}^{2}} \\
&\leq C_{*} \left (\nu \norm{ \Stokes \bm{v}_{n}}_{\bm{L}^{2}}^{\frac{1}{2}} + C_{\nu} \norm{\nabla \bm{v}_{n}}_{\bm{L}^{2}}^{\frac{1}{2}} \right ) \norm{ \Stokes \bm{v}_{n}}_{\bm{L}^{2}}^{\frac{3}{2}}  \\
&\leq C_{*} \nu \norm{ \Stokes \bm{v}_{n}}_{\bm{L}^{2}}^{2} + \eps \norm{ \Stokes \bm{v}_{n}}_{\bm{L}^{2}}^{2} + \frac{C_{*}^{2} C_{\nu}^{2}}{4\eps} \norm{\nabla \bm{v}_{n}}_{\bm{L}^{2}}^{2},
\end{align*}
where the constant $C_{*}$ is independent of $\eps$ and $\nu$.  Next, for $R_{13}$, we deduce from \eqref{vphiGN} and \eqref{Apriori:3} that
\begin{align*}
R_{13}& \leq B\norm{\Laplace \varphi_{n}}_{L^{2}} \norm{\nabla \varphi_{n} }_{\bm{L}^{\infty}} \norm{\Stokes \bm{v}_{n}}_{\bm{L}^{2}} \\
& \leq \norm{\Laplace \varphi_{n}}_{L^{2}} \norm{\nabla \varphi_{n} }_{\bm{H}^{3}}^{\frac{1}{3}}\norm{\nabla \varphi_{n} }_{\bm{L}^{2}}^{\frac{2}{3}} \norm{\Stokes \bm{v}_{n}}_{\bm{L}^{2}} \\
& \leq \eps \norm{ \Stokes \bm{v}_{n}}_{\bm{L}^{2}}^{2} + C \norm{\Laplace^{2} \varphi_{n}}_{L^{2}}^{2} + C.
\end{align*}
Finally, thanks to the fact that $\varphi_{n} \in L^{\infty}(0,T;H^{2})$, we observe
\begin{align*}
R_{14} & \leq A \norm{\P_{W_{n}}(\Psi'(\varphi_{n}))}_{L^{2}} \norm{\nabla \varphi_{n}}_{\bm{L}^{\infty}} \norm{ \Stokes \bm{v}_{n}}_{\bm{L}^{2}} \\
&\leq C \left (1+\norm{\varphi_{n}}_{L^{2r+2}}^{r+1} \right ) \norm{\nabla \varphi_{n} }_{\bm{H}^{3}}^{\frac{1}{3}} \norm{\nabla \varphi_{n} }_{\bm{L}^{2}} ^{\frac{2}{3}} \norm{\Stokes \bm{v}_{n}}_{\bm{L}^{2}} \\
& \leq \eps \norm{ \Stokes \bm{v}_{n}}_{\bm{L}^{2}}^{2} + C \norm{\Laplace^{2} \varphi_{n}}_{L^{2}}^{2} + C.
\end{align*}
Collecting the above estimates, we infer from \eqref{vnH1} that
\begin{align*}
& {\frac{1}{2}} \frac{\dd}{\dt} \norm{\nabla \bm{v}_{n}}_{\bm{L}^{2}}^{2} + \left (\eta_{0} - 5 \eps -C_{*} \nu \right )\norm{\Stokes \bm{v}_{n}}_{\bm{L}^{2}}^{2} \\
& \quad \leq \left( C \norm{\nabla \bm{v}_{n}}_{\bm{L}^{2}}^{2} + \frac{ C_{*}^{2} C_{\nu}^{2}}{4\eps} \right) \norm{\nabla \bm{v}_{n}}_{\bm{L}^{2}}^{2} + C\norm{\Laplace^{2} \varphi_{n}}_{L^{2}}^{2} + C.
\end{align*}
In the above inequality, choosing
\begin{align*}
\eps = \frac{\eta_{0}}{20}\quad\text{and}\quad  \nu = \frac{\eta_{0}}{4C_{*}},
\end{align*}
then, by Gronwall's lemma and \eqref{Apriori:3}, we get
\begin{align*}
\sup_{t\in [0,T]} \norm{\nabla \bm{v}_{n}(t)}_{\bm{L}^{2}}^{2} + \norm{\Stokes \bm{v}_{n}}_{L^{2}(0,T;\bm{L}^{2})}^{2} \leq C,
\end{align*}
which together with Lemma \ref{stoo} yields that
\begin{align}\label{Apriori:4}
\norm{\bm{v}_{n}}_{L^{\infty}(0,T; \bm{H}^{1})\cap L^{2}(0,T; \bm{H}^{2})}\leq C.
\end{align}

\paragraph{Third estimate.}
We now obtain further estimates for the time derivatives $\pd_{t} \varphi_{n}$, $\pd_{t} \sigma_{n}$ and $\pd_{t} \bm{v}_{n}$.
Thanks to
\begin{align*}
\norm{\P_{W_{n}}(\bm{v}_{n} \cdot \nabla \varphi_{n})}_{L^{2}} & \leq \norm{\bm{v}_{n} \cdot \nabla \varphi_{n}}_{L^{2}} \leq \norm{\bm{v}_{n}}_{\bm{L}^{4}} \norm{\nabla \varphi_{n}}_{\bm{L}^{4}},
\end{align*}
and
\begin{align*}
&\norm{\P_{W_{n}}(\div (m(\varphi_{n}) \nabla \mu_{n}))}_{L^{2}}
\leq \norm{\div (m(\varphi_{n}) \nabla \mu_{n})}_{L^{2}} \\
&\quad \leq \norm{m'(\varphi_{n})}_{L^{\infty}}\norm{\nabla \varphi_{n}}_{\bm{L}^{4}}\norm{\nabla \mu_{n}}_{\bm{L}^{4}} + \norm{m(\varphi_{n})}_{L^{\infty}} \norm{\Delta \mu_{n}}_{\bm{L}^{2}},
\end{align*}
we infer from the equation \eqref{Galerkin:varphi} and the estimates \eqref{Apriori:3}, \eqref{munL2H2} and \eqref{Apriori:4} that
\begin{align}\label{vphitL2L2}
\norm{\pd_t \varphi_{n}}_{L^{2}(0,T; L^{2})}\leq C.
\end{align}
In a similar manner, from
\begin{align*}
\norm{\P_{W_{n}}(\bm{v}_{n} \cdot \nabla \sigma_{n})}_{L^{2}} &\leq \norm{\bm{v}_{n} \cdot \nabla \sigma_{n}}_{L^{2}} \leq \norm{\bm{v}_{n}}_{\bm{L}^{4}} \norm{\nabla \sigma_{n}}_{\bm{L}^{4}},
\end{align*}
and
\begin{align*}
&\norm{\P_{W_{n}}(\div(n(\varphi_{n}) \nabla (\sigma_{n} - \chi \varphi_{n})))}_{L^{2}}\nonumber\\
&\quad \leq \norm{\div(n(\varphi_{n}) \nabla (\sigma_{n} - \chi \varphi_{n}))}_{L^{2}} \\
&\quad \leq \norm{n'(\varphi_{n})}_{L^{\infty}} \norm{\nabla \varphi_{n}}_{\bm{L}^{4}} \left ( \norm{\nabla \sigma_{n}}_{\bm{L}^{4}}+ \chi\norm{\nabla \varphi_{n}}_{\bm{L}^{4}} \right ) \nonumber\\
& \qquad + \norm{n(\varphi_{n})}_{L^{\infty}} \left (\norm{\Laplace \sigma_{n}}_{L^{2}}+ \chi \norm{\Laplace \varphi_{n}}_{L^{2}} \right )
\end{align*}
as well as \eqref{Galerkin:sigma}, we conclude
\begin{align}\label{sigmatL2L2}
\norm{\pd_t \sigma_{n}}_{L^{2}(0,T; L^{2})}\leq C.
\end{align}
Finally, testing \eqref{Galerkin:velo} with $\pd_{t} \bm{v}_{n}$ (again $\pd_{t} \bm{v}_{n} \in \bm{Y}_{n}$ and thus the projection operator $\P_{\bm{Y}_{n}}$ can be shifted to $\pd_{t} \bm{v}_{n}$) and using \eqref{Apriori:3}, \eqref{Apriori:4} we get
\begin{equation}\label{Reg:NS}
\begin{aligned}
& \norm{\pd_{t} \bm{v}_{n}}_{\bm{L}^{2}}^{2} \\
& \quad = \int_{\Omega} \div(2\eta(\varphi_{n}) \der \bm{v}_{n})\cdot \pd_{t} \bm{v}_{n}- (\bm{v}_{n} \cdot \nabla ) \bm{v}_{n} \cdot \pd_{t} \bm{v}_{n} + (\mu_{n} + \chi \sigma_{n}) \nabla \varphi_{n} \cdot \pd_{t} \bm{v}_{n} \dx \\
&\quad \leq 2\norm{\eta(\varphi_{n})}_{L^{\infty}} \norm{\Laplace \bm{v}_{n}}_{\bm{L}^{2}} \norm{\pd_{t} \bm{v}_{n}}_{\bm{L}^{2}} + 2 \norm{\eta'(\varphi_{n})}_{L^{\infty}} \norm{\nabla \varphi_{n}}_{\bm{L}^{\infty}} \norm{\nabla \bm{v}_{n}}_{\bm{L}^{2}} \norm{\pd_{t} \bm{v}_{n}}_{\bm{L}^{2}}\\
&\qquad + \norm{\bm{v}_{n}}_{\bm{L}^{4}} \norm{\nabla \bm{v}_{n}}_{\bm{L}^{4}} \norm{\pd_{t} \bm{v}_{n}}_{\bm{L}^{2}} +  \norm{(\mu_{n} + \chi \sigma_{n}) \nabla \varphi_{n}}_{\bm{L}^{2}} \norm{\pd_{t} \bm{v}_{n}}_{\bm{L}^{2}} \\
& \quad \leq {\frac{1}{2}} \norm{\pd_{t} \bm{v}_{n}}_{\bm{L}^{2}}^{2} + C \norm{\Laplace \bm{v}_{n}}_{\bm{L}^{2}}^{2} + C\norm{\nabla \varphi_{n}}_{\bm{L}^{\infty}}^{2} \norm{\nabla \bm{v}_{n}}_{\bm{L}^{2}}^{2} + C \norm{\bm{v}_{n}}_{\bm{L}^{4}}^{2} \norm{\nabla \bm{v}_{n}}_{\bm{L}^{4}}^{2}  \\
& \qquad + C \norm{(\mu_{n} + \chi \sigma_{n})}_{\bm{L}^{2}}^{2} \norm{\nabla \varphi_{n}}_{\bm{L}^{\infty}}^{2}.
\end{aligned}
\end{equation}
Using that $\mu_{n} + \chi \sigma_{n} \in L^{\infty}(0,T;L^{2})$ and $\nabla \varphi \in L^{2}(0,T;\bm{H}^{2}) \subset L^{2}(0,T;\bm{L}^{\infty})$, the last term on the right-hand side of \eqref{Reg:NS} belongs to $L^{1}(0,T)$.  Furthermore, by the Gagliardo--Nirenberg inequality
\begin{align*}
 \norm{\bm{v}_{n}}_{\bm{L}^{4}}^{2} \norm{\nabla \bm{v}_{n}}_{\bm{L}^{4}}^{2} \leq C\norm{\bm{v}_{n}}_{\bm{L}^{2}} \norm{\nabla \bm{v}_{n}}_{\bm{L}^{2}}^{2} \norm{\nabla \bm{v}_{n}}_{\bm{H}^{1}}
\end{align*}
and by \eqref{Apriori:4}, we can conclude that
\begin{align}
\norm{\pd_{t} \bm{v}_{n}}_{L^{2}(0,T; \bm{L}^{2})}\leq C.
\end{align}

\paragraph{Fourth estimate.}

Taking the time derivative of \eqref{Galerkin:mu} leads to
\begin{align}\label{Galerkin:pdtmu}
\pd_{t} \mu_{n} = A \P_{W_{n}}(\Psi''(\varphi_{n}) \pd_{t} \varphi_{n}) - B \Laplace \pd_{t} \varphi_{n} - \chi \pd_{t} \sigma_{n}.
\end{align}
Then testing \eqref{Galerkin:pdtmu} with $- \Laplace \mu_{n}$ and using integration by parts and \eqref{Galerkin:varphi} yields that
\begin{align}\label{mu:t:mu}
\notag &\frac{\dd}{\dt} \int_{\Omega} \frac{1}{2} \abs{\nabla \mu_{n}}^{2} \dx + B \int_{\Omega} m(\varphi_{n}) \abs{\nabla \Laplace \mu_{n}}^{2} \dx \\
\notag &\quad = -A\int_{\Omega} \Psi''(\varphi_{n}) \pd_{t}\varphi_{n} \Laplace \mu_{n} \dx + \chi \int_{\Omega} \pd_{t}\sigma_{n} \Laplace \mu_{n} \dx - B \int_{\Omega} \nabla \Gamma \cdot \nabla \Laplace \mu_{n} \dx \\
\notag & \qquad - B \int_{\Omega} m'(\varphi_{n}) \Laplace \mu_{n} \nabla \varphi_{n} \cdot \nabla \Laplace \mu_{n} \dx - B \int_{\Omega} \nabla \left ( m'(\varphi_{n}) \nabla \mu_{n} \cdot  \nabla \varphi_{n} \right) \cdot \nabla \Laplace \mu_{n} \dx \\
\notag & \qquad + B  \int_{\Omega} \nabla \left ( \bm{v}_{n} \cdot \nabla \varphi_{n} \right ) \cdot \nabla \Laplace \mu_{n} \dx \\
&\quad =: R_{15} + R_{16} + R_{17} + R_{18} + R_{19} + R_{20},
\end{align}
where we have used the fact that $\Laplace \mu_{n} \in W_{n}$ and thus
\begin{align*}
(\P_{W_{n}}(\Psi''(\varphi_{n}) \pd_{t}\varphi_{n}), \Laplace \mu_{n} ) = ( \Psi''(\varphi_{n}) \pd_{t}\varphi_{n}, \P_{W_{n}}(\Laplace \mu_{n})) = (\Psi''(\varphi_{n}) \pd_{t}\varphi_{n}, \Laplace \mu_{n}).
\end{align*}
Here we point out that the assumption $m \in C^{2}_{b}(\R)$ is used for this estimate.  It follows from \eqref{Apriori:3} and Young's inequality that
\begin{align*}
R_{15}&\leq  A\norm{\Psi''(\varphi_{n})}_{L^{\infty}}\norm{\pd_{t}\varphi_{n}}_{L^{2}}\norm{\Laplace \mu_{n}}_{L^{2}}\nonumber\\
&\leq C \left (1+\norm{\varphi_{n}}_{L^{\infty}}^{r} \right ) \norm{\pd_{t}\varphi_{n}}_{L^{2}} \norm{\nabla \mu_{n}}_{\bm{L}^{2}}^{\frac{1}{2}} \norm{\nabla \Laplace \mu_{n}}_{\bm{L}^{2}}^{\frac{1}{2}} \\
&\leq \eps \norm{\nabla \Laplace \mu_{n}}_{\bm{L}^{2}}^{2} + C\norm{\pd_{t}\varphi_{n}}_{L^{2}}^{2} + C\norm{\nabla \mu_{n}}_{\bm{L}^{2}}^{2},
\end{align*}
\begin{align*}
R_{16} &\leq  \chi \norm{\pd_{t}\sigma_{n}}_{L^{2}}\norm{\Laplace \mu_{n}}_{L^{2}} \leq \eps \norm{\nabla \Laplace \mu_{n}}_{\bm{L}^{2}}^{2} + C\norm{\pd_{t}\sigma_{n}}_{L^{2}}^2+ C\norm{\nabla \mu_{n}}_{\bm{L}^{2}}^{2},
\end{align*}
\begin{align*}
R_{17} & \leq B \norm{\nabla \Gamma}_{\bm{L}^{2}} \norm{\nabla \Laplace \mu_{n}}_{\bm{L}^{2}} \leq \eps \norm{\nabla \Laplace \mu_{n}}_{\bm{L}^{2}}^{2} + C \norm{\nabla \Gamma}_{\bm{L}^{2}}^{2},
\end{align*}
where we used the fact that $\pdnu (\Laplace \mu_{n}) = 0$ on $\pd \Omega$ and integration by parts to deduce that
\begin{align*}
\norm{\Laplace \mu_{n}}_{L^{2}}^{2} = \int_{\Omega} \nabla \Laplace \mu_{n} \cdot \nabla \mu_{n} \dx \leq \norm{\nabla \Laplace \mu_{n}}_{\bm{L}^{2}} \norm{\nabla \mu_{n}}_{\bm{L}^{2}}.
\end{align*}
By Agmon's inequality $\norm{f}_{L^{\infty}} \leq C \norm{f}_{H^{2}}^{\frac{1}{2}} \norm{f}_{L^{2}}^{\frac{1}{2}}$, we see that
\begin{align*}
R_{18} & \leq B\norm{ m'(\varphi_{n})}_{L^{\infty}} \norm{\Laplace \mu_{n}}_{L^{2}} \norm{ \nabla \varphi_{n}}_{\bm{L}^{\infty}} \norm{ \nabla \Laplace \mu_{n} }_{\bm{L}^{2}} \\
&\leq C \norm{ \nabla \varphi_{n}}_{\bm{L}^{2}}^{\frac{1}{2}} \norm{\varphi_{n}}_{H^{3}}^{\frac{1}{2}} \norm{\nabla \mu_{n}}_{\bm{L}^{2}}^{\frac{1}{2}} \norm{ \nabla \Laplace \mu_{n} }_{\bm{L}^{2}}^{\frac{3}{2}} \\
&\leq \eps \norm{\nabla \Laplace \mu_{n}}_{\bm{L}^{2}}^{2} + C \norm{\varphi_{n}}_{H^{3}}^{2} \norm{\nabla  \mu_{n}}_{\bm{L}^{2}}^{2}.
\end{align*}
In a similar manner, expanding the gradient term in $R_{19}$ leads to
\begin{align*}
R_{19} & \leq  B \norm{m''(\varphi_{n})}_{L^{\infty}} \norm{\nabla \varphi_{n}}_{\bm{L}^{\infty}}^{2} \norm{\nabla \mu_{n}}_{\bm{L}^{2}} \norm{\nabla \Laplace \mu_{n}}_{\bm{L}^{2}} \\
&\quad + B \norm{m'(\varphi_{n})}_{L^{\infty}} \norm{\nabla^{2} \varphi_{n}}_{\bm{L}^{\infty}} \norm{\nabla \mu_{n}}_{\bm{L}^{2}} \norm{\nabla \Laplace \mu_{n}}_{\bm{L}^{2}} \\
&\quad + B \norm{m'(\varphi_{n})}_{L^{\infty}} \norm{\nabla \varphi_{n}}_{\bm{L}^{\infty}} \norm{ \nabla^{2} \mu_{n}}_{\bm{L}^{2}} \norm{\nabla \Laplace \mu_{n}}_{\bm{L}^{2}} \\
&\leq C \norm{\nabla \varphi_{n}}_{\bm{L}^{2}}\norm{\varphi_{n}}_{H^{3}} \norm{\nabla \mu_{n}}_{\bm{L}^{2}} \norm{\nabla \Laplace \mu_{n}}_{\bm{L}^{2}} \\
&\quad + C \norm{\varphi_{n}}_{H^{2}}^{\frac{1}{2}} \norm{\varphi_{n}}_{H^{4}}^{\frac{1}{2}} \norm{\nabla \mu_{n}}_{\bm{L}^{2}} \norm{\nabla \Laplace \mu_{n}}_{\bm{L}^{2}} \\
&\quad + C\norm{\nabla \varphi_{n}}_{\bm{L}^{2}}^{\frac{1}{2}} \norm{\varphi_{n}}_{H^{3}}^{\frac{1}{2}} \norm{\nabla \mu_{n}}_{\bm{L}^{2}}^{\frac{1}{2}}\norm{\nabla \mu_{n}}_{\bm{H}^{2}}^{\frac{1}{2}} \norm{\nabla \Laplace \mu_{n}}_{\bm{L}^{2}} \\
&\leq \eps \norm{\nabla \Laplace \mu_{n}}_{\bm{L}^{2}}^{2} + C \left (1+ \norm{\varphi_{n}}_{H^{4}}^{2} \right )\norm{\nabla \mu_{n}}_{\bm{L}^{2}}^{2}.
\end{align*}
Lastly, for $R_{20}$ we have
\begin{align*}
R_{20} & \leq  C \norm{\nabla \bm{v}_{n}}_{\bm{L}^{2}} \norm{\nabla \varphi_{n}}_{\bm{L}^{\infty}} \norm{\nabla \Laplace \mu_{n}}_{\bm{L}^{2}} + C \norm{ \bm{v}_{n}}_{\bm{L}^{4}} \norm{\nabla^{2} \varphi_{n}}_{\bm{L}^{4}} \norm{\nabla \Laplace \mu_{n}}_{\bm{L}^{2}} \\
&\leq C \norm{\nabla \bm{v}_{n}}_{\bm{L}^{2}} \norm{\varphi_{n}}_{H^{3}} \norm{\nabla \Laplace \mu_{n}}_{\bm{L}^{2}} \nonumber\\
&\quad + C \norm{\bm{v}_{n}}_{\bm{L}^{2}}^{\frac{1}{2}} \norm{\bm{v}_{n}}_{\bm{H}^{1}}^{\frac{1}{2}} \norm{\varphi_{n}}_{H^{2}}^{\frac{1}{2}} \norm{\varphi_{n}}_{H^{3}}^{\frac{1}{2}} \norm{\nabla \Laplace \mu_{n}}_{\bm{L}^{2}} \\
&\leq \eps  \norm{\nabla \Laplace \mu_{n}}_{\bm{L}^{2}}^{2} +C\norm{\bm{v}_{n}}_{\bm{H}^{1}}^{2} + C \norm{\varphi_{n}}_{H^{3}}^{2}.
\end{align*}
Collecting the above estimates together, we infer from \eqref{mu:t:mu} that
\begin{align*}
&{\frac{1}{2}} \frac{\dd}{\dt} \norm{\nabla \mu_{n}}_{\bm{L}^{2}}^{2} + \left(B m_{0} - 6 \eps \right)\norm{\nabla \Laplace \mu_{n}}_{\bm{L}^{2}}^{2}\\
&\quad \leq C \left (1+ \norm{\varphi_{n}}_{H^{4}}^{2} \right )\norm{\nabla \mu_{n}}_{\bm{L}^{2}}^{2} + C\norm{\pd_{t}\varphi_{n}}_{L^{2}}^{2} + C\norm{\pd_{t}\sigma_{n}}_{L^{2}}^{2} \\
&\qquad +C \norm{\bm{v}_{n}}_{\bm{H}^{1}}^{2} + C  \norm{ \varphi_{n}}_{H^{3}}^{2} + C\norm{\nabla \Gamma}_{\bm{L}^{2}}^{2}.
\end{align*}
Taking $\eps =\frac{Bm_{0}}{12}$, then it follows from \eqref{Apriori:3}, \eqref{Apriori:4} and Gronwall's lemma that
\begin{align*}
\sup_{t\in [0,T]} \norm{\nabla \mu_{n}(t)}_{\bm{L}^{2}}^{2} + \norm{\nabla \Laplace \mu_{n}}_{L^{2}(0,T; \bm{L}^{2})}^{2} \leq C.
\end{align*}
The above estimate combined with \eqref{mean:mu} and Poincar\'{e}'s inequality yields
\begin{align}\label{Apriori:5}
\norm{\mu_{n}}_{L^{\infty}(0,T; H^{1})\cap L^{2}(0,T; H^{3})}\leq C.
\end{align}
Going back to \eqref{Galerkin:varphi} we further deduce that
\begin{align}\label{Apriori:6}
\norm{\pd_{t}\varphi_{n}}_{L^{2}(0,T;H^{1})} \leq C£¬
\end{align}
and by \eqref{Apriori:3}, it follows
\begin{align}\label{Apriori:7}
\norm{\varphi_{n}}_{L^{\infty}(0,T;H^{3})} \leq C.
\end{align}

\paragraph{Passing to the limit as $n\to+\infty$.}
By the above higher order a priori estimates, there exists a quadruple of functions $(\varphi, \mu, \sigma, \bm{v})$ that satisfies the regularities stated in Theorem \ref{thm:2D:strong} and
\begin{subequations}
\begin{align}
0 & = \pd_{t} \varphi + \bm{v} \cdot \nabla \varphi - \div(m(\varphi) \mu) - \Gamma, \\
0 & = \mu - A \Psi'(\varphi) + B \Laplace \varphi + \chi \sigma, \\
0 & = \pd_{t} \sigma + \bm{v} \cdot \nabla \sigma - \div (n(\varphi)\nabla (\sigma + \chi (1-\varphi))) - S, \label{StrSoln:sigma} \\
\bm{0} & = \pd_{t} \bm{v} + (\bm{v} \cdot \nabla) \bm{v} - \div (2 \eta(\varphi) \der \bm{v}) - (\mu + \chi \sigma) \nabla \varphi, \label{StrSoln:velo}
\end{align}
\end{subequations}
hold a.e in $Q$.  It remains to derive estimates for the pressure.  From \eqref{StrSoln:velo}, there exists a function $q \in L^{2}(0,T;L^{2}_{0})$ such that
\begin{align*}
\pd_{t} \bm{v}(t) + (\bm{v}(t) \cdot \nabla) \bm{v}(t) - \div (2 \eta(\varphi(t)) \der \bm{v}(t)) - (\mu(t) + \chi \sigma(t)) \nabla \varphi(t) = - \nabla q(t)
\end{align*}
holds as an equality in the sense of distribution for a.e. $t \in (0,T)$.   Since the left-hand side belongs to $\bm{L}^{2}$ for a.e. $t \in (0,T)$, we obtain that $\nabla q \in L^{2}(0,T;\bm{L}^{2})$.
\begin{remark}[Further H\"{o}lder regularity for $\sigma$]
In the case $\chi = 0$ and $S = 0$, the equation \eqref{StrSoln:sigma} reduces to
\begin{align}\label{sigma:Holder:reg}
\pd_{t}\sigma + \bm{v} \cdot \nabla \sigma - \div (n(\varphi) \nabla \sigma)=0.
\end{align}
If in addition, $\sigma_{0} \in L^{\infty}(\Omega)$, then by a weak comparison principle, i.e., testing \eqref{sigma:Holder:reg} with $(\sigma - l)_{-} = \max (l - \sigma, 0)$, $l = \inf_{x \in \Omega} \sigma_{0}(x)$, and testing \eqref{sigma:Holder:reg} with $(\sigma - m)_{+} = \max(\sigma - m, 0)$, $m = \sup_{x \in \Omega} \sigma_{0}(x)$, it follows that $\sigma \in L^{\infty}(0,T;L^{\infty})$ and $\norm{\sigma}_{L^{\infty}(0,T; L^\infty)} \leq \norm{\sigma_{0}}_{L^{\infty}}$. For more details, see for instance \cite{GLRome,Lorca}. We just remark that the convection term can be handled as follows:
\begin{align*}
(\bm{v} \cdot \nabla (\sigma - l), (\sigma - l)_{-}) = -\frac{1}{2} (\bm{v} , \nabla \abs{(\sigma - l)_{-}}^{2}) = \frac{1}{2}(\div \bm{v}, \abs{(\sigma - l)_{-}}^{2}) = 0.
\end{align*}
Then the divergence-free property of $\bm{v}$ and the fact that $\bm{v} \in L^{4}(Q)$ allow us to deduce the existence of constants $C > 0$ and $\alpha \in (0,1)$, depending on $\norm{\sigma}_{L^{\infty}(0,T;L^{\infty})}$ and $\norm{\bm{v}}_{L^{4}(Q)}$ such that
\begin{align*}
\abs{\sigma(x,t) - \sigma(y,s)} \leq C \left (\abs{x-y}^{\alpha} + \abs{t-s}^\frac{\alpha}{2} \right), \quad \forall\, (x,t),(y,s) \in \overline{Q},\ x \neq y,\ s \neq t.
\end{align*}
The above estimate follows from similar arguments to the proof of \cite[Lem. 3.2]{SunZhang} (see also \cite[Lem. 2]{FGG}).
\end{remark}

\section{Proof of Theorem \ref{thm:2D:ctsdep}: continuous dependence in two dimensions}\label{sec:ctsdep2D}

Let $\{\varphi_{i}, \mu_{i}, \sigma_{i}, \bm{v}_{i}, q_{i}\}_{i=1,2}$ denote two global strong solutions to problem \eqref{ZEED}--\eqref{ZEED:bdy} corresponding to initial data $\{\varphi_{0,i}, \sigma_{0,i}, \bm{v}_{0,i}\}_{i=1,2}$ and source terms $\{\Gamma_{i}, S_{i}\}_{i=1,2}$.
Denoting the differences as $\hat{f}$ for $f \in \{\varphi, \mu, \sigma, \bm{v}, q, \varphi_{0}, \sigma_{0}, \bm{v}_{0}, \Gamma, S\}$, we obtain the system satisfied by the difference of solutions
\begin{subequations}
\begin{align}
0 & = \pd_{t}\hat{\varphi} + \hat{\bm{v}} \cdot \nabla \varphi_{1} + \bm{v}_{2} \cdot \nabla \hat{\varphi} - \div( (m(\varphi_{1}) - m(\varphi_{2})) \nabla \mu_{1} + m(\varphi_{2}) \nabla \hat{\mu}) - \hat{\Gamma}, \label{CtsDep:varphi} \\
0 & = \hat{\mu} - A (\Psi'(\varphi_{1}) - \Psi'(\varphi_{2})) + B \Laplace \hat{\varphi} + \chi \hat{\sigma}, \label{CtsDep:mu} \\
0 & = \pd_{t} \hat{\sigma} + \hat{\bm{v}} \cdot \nabla \sigma_{1} + \bm{v}_{2} \cdot \nabla \hat{\sigma}  - \hat{S} \label{CtsDep:sigma} \\
\notag & \quad - \div (n(\varphi_{1}) \nabla (\hat{\sigma} - \chi \hat{\varphi}) + (n(\varphi_{1}) - n(\varphi_{2})) \nabla (\sigma_{2} + \chi (1-\varphi_{2}))), \\
\mathbf{0} & = \pd_{t}\hat{\bm{v}} + (\hat{\bm{v}} \cdot \nabla) \bm{v}_{1} + (\bm{v}_{2} \cdot \nabla) \hat{\bm{v}} + \nabla \hat{q} - (\hat{\mu} + \chi \hat{\sigma}) \nabla \varphi_{1} - (\mu_{2} + \chi \sigma_{2}) \nabla \hat{\varphi} \label{CtsDep:velo} \\
\notag & \quad - 2 \div (  \eta(\varphi_{1}) \der \hat{\bm{v}} + (\eta(\varphi_{1})-\eta(\varphi_{2})) \der \bm{v}_{2}),
\end{align}
\end{subequations}
with
\begin{align*}
\hat{\varphi}(0) = \hat{\varphi}_{0}, \quad \hat{\sigma}(0) = \hat{\sigma}_{0}, \quad \hat{\bm{v}}(0) = \hat{\bm{v}}_{0}.
\end{align*}
Testing \eqref{CtsDep:varphi} with $B\hat{\varphi}$, \eqref{CtsDep:mu} with $m(\varphi_{2}) \hat{\mu}$, \eqref{CtsDep:sigma} with $\hat{\sigma}$, \eqref{CtsDep:velo} with $\hat{\bm{v}}$, and \eqref{CtsDep:mu} with $K \hat{\varphi}$ for some positive constant $K$ yet to be determined, then integrate by parts and summing leads to
\begin{align}
& \frac{1}{2} \frac{\dd}{\dt} \left ( B \norm{\hat{\varphi}}_{L^{2}}^{2} + \norm{\hat{\bm{v}}}_{\bm{L}^{2}}^{2} + \norm{\hat{\sigma}}_{L^{2}}^{2} \right ) \nonumber \\
& \qquad  + \norm{m^{\frac{1}{2}}(\varphi_{2})\hat{\mu}}_{L^{2}}^{2} + 2  \norm{ \eta^{\frac{1}{2}}(\varphi_{1}) \der \hat{\bm{v}}}_{\bm{L}^{2}}^{2} + \norm{n^{\frac{1}{2}}(\varphi_{1}) \nabla \hat{\sigma}}_{\bm{L}^{2}}^{2} + BK \norm{\nabla \hat{\varphi}}_{\bm{L}^{2}}^{2} \nonumber \\
& \quad = \int_{\Omega} - B \hat{\varphi} \hat{\bm{v}} \cdot \nabla \varphi_{1} + B \hat{\Gamma} \hat{\varphi} + m(\varphi_{2}) \left (A (\Psi(\varphi_{1})-\Psi(\varphi_{2})) \hat{\mu} - \chi \hat{\sigma} \hat{\mu} \right )  \dx \nonumber \\
& \qquad + \int_{\Omega} B m'(\varphi_{2}) \hat{\mu} \nabla \varphi_{2} \cdot \nabla \hat{\varphi} \dx +\int_\Omega (\hat{\mu} + \chi \hat{\sigma}) \nabla \varphi_{1} \cdot \hat{\bm{v}} + (\mu_{2} + \chi \sigma_{2}) \nabla \hat{\varphi} \cdot \hat{\bm{v}} \dx\nonumber\\
&\qquad + \int_\Omega 2(\eta(\varphi_{2})-\eta(\varphi_{1})) \der \bm{v}_{2} \cdot \der \hat{\bm{v}} - (\hat{\bm{v}} \cdot \nabla) \bm{v}_{1} \cdot \hat{\bm{v}} \dx\nonumber  \\
& \qquad + \int_{\Omega} \hat{S} \hat{\sigma} + (n(\varphi_{2})-n(\varphi_{1})) \nabla (\sigma_{2} - \chi \varphi_{2}) \cdot \nabla \hat{\sigma} \dx \nonumber \\
& \qquad - \int_{\Omega} \chi n(\varphi_{1}) \nabla \hat{\varphi} \cdot \nabla \hat{\sigma} + \hat{\sigma} \hat{\bm{v}} \cdot \nabla \sigma_{1} \dx -\int_\Omega B (m(\varphi_{1}) - m(\varphi_{2})) \nabla \mu_{1} \cdot \nabla \hat{\varphi} \dx \nonumber \\
& \qquad + K \int_{\Omega} (\hat{\mu} + \chi \hat{\sigma}) \hat{\varphi} - A (\Psi'(\varphi_{1})-\Psi'(\varphi_{2})) \hat{\varphi} \dx\nonumber\\
& \quad  =: \sum_{j=1}^8 I_j,\label{difference}
\end{align}
where we have used the following facts
\begin{align*}
&(\hat{\varphi},\bm{v}_{2} \cdot \nabla \hat{\varphi}) = \frac{1}{2}(\bm{v}_{2}, \nabla \abs{\hat{\varphi}}^{2}) = -\frac{1}{2} (\div \bm{v}_{2}, \abs{\hat{\varphi}}^{2}) = 0,\nonumber\\
&((\bm{v}_{2} \cdot \nabla) \hat{\bm{v}}, \hat{\bm{v}}) = - \frac{1}{2}( \div \bm{v}_{2}, \abs{\hat{\bm{v}}}^{2}) = 0,\\
&(\nabla \hat{q}, \hat{\bm{v}}) = -(\div \hat{\bm{v}}, q) = 0, \quad (\bm{v}_{2} \cdot \nabla \hat{\sigma}, \hat{\sigma}) = 0.
\end{align*}
We now proceed to estimate the right-hand side of \eqref{difference} term by term.
 First, assumption \eqref{assump:Psi:ctsdep} yields that
\begin{align*}
I_{1}+I_2
 & \leq B \norm{\hat{\varphi}}_{L^{2}} \norm{\hat{\bm{v}}}_{\bm{L}^{2}} \norm{\nabla \varphi_{1}}_{\bm{L}^{\infty}} + B \norm{\hat{\Gamma}}_{L^{2}} \norm{\hat{\varphi}}_{L^{2}} + m_{1} \chi \norm{\hat{\sigma}}_{L^{2}} \norm{\hat{\mu}}_{L^{2}} \\
& \quad + m_{1} A \left ( 1 + \norm{\varphi_{1}}_{L^{\infty}}^{r} + \norm{\varphi_{2}}_{L^{\infty}}^{r} \right ) \norm{\hat{\varphi}}_{L^{2}} \norm{\hat{\mu}}_{L^{2}} \nonumber\\
&\quad + B \norm{m'}_{L^{\infty}} \norm{\nabla \varphi_{2}}_{\bm{L}^{\infty}} \norm{\hat{\mu}}_{L^{2}} \norm{\nabla \hat{\varphi}}_{\bm{L}^{2}} \\
& \leq \delta \norm{\hat{\mu}}_{L^{2}}^{2} + C_\delta \norm{\hat{\sigma}}_{L^{2}}^{2} + C \norm{\hat{\Gamma}}_{L^{2}}^{2} + \norm{\hat{\bm{v}}}_{\bm{L}^{2}}^{2}
+ C \norm{\nabla \varphi_{2}}_{\bm{L}^{\infty}}^{2} \norm{\nabla \hat{\varphi}}_{\bm{L}^{2}}^{2}\\
& \quad + C \left ( 1 + \norm{\varphi_{1}}_{L^{\infty}}^{2r} +\norm{\varphi_{2}}_{L^{\infty}}^{2r}+ \norm{\nabla \varphi_{1}}_{\bm{L}^{\infty}}^{2} \right ) \norm{\hat{\varphi}}_{L^{2}}^{2},
\end{align*}
for some $\delta > 0$ yet to be determined.  Similarly,  we have
\begin{align*}
I_{8} & \leq K \left (\norm{\hat{\mu}}_{L^{2}} + \chi \norm{\hat{\sigma}}_{L^{2}} \right ) \norm{\hat{\varphi}}_{L^{2}} + K A \left ( 1 + \norm{\varphi_{1}}_{L^{\infty}}^{r} + \norm{\varphi_{2}}_{L^{\infty}}^{r} \right ) \norm{\hat{\varphi}}_{L^{2}}^{2} \\
& \leq \delta \norm{\hat{\mu}}_{L^{2}}^{2} + \norm{\hat{\sigma}}_{L^{2}}^{2} + \left (C K^{2} + CK \left (1 + \norm{\varphi_{1}}_{L^{\infty}}^{r}+\norm{\varphi_{2}}_{L^{\infty}}^{r} \right ) \right ) \norm{\hat{\varphi}}_{L^{2}}^{2}.
\end{align*}
Next, thanks to the Lipschitz continuity of $m$, we have
\begin{align*}
I_{7} \leq C \norm{\nabla \mu_{1}}_{\bm{L}^{\infty}} \norm{\hat{\varphi}}_{L^{2}} \norm{\nabla \hat{\varphi}}_{\bm{L}^{2}} \leq  \norm{\nabla \hat{\varphi}}_{\bm{L}^{2}}^{2} + C \norm{\mu_{1}}_{H^{3}}^{2} \norm{\hat{\varphi}}_{L^{2}}^{2}.
\end{align*}
Similarly, the Lipschitz continuity of $n$, Young's inequality and Korn's inequality allow us to deduce that
\begin{align*}
I_5+I_6 & \leq \norm{\hat{S}}_{L^{2}} \norm{\hat{\sigma}}_{L^{2}}
+ C\norm{\nabla (\sigma_{2} - \chi \varphi_{2})}_{\bm{L}^{4}} \norm{\hat{\varphi}}_{L^{4}} \norm{\nabla \hat{\sigma}}_{\bm{L}^{2}} \\
& \quad + C \norm{\nabla \hat{\varphi}}_{\bm{L}^{2}}\norm{\nabla \hat{\sigma}}_{\bm{L}^{2}}
+ \norm{\hat{\bm{v}}}_{\bm{L}^{4}} \norm{\hat{\sigma}}_{L^{2}} \norm{\nabla \sigma_{1}}_{\bm{L}^{4}} \\
& \leq \norm{\hat{\sigma}}_{L^{2}}^{2} + C \norm{\hat{S}}_{L^{2}}^{2}
+ \eps \norm{\nabla \hat{\sigma}}_{\bm{L}^{2}}^{2} + C_\eps \left (\norm{\nabla (\sigma_{2} - \chi \varphi_{2})}_{\bm{L}^{4}}^{2} \norm{\hat{\varphi}}_{L^{4}}^{2} + \norm{\nabla \hat{\varphi}}_{\bm{L}^{2}}^{2} \right ) \\
& \quad + \theta \norm{\der \hat{\bm{v}}}_{\bm{L}^{2}}^{2}
+ C_\theta \norm{\hat{\bm{v}}}_{\bm{L}^{2}}^{2}
+ C_\theta \norm{\nabla \sigma_{1}}_{\bm{L}^{4}}^{2} \norm{\hat{\sigma}}_{L^{2}}^2,
\end{align*}
for some positive constants $\eps, \theta$ yet to be determined.  By the Gagliardo--Nirenberg inequality $\norm{f}_{L^{4}} \leq C \norm{f}_{L^{2}}^{\frac{1}{2}} \norm{\nabla f}_{L^{2}}^{\frac{1}{2}} + C \norm{f}_{L^{2}}$ in two dimensions, we have
\begin{align*}
\norm{\nabla (\sigma_{2} - \chi \varphi_{2})}_{\bm{L}^{4}}^{2} \norm{\hat{\varphi}}_{L^{4}}^{2} & \leq C \norm{\nabla (\sigma_{2} - \chi \varphi_{2})}_{\bm{L}^{4}}^{2} \left ( \norm{\hat{\varphi}}_{L^{2}} \norm{\nabla \hat{\varphi}}_{L^{2}} + \norm{\hat{\varphi}}_{L^{2}}^{2} \right ) \\
& \leq \norm{\nabla \hat{\varphi}}_{L^{2}}^{2} + C \left ( 1 + \norm{\nabla (\sigma_{2} - \chi \varphi_{2})}_{\bm{L}^{4}}^{4} \right ) \norm{\hat{\varphi}}_{L^{2}}^{2}.
\end{align*}
Then it follows  that
\begin{align*}
I_5+I_6 & \leq  \eps \norm{\nabla \hat{\sigma}}_{L^{2}}^{2} + C_{\eps} \norm{\nabla \hat{\varphi}}_{\bm{L}^{2}}^{2} +  \theta \norm{\der \hat{\bm{v}}}_{\bm{L}^{2}}^{2}
+ \left (1 + C_\theta \norm{\nabla \sigma_{1}}_{\bm{L}^{4}}^{2} \right ) \norm{\hat{\sigma}}_{L^{2}}^{2}  \\
& \quad + C_{\eps} \left (1 + \norm{\nabla (\sigma_{2} - \chi \varphi_{2})}_{\bm{L}^{4}}^{4} \right ) \norm{\hat{\varphi}}_{L^{2}}^{2}
+ C_\theta \norm{\hat{\bm{v}}}_{\bm{L}^{2}}^{2} + C \norm{\hat{S}}_{L^{2}}^{2}.
\end{align*}
Lastly, thanks to the Lipschitz continuity of $\eta$, it holds that
\begin{align*}
I_3+I_4
 & \leq \norm{\nabla \varphi_{1}}_{\bm{L}^{\infty}} (\norm{\hat{\mu}}_{L^{2}}
        + \chi \norm{\hat{\sigma}}_{L^{2}}) \norm{\hat{\bm{v}}}_{\bm{L}^{2}}
        + \norm{\mu_{2} + \chi \sigma_{2}}_{L^{\infty}} \norm{\hat{\bm{v}}}_{\bm{L}^{2}} \norm{\nabla \hat{\varphi}}_{\bm{L}^{2}} \\
& \quad + C \norm{\der \bm{v}_{2}}_{\bm{L}^{4}} \norm{\der \hat{\bm{v}}}_{\bm{L}^{2}} \norm{\hat{\varphi}}_{L^{4}}
        + \norm{\nabla \bm{v}_{1}}_{\bm{L}^{2}} \norm{\hat{\bm{v}}}_{\bm{L}^{4}}^{2}  \\
& \leq \delta \norm{\hat{\mu}}_{L^{2}}^{2} + \norm{\hat{\sigma}}_{L^{2}}^{2}
        +  \norm{\nabla \hat{\varphi}}_{\bm{L}^{2}}^{2}
        + \left (C_{\delta} \norm{\nabla \varphi_{1}}_{\bm{L}^{\infty}}^{2}
        + C \norm{\mu_{2} + \chi \sigma_{2}}_{L^{\infty}}^{2} \right ) \norm{\hat{\bm{v}}}_{\bm{L}^{2}}^{2} \\
& \quad + \theta \norm{\der \hat{\bm{v}}}_{\bm{L}^{2}}^{2}
        + C_{\theta} \norm{\der \bm{v}_{2}}_{\bm{L}^{4}}^{2} \norm{\hat{\varphi}}_{L^{4}}^{2}
        + \norm{\nabla \bm{v}_{1}}_{\bm{L}^{2}} \norm{\hat{\bm{v}}}_{\bm{L}^{4}}^{2}.
\end{align*}
By the Gagliardo--Nirenberg inequality $\norm{f}_{L^{4}} \leq C \norm{f}_{L^{2}}^{\frac{1}{2}} \norm{\nabla f}_{L^{2}}^{\frac{1}{2}} + C \norm{f}_{L^{2}}$ in two dimensions, and Korn's inequality, we have
\begin{align*}
\norm{\nabla \bm{v}_{1}}_{\bm{L}^{2}} \norm{\hat{\bm{v}}}_{\bm{L}^{4}}^{2}
& \leq C \norm{\nabla \bm{v}_{1}}_{\bm{L}^{2}} \left ( \norm{\hat{\bm{v}}}_{\bm{L}^{2}} \norm{\der \hat{\bm{v}}}_{\bm{L}^{2}}
       + \norm{\hat{\bm{v}}}_{\bm{L}^{2}}^{2} \right ) \\
& \leq \theta \norm{\der \hat{\bm{v}}}_{\bm{L}^{2}}^{2} + (C_{\theta} \norm{\nabla \bm{v}_{1}}_{\bm{L}^{2}}^{2}
       + C\norm{\nabla \bm{v}_{1}}_{\bm{L}^{2}})\norm{\hat{\bm{v}}}_{\bm{L}^{2}}^2,
\end{align*}
\begin{align*}
C_{\theta} \norm{\der \bm{v}_{2}}_{\bm{L}^{4}}^{2} \norm{\hat{\varphi}}_{L^{4}}^{2} & \leq C_{\theta} \norm{\der \bm{v}_{2}}_{\bm{L}^{4}}^{2} (\norm{\hat{\varphi}}_{L^{2}} \norm{\nabla \hat{\varphi}}_{\bm{L}^{2}} + \norm{\hat{\varphi}}_{L^{2}}^{2} )\\
& \leq \norm{\nabla \hat{\varphi}}_{\bm{L}^{2}}^{2} + C_{\theta} (1 + \norm{\der \bm{v}_{2}}_{\bm{L}^{4}}^{4} )\norm{\hat{\varphi}}_{L^{2}}^{2},
\end{align*}
As a consequence,
\begin{align*}
I_3+I_4
 & \leq 2 \theta \norm{\der \hat{\bm{v}}}_{\bm{L}^{2}}^{2}
 + \delta \norm{\hat{\mu}}_{L^{2}}^{2}
 + 2 \norm{\nabla \hat{\varphi}}_{\bm{L}^{2}}^{2} + \norm{\hat{\sigma}}_{L^{2}}^{2}
 + C_{\theta} \left ( 1 + \norm{\der \bm{v}_{2}}_{\bm{L}^{4}}^{4} \right ) \norm{\hat{\varphi}}_{\bm{L}^{2}}^{2}  \\
& \quad + \left (C_{\delta} \norm{\nabla \varphi_{1}}_{\bm{L}^{\infty}}^{2}
        + C\norm{\mu_{2} + \chi \sigma_{2}}_{L^{\infty}}^{2}
        + C_{\theta} \norm{\nabla \bm{v}_{1}}_{\bm{L}^{2}}^{2}
        + C \right )\norm{\hat{\bm{v}}}_{\bm{L}^{2}}^{2}.
\end{align*}
In summary, we infer from \eqref{difference} and the above estimates that
\begin{equation}\label{2D:CtsDep:Est}
\begin{aligned}
& \frac{1}{2} \frac{\dd}{\dt} \left ( B \norm{\hat{\varphi}}_{L^{2}}^{2} + \norm{\hat{\bm{v}}}_{\bm{L}^{2}} + \norm{\hat{\sigma}}_{L^{2}}^{2} \right )   \\
& \quad + \left (m_{0} -3 \delta \right )\norm{\hat{\mu}}_{L^{2}}^{2}
+ (n_{0} - \eps) \norm{\nabla \hat{\sigma}}_{\bm{L}^{2}}^{2}
+ \left (2\eta_{0} - 3\theta \right ) \norm{\der \hat{\bm{v}}}_{\bm{L}^{2}}^{2} \\
&\quad + \left (BK - 3 - C_{\eps} - C_{\delta} \norm{\nabla \varphi_{2}}_{L^{\infty}(0,T;\bm{L}^{\infty})}^{2} \right ) \norm{\nabla \hat{\varphi}}_{\bm{L}^{2}}^{2} \\
& \, \leq C \big(\norm{\hat{\Gamma}}_{L^{2}}^{2}
+  \norm{\hat{S}}_{L^{2}}^{2}\big) + h_{0}(t) \norm{\hat{\sigma}}_{L^{2}}^{2}
+ h_{1}(t) \norm{\hat{\bm{v}}}_{\bm{L}^{2}}^{2}
+ h_{2}(t) \norm{\hat{\varphi}}_{L^{2}}^{2}
\end{aligned}
\end{equation}
where
\begin{align*}
h_{0}(t) & := C \left (1 + \norm{\nabla \sigma_{1}}_{\bm{L}^{4}}^{2} \right ), \\
h_{1}(t) & := C \left (1 + \norm{\nabla \varphi_{1}}_{\bm{L}^{\infty}}^{2} + \norm{\mu_{2} + \chi \sigma_{2}}_{L^{\infty}}^{2} + \norm{\nabla \bm{v}_{1}}_{\bm{L}^{2}}^{2} \right  ), \\
h_{2}(t) & :=  C \left ( 1 + \norm{\varphi_{1}}_{L^{\infty}}^{2r} + \norm{\varphi_{2}}_{L^{\infty}}^{2r} + \norm{\nabla \varphi_{1}}_{\bm{L}^{\infty}}^{2}\right)\nonumber\\
 &\qquad + C\left( \norm{\nabla (\sigma_{2} - \chi \varphi_{2})}_{\bm{L}^{4}}^{4} + \norm{\der \bm{v}_{2}}_{\bm{L}^{4}}^{4} + \norm{\mu_{1}}_{H^{3}}^{2} \right ),
\end{align*}
where in the expression of $h_0$, $h_1$ and $h_2$, the constant $C$ may depend on the parameters $\delta$, $\theta$, $\eps$ and $K$.

Since $\norm{\varphi_{2}}_{L^{\infty}(0,T;H^3)}$ is bounded (cf. \eqref{Apriori:7}), so is $\norm{\nabla \varphi_{2}}_{L^{\infty}(0,T;\bm{L}^{\infty})}$ by the Sobolev embedding theorem. Then in \eqref{2D:CtsDep:Est}, we can choose
 $$ \delta =\frac{ m_{0}}{6},\quad \eps =\frac{ n_{0}}{2},\quad  \theta = \frac{ \eta_{0}}{3},\quad K=  \frac{1}{B} \left (4 + C_{\eps} + C_{\delta} \norm{\nabla \varphi_{2}}_{L^{\infty}(0,T;\bm{L}^{\infty})}^{2} \right ).$$
 On the other hand, by regularities of the global strong solutions $\{ \varphi_{i}, \mu_{i}, \sigma_{i}, \bm{v}_{i}\}_{i=1,2}$ stated in Theorem \ref{thm:2D:strong} and the continuous embedding $L^{\infty}(0,T;L^{2}) \cap L^{2}(0,T;H^{1}) \subset L^{4}(Q)$ in two dimensions, it holds that $h_{0}, h_{1}, h_{2} \in L^{1}(0,T)$.  Then, applying Gronwall's lemma to \eqref{2D:CtsDep:Est}, and then using Korn's inequality, we have
\begin{align*}
& \sup_{t \in (0,T]} \left ( \norm{\hat{\varphi}(t)}_{L^{2}}^{2} + \norm{\hat{\bm{v}}(t)}_{\bm{L}^{2}}^{2} + \norm{\hat{\sigma}(t)}_{L^{2}}^{2} \right )\nonumber\\
&\qquad + \int_{0}^{T} \norm{\hat{\mu}}_{L^{2}}^{2} + \norm{\nabla \hat{\sigma}}_{\bm{L}^{2}}^{2} + \norm{\nabla \hat{\bm{v}}}_{\bm{L}^{2}}^{2} + \norm{\nabla \hat{\varphi}}_{\bm{L}^{2}}^{2} \dt \\
& \quad \leq C \left ( \norm{\hat{\Gamma}}_{L^{2}(Q)}^{2} + \norm{\hat{S}}_{L^{2}(Q)}^{2} + \norm{\hat{\sigma}_{0}}_{L^{2}}^{2} + \norm{\hat{\bm{v}}_{0}}_{\bm{L}^{2}}^{2} + \norm{\hat{\varphi}_{0}}_{L^{2}}^{2} \right )=: C\mathcal{Y}.
\end{align*}
To obtain continuous dependence of $\hat{\varphi}$ in the $L^{2}(0,T;H^{2})$-norm, first we see that by \eqref{assump:Psi:ctsdep}
\begin{align*}
\norm{\Psi'(\varphi_{1}) - \Psi'(\varphi_{2})}_{L^{2}(Q)}^{2} \leq C \left (  1 + \norm{\varphi_{1}}_{L^{\infty}(Q)}^{2r} + \norm{\varphi_{2}}_{L^{\infty}(Q)}^{2r} \right ) \norm{\hat{\varphi}}_{L^{2}(Q)}^{2} \leq C \mathcal{Y},
\end{align*}
then viewing \eqref{CtsDep:mu} as an elliptic problem for $\hat{\varphi}$ and applying the elliptic regularity theory leads to
\begin{align*}
\norm{\hat{\varphi}}_{L^{2}(0,T;H^{2})} \leq C \left (\norm{\hat{\mu}}_{L^{2}(Q)} + \norm{\Psi'(\varphi_{1})-\Psi'(\varphi_{2})}_{L^{2}(Q)} + \norm{\hat{\sigma}}_{L^{2}(Q)} \right ) \leq C\mathcal{Y}.
\end{align*}
The proof of Theorem \ref{thm:2D:ctsdep} is complete.

\section{Conclusion}
In this work, we derived a class of thermodynamically consistent Navier--Stokes--Cahn--Hilliard system for two-phase fluid flows with density contrast, based on a volume-averaged velocity in the spirit of \cite{AGG}, which also allows for mass transfer between the fluids and chemotactic response to a chemical species present in the physical domain. Besides, we present various simplified versions of the general model and state its sharp interface limit. As far as the mathematical analysis is concerned, for the general model some difficulties will be encountered in deriving a priori estimates
due to the mass transfer terms, especially when the volume-averaged velocity is not divergence-free.  Hence, we considered the simplest model variant as a starting point, which has a divergence-free velocity and no density contrast, but allows the fluid viscosity and mobilities to be dependent on the order parameter $\varphi$.  By a suitable Galerkin approximation, we establish global weak existence in both two and three dimensions for prescribed mass transfer terms, and then under additional assumptions we show global strong well-posedness as well as continuous dependence in two dimensions.  We believe that the a priori estimates derived here will also be useful to study long-time behavior of the system.  Furthermore, we expect that global weak existence to the more general model variant \eqref{mom}, \eqref{mu}, \eqref{sigma} and \eqref{Scaled:ZeroExcess} that has a solenodial velocity but allows for density contrast can be shown by employing the approach of \cite{ADG1}.

\section*{Acknowledgments}
K.F. Lam expresses his gratitude to Fudan University for its hospitality during his visit in which part of this research was completed.  H. Wu is partially supported by NNSFC under the grant No. 11631011 and Shanghai Center for Mathematical Sciences of Fudan University.

\bibliographystyle{plain}
\bibliography{NSCHTumor}
\end{document}